\documentclass[a4paper,11pt]{amsart}
\usepackage{mathrsfs}\usepackage{bbm}
\usepackage{latexsym, amssymb,amsmath,amsthm}

\def \dim{\operatorname{dim}}

\def \gr{\operatorname{gr}}
\def \Gr{\operatorname{Gr}}
\def \Hom{\operatorname{Hom}}
\def \Ext{\operatorname{Ext}}
\def \T{\operatorname{T}}

\def \C{\textrm{C}}

\def \End{\operatorname{End}}
\def \H{\operatorname{H}}
\def \Ker{\operatorname{Ker}}
\def \Im{\operatorname{Im}}
\def \deg{\operatorname{deg}}

\numberwithin{equation}{section}

\newtheorem{theorem}{Theorem}[section]
\newtheorem{lemma}[theorem]{Lemma}
\newtheorem{proposition}[theorem]{Proposition}
\newtheorem{corollary}[theorem]{Corollary}
\newtheorem{definition}[theorem]{Definition}

\newtheorem{remark}[theorem]{Remark}
\newtheorem{conjecture}[theorem]{Conjecture}

\begin{document}
\title[Support varieties of Lie superalgebras]{Support varieties and representations of tame basic classical Lie superalgebras }

\author{Gongxiang Liu}
\address{Department of Mathematics, Nanjing University, Nanjing 210093, China}
\email{gxliu@nju.edu.cn}
\date{}
\maketitle
\begin{abstract}
Let $\kappa$ be an algebraically closed field of characteristic
$p>3$ and $\mathfrak{g}$ a restricted Lie superalgebra over
$\kappa$. We introduce the definition of restricted cohomology for
$\mathfrak{g}$ and show its cohomology ring is finitely generated
provided $\mathfrak{g}$ is a basic classical Lie superalgebra. As a
consequence, we show that the restricted enveloping algebra of a
basic classical Lie superalgebra $\mathfrak{g}$ is always wild
except $\mathfrak{g}=\mathfrak{sl}_{2}$ or
$\mathfrak{g}=\mathfrak{osp}(1|2)$ or $\mathfrak{g}=\mathbf{C}(2)$.
 All finite dimensional indecomposable restricted representations
of $\mathbf{u}(\mathfrak{osp}(1|2))$, the restricted enveloping
algebra of Lie superalgebra $\mathfrak{osp}(1|2)$, are determined.
\vskip 5pt

\noindent{\bf 2000 MR Subject Classification} \ \ Primary 17B50,
17B56, Secondary 16G60, 16W30
\end{abstract}

\section{introduction}
\subsection{} As  generalizations and deep continuations of classical
Lie theory, Lie superalgebras, supergroups and their representation
theory over the field of complex numbers $\mathbb{C}$ have been
studied extensively since the classification of finite dimensional
complex simple Lie superalgebras by Kac \cite{Kac1}. More on
supergroups, supergeometry and supersymmetric theory can be found in
\cite{DM,Man}. In recent years, there has been increasing interest
in modular representation theory of algebraic supergroups.
Especially, the modular representations of $GL(m|n)$, $Q(n)$ and
ortho-symplectic supergroups have been initiated by Brundan,
Kleshchev, Kujawa \cite{B,BK1,BK2,BKu,Ku}, and Shu-Wang \cite{SW}. A
systemic research of modular Lie superalgebras has been started
\cite{WZ,WZ2}. In \cite{WZ}, the super version of the celebrated
Kac-Weisfeiler Property is shown to be hold for the basic classical
Lie superalgebras, which by definition admit an even nondegenerate
supersymmetric bilinear form and whose even subalgebras are
reductive. Actually, the modular representation theory of
supergroups and Lie superalgebras not only is of intrinsic interest
in its own right, but also has found remarkable applications to
classical mathematics. See \cite{SW} for some historical remarks.

Support varieties were introduced in the pioneering work of Alperin
\cite{Alp} and Carlson \cite{Ca1,Ca2} nearly 30 years ago as a
method to study complexes and resolutions of modules over group
algebras. They open an algebro-geometric gate to linear
representations of finite groups. Since then such ideas have been
extended to restricted Lie algebras \cite{FB}, Steenrod algebra
\cite{NP}, infinitesimal group schemes \cite{SFB}, arbitrary finite
dimensional cocommutative Hopf algebras \cite{FS} and even to finite
dimensional algebras \cite{SS}. See \cite{So} for a nice survey on
the theory of support varieties.

\subsection{} Up to now, we are lack of this algebro-geometric tool
for modular Lie superalgebras, perhaps due to the representation
theory of simple Lie superalgebras over $\mathbb{C}$ is already very
difficult and remains to be better understood. Recently, such tools
were introduced for Lie superalgebras over $\mathbb{C}$ in
\cite{BKN} by using so-called relative cohomology. It seems that the
methods used in \cite{BKN} can not be applied to positive
characteristic case. The main aim of this paper is to establish a
kind of definition for a support variety, which is suitable for our
purpose, and give an application. At first, we realize that for any
restricted Lie superalgebra $\mathfrak{g}$ one can relate it with an
ordinary Hopf algebra $\mathbf{u}(\mathfrak{g})\rtimes
\kappa\mathbb{Z}_{2}$ possessing equivalent representation theory as
$\mathbf{u}(\mathfrak{g})$. So we can pass from ``super world" to
the ``usual world" without losing information. Using this ordinary
Hopf algebra, we can define its cohomology algebra naturally.

It is known that support varieties can be defined once the finite
generation of cohomology is established, which is hard to prove in
general. In this paper, we prove this finite generation property for
the class of basic classical Lie superalgebras. It consists of
several infinite series and 3 exceptional ones. We divide our proof
in two different case $\mathfrak{g}\neq \mathbf{A}(1,1)$ or
$\mathfrak{g}=\mathbf{A}(1,1)$. In the first case, we give a
two-step filtration to reduce $\mathbf{u}(\mathfrak{g})$ to a
familiar algebra whose cohomology ring is known and each of
filtration involves a convergent spectral sequence.  We find some
permanent cycles in such spectral sequences and apply a lemma cited
from \cite{MPSW} to conclude finite generation. To give the
filtration, a new kind of PBW basis are developed. We put the case
$\mathfrak{g}=\mathbf{A}(1,1)$ in a  bigger context, in which all
$\mathbf{u}(\mathfrak{g})$ are equipped with a nice filtration
similar to the coradical filtration of a coalgebra. Through this
one-step filtration, we can reduce $\mathbf{u}(\mathbf{A}(1,1))$ to
a familiar algebra already. Then the same idea developed in the
first case can be applied.

One central question in the modern representation theory of algebras
is the determination of the representation type. By Drozd's
fundamental trichotomy \cite{Dr}, finite dimensional algebras over
an algebraically closed field may be subdivided into the disjoint
classes of representation finite, tame and wild algebras. As an
application of support varieties we built, we will prove all
$\mathbf{u}(\mathfrak{g})$ are wild with only three exceptions:
$\mathfrak{g}=\mathfrak{sl}_{2}, \mathfrak{osp}(1|2),
\mathbf{C}(2)$. The case $\mathbf{C}(2)$ is conjectured to be wild
and we have known $\mathbf{u}(\mathfrak{sl}_{2})$ and
$\mathbf{u}(\mathfrak{osp}(1|2))$ are tame. Inspired by the
similarity between $\mathfrak{sl}_{2}$ and $\mathfrak{osp}(1|2)$ and
for further understanding of the representations of modular Lie
superalgebras, all finite dimensional  restricted  indecomposable
$\mathfrak{osp}(1|2)$-modules are also characterized.

The paper is organized as follows. All subsidiary results to prove
the finite generation of cohomology algebras are builded in Section
2. Especially, a new kind of PBW basis suitable for our purpose and
some filtrations are given. Section 3 is to give the proof of finite
generation. The definition of a support variety is given in Section
4. Moreover, its connections with complexity and representation type
are established. As the final conclusion of this section, the
representation type of any $\mathbf{u}(\mathfrak{g})$ is determined
except the case $\mathbf{C}(2)$, which is conjectured to be a  wild
algebra. In the last section of this paper, a complete list of  all
finite dimensional  restricted  indecomposable
$\mathfrak{osp}(1|2)$-modules up to isomorphism is formulated.

\section{preliminaries} Throughout of this paper, $\kappa$ is an
algebraically closed field of characteristic $p\neq 0$ and $p>3$ is
always assumed unless stated otherwise. All spaces are
$\kappa$-spaces. All modules are left modules.

\subsection{Hopf algebras in Yetter-Drinfeld categories} Let $J$ be a Hopf algebra
with bijective antipode and $^{J}_{J}\mathscr{YD}$ the category of
the Yetter-Drinfeld modules with left $J$-module action and left
$J$-comodule coaction. It is naturally forms a braided monoidal
category with the braiding
$$c_{M,N}:\;M\otimes N\to N\otimes M,\;\;m\otimes n\mapsto \sum n_{0}\otimes n_{-1}\cdot m,$$
where $n\mapsto \sum n_{-1}\otimes n_{0},\;N\to J\otimes N$ denotes
the comodule structure, as usual. Let $A$ be a \emph{braided Hopf
algebra} in $^{J}_{J}\mathscr{YD}$.  By definition, it is an algebra
as well as coalgebra in $^{J}_{J}\mathscr{YD}$ such that its
comultiplication and counit are algebra morphism, and such that the
identity morphism has a convolution inverse in
$^{J}_{J}\mathscr{YD}$. When we say that the comultiplication
$\Delta:\; A\to A\otimes A$ should be an algebra morphism, the
braiding defined as above arises in the definition of the algebra
structure of $A\otimes A$ and so $A$ is not an ordinary Hopf algebra
in general. Through the Radford-Majid bosonization
\cite{Majid,Radford}, it gives rise to an ordinary Hopf algebra
$A\rtimes J$. As an algebra, this is the smash product $A\# J$, and
it is the smash coproduct as a coalgebra.

\begin{lemma} Let $J$ be a Hopf algebra
with bijective antipode and $A$ a braided Hopf algebra in
$^{J}_{J}\mathscr{YD}$. Then the cohomology ring
$\H^{\ast}(A,\kappa):=\bigoplus_{i\geq
0}\Ext^{i}_{A}(\kappa,\kappa)$ is a braided graded commutative
algebra in $^{J}_{J}\mathscr{YD}$.
\end{lemma}
\begin{proof} By Theorem 3.12 in \cite{MPSW}, the Hochschild
cohomology ring $$\textrm{HH}^{\ast}(A,\kappa):=\bigoplus_{i\geq
0}\Ext^{i}_{A\otimes A^{op}}(A,\kappa)$$ is a braided graded
commutative algebra in $^{J}_{J}\mathscr{YD}$. By the standard bar
resolution for computing these extension groups, one can see that
$\Ext^{i}_{A}(\kappa,\kappa)\cong \Ext^{i}_{A\otimes
A^{op}}(A,\kappa)$ for $i\geq 0$ (see also subsection 2.4 in
\cite{MPSW}). The proof is complete.
\end{proof}

\subsection{Cohomology of restricted Lie superalgebras}

We fix some notions at first. By definition, a superalgebra is
nothing but a $\mathbb{Z}_{2}$-graded algebra. By forgetting the
grading we may consider any superalgebra  $A$ as a usual algebra and
this algebra will be denoted by $|A|$. For any two
$\mathbb{Z}_{2}$-graded vector spaces $V,W$, we use
$\Hom_{\kappa}(V,W)$ to represent the set of all linear maps from
$V$ to $W$ and $\underline{\Hom}_{\kappa}(V,W)$ to denote that of
all even linear maps.

Now let $A=A_{0}\oplus A_{1}$ be a superalgebra. Then there is a
natural action of $\mathbb{Z}_{2}=\langle g|g^{2}=1\rangle$ on $A$
given by
$$g\cdot a=a,\;\;g\cdot b=-b,\;\;\;\textrm{for}\; a\in A_{0}, \;b\in A_{1}.$$
Note that this definition makes sense as stated only for homogeneous
elements, it should be interpreted via linearity in the general
case. Thus $A$ is a $\kappa\mathbb{Z}_{2}$-module algebra (for
definition, see Section 4.1 in \cite{Mon}) and the smash product
$A\# \kappa\mathbb{Z}_{2}$ is a usual algebra. We use $A$-smod to
denote the category of all finitely generated left $A$-supermodules
with even homomorphisms and $A\# \kappa\mathbb{Z}_{2}$-mod the usual
finitely generated left $A\# \kappa\mathbb{Z}_{2}$-modules category.

\begin{lemma} Let $A$ be a superalgebra. Then $A$-\emph{smod} is equivalent
to $A\# \kappa\mathbb{Z}_{2}$-\emph{mod}.
\end{lemma}
\begin{proof} Let $M=M_{0}\oplus M_{1}$ be an $A$-supermodule and $g\in
\mathbb{Z}_{2}$ the generator of $\mathbb{Z}_{2}$. Through assigning
$g\cdot m_{0}:=m_{0},\;g\cdot m_{1}:=-m_{1}$ for $m_{0}\in
M_{0},\;m_{1}\in M_{1}$, $M$ is a $\kappa\mathbb{Z}_{2}$-module. Now
just define the action of $A\# \kappa\mathbb{Z}_{2}$ on $M$ through
$(a\otimes g)\cdot m:=a\cdot(g\cdot m)$ for $a\in A$ and $m\in M$.
To show it is indeed an $A\# \kappa\mathbb{Z}_{2}$-module, one need
verify the equality
$$(1\otimes g)(a\otimes 1)\cdot m=((g\cdot a)\otimes g)\cdot m,\;\;(\ast)$$
for $a\in A$ and $m\in M$. It is not hard to see that this is
equivalent to the fact $A_{i}M_{j}\subseteq M_{i+j}$ for $i,j\in
\mathbb{Z}_{2}$.

Conversely, let $M$ be an $A\# \kappa\mathbb{Z}_{2}$-module. Since
the characteristic of $\kappa$ is not equal to $2$,
$\kappa\mathbb{Z}_{2}$ is semisimple. Therefore, $M=M_{0}\oplus
M_{1}$ with $M_{0}=\{m\in M|g\cdot m=m\}$ and $M_{1}=\{m\in M|g\cdot
m=-m\}$.  Also, the equality $(\ast)$ implies that
$A_{i}M_{j}\subseteq M_{i+j}$ for $i,j\in \mathbb{Z}_{2}$. Thus $M$
is an $A$-supermodule.

At last, it is clear that $\underline{\Hom}_{A}(-,-)=\Hom_{A\#
\kappa\mathbb{Z}_{2}}(-,-)$. The lemma is proved.
\end{proof}

Now we specialize this simple observation to the case of restricted
enveloping algebras of restricted Lie superalgebras.

\begin{definition}\emph{ A Lie superalgebra $\mathfrak{g}=\mathfrak{g}_{0}\oplus
\mathfrak{g}_{1}$ is called a \emph{restricted Lie superalgebra}, if
there is a $p$th map $\mathfrak{g}_{0}\to \mathfrak{g}_{0}$, denoted
as $^{[p]}$, satisfying}

\emph{(a) } \emph{$(cx)^{[p]}=c^{p}x^{[p]}$ for all $c\in k$ and
$x\in \mathfrak{g}_{0}$,}

\emph{(b) }\emph{$[x^{[p]},y]=(adx)^{p}(y)$ for all $x\in
\mathfrak{g}_{0}$ and $y\in \mathfrak{g}$,}

\emph{(c)
}\emph{$(x+y)^{[p]}=x^{[p]}+y^{[p]}+\sum_{i=1}^{p-1}s_{i}(x,y)$ for
all $x,y\in \mathfrak{g}_{0}$ where $is_{i}$ is the coefficient of
$\lambda^{i-1}$ in $(ad(\lambda x+y))^{p-1}(x)$.}
\end{definition}

In short, a restricted Lie superalgebra is a Lie superalgebra whose
even subalgebra is a restricted Lie algebra and the odd part is a
restricted module by the adjoint action of the even subalgebra. All
the Lie (super)algebras in this paper will be assumed to be
restricted. For a restricted Lie superalgebra $\mathfrak{g}$,
$U(\mathfrak{g})$ is denoted to be its universal enveloping algebra
and $\mathbf{u}(\mathfrak{g})=U(\mathfrak{g})/(x^{p}-x^{[p]}|x\in
\mathfrak{g}_{0})$ its restricted enveloping algebra. The following
is a consequence of PBW theorem for $U(\mathfrak{g})$ and
$\mathbf{u}(\mathfrak{g})$.

\begin{lemma} Let $\mathfrak{g}=\mathfrak{g}_{0}\oplus
\mathfrak{g}_{1}$ be a Lie superlagebra and $x_{1},\ldots,x_{s}$ a
basis of $\mathfrak{g}_{1}$, $y_{1},\ldots,y_{t}$ a basis of
$\mathfrak{g}_{0}$. Then

\emph{(1) }  $U(\mathfrak{g})$ has a basis
$$\{x_{1}^{a_{1}}\cdots x_{s}^{a_{s}}y_{1}^{b_{1}}\cdots y_{t}^{b_{t}}|b_{i}\in \mathbb{N},\; a_{j}=0,1
\;\emph{\textrm{for all}} \;i,j\}.$$

\emph{(2)}   $\mathbf{u}(\mathfrak{g})$ has a basis
$$\{x_{1}^{a_{1}}\cdots x_{s}^{a_{s}}y_{1}^{b_{1}}\cdots y_{t}^{b_{t}}|0\leq b_{i}< p,\; a_{j}=0,1
\;\emph{\textrm{for all}} \;i,j\}.$$
\end{lemma}

The following proposition gives  new kinds of PBW basis, which are
suitable for our purpose.

\begin{proposition} Let $\mathfrak{g}=\mathfrak{g}_{0}\oplus \mathfrak{g}_{1}$
 be a Lie superalgebra  and  $x_{1},\ldots,x_{s}$ a
basis of $\mathfrak{g}_{1}$ in which we assume $[x_{i},x_{i}]=0$ for
$i\leq s_{1}$ and $z_{j}:=[x_{j},x_{j}]\neq 0$ for $s_{1}<j\leq s$.
Assume that $z_{s_{1}+1},\ldots, z_{s}$ are linear independent and
denote the subspace of $\mathfrak{g}_{0}$ spanned by them by $V$.
Let $W$ be a subspace of $\mathfrak{g}_{0}$ such that
$\mathfrak{g}_{0}=W\oplus V$ and $y_{1},\ldots,y_{t_{1}}$ be a basis
of $W$. Then

\emph{(1)}  $U(\mathfrak{g})$ has a basis consisting of
$$x_{1}^{a_{1}}\cdots x_{s_{1}}^{a_{s}}x_{s_{1}+1}^{b_{1}}\cdots
x_{s}^{b_{s-s_{1}}}y_{1}^{c_{1}}\cdots y_{t_{1}}^{c_{t_{1}}}$$ where
$0\leq a_{i}< 2,\; b_{j}, c_{k}\in \mathbb{N}\;\textrm{for all}
\;i,j,k.$

\emph{(2)}   $\mathbf{u}(\mathfrak{g})$ has a basis consisting of
$$x_{1}^{a_{1}}\cdots x_{s_{1}}^{a_{s}}x_{s_{1}+1}^{b_{1}}\cdots
x_{s}^{b_{s-s_{1}}}y_{1}^{c_{1}}\cdots y_{t_{1}}^{c_{t_{1}}}$$ where
$0\leq a_{i}< 2,\; 0\leq b_{j}<2p \;,\;0\leq c_{k}< p\;\textrm{for
all} \;i,j,k.$
\end{proposition}
\begin{proof} We only prove (2) since (1) can be proved similarly.
By assumption  the set $\{z_{i},y_{j}|s_{1}<i\leq s, 0\leq j\leq
t_{1}\}$ is a basis of $\mathfrak{g}_{0}$. Owing to Lemma 2.4 (2),
$$\{x_{1}^{a_{1}}\cdots x_{s}^{a_{s}}y_{1}^{b_{1}}\cdots
y_{t_{1}+s-s_{1}}^{b_{t_{1}+s-s_{1}}}|0\leq b_{i}< p,\; a_{j}=0,1
\;\textrm{for all} \;i,j\}$$ is a basis of
$\mathbf{u}(\mathfrak{g})$ where we set
$y_{t_{1}+i}:=z_{i}\;(s_{1}+1\leq i \leq s)$ for consistence. By the
proof of the PBW theorem, there is no any restriction on the order
of elements we choose and thus the following elements also form a
basis of $\mathbf{u}(\mathfrak{g})$:
\begin{equation} x_{1}^{a_{1}}\cdots
x_{s_{1}}^{a_{s_{1}}}
x_{s_{1}+1}^{a_{s_{1}+1}}z_{s_{1}+1}^{b_{s_{1}+1}}\cdots
x_{s}^{a_{s}}z_{s}^{b_{s}}y_{1}^{b_{1}}\cdots y_{t_{1}}^{b_{t_{1}}}
\end{equation}
where $0\leq b_{i}< p,\; a_{j}=0,1$ for all $i,j$. Since
$$z_{i}=[x_{i},x_{i}]=2x_{i}^{2}$$ in $\mathbf{u}(\mathfrak{g})$ for
$s_{1}+1\leq i\leq s$, the set
$\{x_{i}^{a_{i}}z_{i}^{b_{i}}|a_{i}=0,1,\;0\leq b_{i}<
p\}=\{a(m_{i})x_{i}^{m_{i}}|0\leq m_{i}< 2p, \;\textrm{some}\; 0\neq
a(m_{i})\in \kappa\}$. So we can abbreviate elements of (2.1) and
get the ones described in the proposition. The conclusion is proved.
\end{proof}

 Both $U(\mathfrak{g})$ and $\mathbf{u}(\mathfrak{g})$
are super cocommutative Hopf algebras. Thus they are braided Hopf
algebras in
$^{\kappa\mathbb{Z}_{2}}_{\kappa\mathbb{Z}_{2}}\mathscr{YD}$. In
particular, $\mathbf{u}(\mathfrak{g})\# \kappa\mathbb{Z}_{2}$ is an
ordinary algebra. Actually, it is a Hopf algebra by above
subsection. Let $M,N$ be two $\mathbf{u}(\mathfrak{g})\#
\kappa\mathbb{Z}_{2}$-modules and $P_{\bullet}\to M$ be a projective
resolution of $M$. Define
$$\H^{i}_{\mathbf{u}(\mathfrak{g})}(M,N):=\Ext^{i}_{\mathbf{u}(\mathfrak{g})\#
\kappa\mathbb{Z}_{2}}(M,N)=\H^{i}(\Hom_{\mathbf{u}(\mathfrak{g})\#
\kappa\mathbb{Z}_{2}}(P_{\bullet},N)),$$
$$\H^{i}({\mathbf{u}(\mathfrak{g})},M):=\Ext^{i}_{\mathbf{u}(\mathfrak{g})\#
\kappa\mathbb{Z}_{2}}(\kappa,M)\;\;\textrm{and}$$
$$\H^{i}({\mathbf{u}(\mathfrak{g})},\kappa):=\Ext^{i}_{\mathbf{u}(\mathfrak{g})\#
\kappa\mathbb{Z}_{2}}(\kappa,\kappa)$$ for $i\geq 0$, where $\kappa$
is the trivial $\mathbf{u}(\mathfrak{g})\#
\kappa\mathbb{Z}_{2}$-module with the action gotten through the
counit $\varepsilon:\; \mathbf{u}(\mathfrak{g})\#
\kappa\mathbb{Z}_{2}\to \kappa$.

\begin{remark} \emph{By Lemma $2.2$, this is equivalent to say that we consider the
restricted cohomology of a  restricted Lie superalgeba
$\mathfrak{g}$ exactly in the category
$\mathbf{u}(\mathfrak{g})$-smod. That is, we only consider even
homomorphisms. This is totally different with the relative
cohomology defined in \cite{BKN}, where the authors indeed bring all
homomorphisms into consideration.}
\end{remark}

For any coalgebra $C$, we denote $\Ker \varepsilon$ by $C^{+}$ as
usual. Also, as a usual algebra $|\mathbf{u}(\mathfrak{g})|$ has its
usual cohomology $\H^{i}(|\mathbf{u}(\mathfrak{g})|,N)$ for any
$|\mathbf{u}(\mathfrak{g})|$-module $N$. For any Hopf algebra $H$
and $H$-module $M$, we define $M^{H}:=\{m\in M|h\cdot
m=\varepsilon(h)m,\;\textrm{for all}\;h\in H\}$.

\begin{lemma} Let $N$ be a $\mathbf{u}(\mathfrak{g})$-supermodule.
Then for any natural number $i$,
$$\H^{i}(\mathbf{u}(\mathfrak{g}),N)\cong
\H^{i}(|\mathbf{u}(\mathfrak{g})|,N)^{\kappa\mathbb{Z}_{2}}.$$
\end{lemma}
\begin{proof} At first, we prove the conclusion in the case $N=\kappa$.
Note that $|\mathbf{u}(\mathfrak{g})|^{+}$ is the augmentation ideal
of $|\mathbf{u}(\mathfrak{g})|$. Now consider the bar resolution of
$\kappa$
\begin{equation}\cdots\to |\mathbf{u}(\mathfrak{g})|\otimes(|\mathbf{u}(\mathfrak{g})|^{+})^{\otimes 2}
\stackrel{d_{2}}{\to} |\mathbf{u}(\mathfrak{g})|\otimes
|\mathbf{u}(\mathfrak{g})|^{+}\stackrel{d_{1}}{\to}|\mathbf{u}(\mathfrak{g})|\stackrel{\varepsilon}{\to}\kappa\to
0,
\end{equation}
where $d_{i}(a_{0}\otimes\cdots\otimes
a_{i})=\sum_{j=0}^{i-1}(-1)^{j}a_{0}\otimes\cdots\otimes
a_{j}a_{j+1}\otimes\cdots \otimes a_{i}$. Thus every differential
map $d_{i}$ is indeed an even homomorphism. Applying
$\Hom_{|\mathbf{u}(\mathfrak{g})|}(-,\kappa)$, one get
\begin{equation} 0\to \Hom_{\kappa}(\kappa,\kappa)
\stackrel{\delta_{0}}{\to}
\Hom_{\kappa}(|\mathbf{u}(\mathfrak{g})^{+}|,\kappa)\stackrel{\delta_{1}}{\to}\Hom_{\kappa}(|\mathbf{u}(\mathfrak{g})^{+}|^{\otimes2},\kappa)
\stackrel{\delta_{2}}{\to} \cdots,
\end{equation}
where $\delta_{i}=d_{i}^{\ast}$. By definition,
$\H^{i}(|\mathbf{u}(\mathfrak{g})|,\kappa)=\Ker \delta_{i}/\Im
\delta_{i-1}$. Meanwhile, $\H^{i}(\mathbf{u}(\mathfrak{g}),\kappa)$
is exactly the $i$th cohomology of the following complex $$ 0\to
\underline{\Hom}_{\kappa}(\kappa,\kappa) \stackrel{\delta_{0}}{\to}
\underline{\Hom}_{\kappa}(\mathbf{u}(\mathfrak{g})^{+},\kappa)\stackrel{\delta_{1}}{\to}\underline{\Hom}_{\kappa}
((\mathbf{u}(\mathfrak{g})^{+})^{\otimes2},\kappa)
\stackrel{\delta_{2}}{\to} \cdots.
$$ Since $\Hom_{\kappa}((\mathbf{u}(\mathfrak{g})^{+})^{\otimes i},\kappa)^{\kappa\mathbb{Z}_{2}}=\Hom_{\kappa\mathbb{Z}_{2}}
((\mathbf{u}(\mathfrak{g})^{+})^{\otimes
i},\kappa)=\underline{\Hom}_{\kappa}((\mathbf{u}(\mathfrak{g})^{+})^{\otimes
i},\kappa)$, $\H^{i}(\mathbf{u}(\mathfrak{g}),\kappa)\cong
\H^{i}(|\mathbf{u}(\mathfrak{g})|,\kappa)^{\kappa\mathbb{Z}_{2}}$.

In general, for any $\mathbf{u}(\mathfrak{g})$-supermodule $N$, one
can apply $\Hom_{|\mathbf{u}(\mathfrak{g})|}(-,N)$ to (2.2) to get
 a similar complex like (2.3). Using totally the same argument as
$\kappa$, one can get the desired conclusion.
\end{proof}

The following result is a direct consequence of Lemma 2.1 by noting
that $\mathbf{u}(\mathfrak{g})\# \kappa\mathbb{Z}_{2}$ is  an
ordinary Hopf algebra.

\begin{corollary} Let $M$ be an $\mathbf{u}(\mathfrak{g})$-supermodule. Then under cup product,
$\H^{ev}({\mathbf{u}(\mathfrak{g})},\kappa):=\bigoplus_{i\geq 0}
\H^{2i}({\mathbf{u}(\mathfrak{g})},\kappa)$ is  a commutative
algebra and
$\H^{\ast}({\mathbf{u}(\mathfrak{g})},M)\\:=\bigoplus_{i\geq 0}
\H^{i}({\mathbf{u}(\mathfrak{g})},M)$ is an
$\H^{ev}({\mathbf{u}(\mathfrak{g})},\kappa)$-module.
\end{corollary}

\subsection{Basic classical Lie superalgebras}

\begin{definition} \emph{A Lie superalgebra is a \emph{basic classical Lie
superalgebra} if it admits an even nondegenerate supersymmetric
bilinear form and its even subalgebra is reductive.}
\end{definition}

In the following, we only deal with basic classical Lie
superalgebras unless we state otherwise. We recall the list of basic
classical Lie superalgebra (see \cite{Kac1,WZ}). They are four
infinite series
$\mathbf{A}(m,n),\;\mathbf{B}(m,n),\;\mathbf{C}(n),\;\mathbf{D}(m,n)$
and three exceptional versions
$\mathbf{D}(2,1;\alpha),\;\mathbf{G}(3),\;\mathbf{F}(4)$ for
$\alpha\in \kappa\backslash\{0,1\}$. They are still simple Lie
superalgebras even the characteristic of base field is not zero. One
merit of a basic classical Lie superalgebra $\mathfrak{g}$ is that
it admits nice root space decompositions:
$$\mathfrak{g}=\mathfrak{h}\oplus \bigoplus_{\alpha\in \Phi}\mathfrak{g}_{\alpha}$$
such that

(i) $\mathfrak{h}$ is a Cartan subalgebra of $\mathfrak{g}$;

(ii) $\dim_{\kappa}\mathfrak{g}_{\alpha}=1$ for $\alpha \in \Phi$
except
 for $\mathbf{A}(1,1)$;

(iii) Except
 for $\mathbf{A}(1,1)$, $[\mathfrak{g}_{\alpha},\mathfrak{g}_{\beta}]\neq 0$ if and
only if $\alpha, \beta, \alpha+\beta\in \Phi$.

See Section 2.5.3 in \cite{Kac1} for details by noting we still can
do such decompositions in positive characteristic case.  In order to
discriminate different root in characteristic $p$ case, we always
assume $p>3$. Also, we fix a root decomposition just as described in
Section 2.5.4 in \cite{Kac1} from now on. $\Phi$ is called a
\emph{root supersystem} of $\mathfrak{g}$. Clearly,
$\Phi=\Phi_{0}\cup \Phi_{1}$, where $\Phi_{0}$ is the root system of
$\mathfrak{g}_{0}$ and $\Phi_{1}$ is the system of weights of the
representation of $\mathfrak{g}_{0}$ on $\mathfrak{g}_{1}$.
$\Phi_{0}$ is called the \emph{even system} and $\Phi_{1}$ the
\emph{odd system}. Define
$$\Phi_{11}:=\{\alpha\in \Phi_{1}|[\mathfrak{g}_{\alpha},\mathfrak{g}_{\alpha}]=0\},\;\;\;\;
\Phi_{12}:=\{\alpha\in
\Phi_{1}|[\mathfrak{g}_{\alpha},\mathfrak{g}_{\alpha}]\neq0\}.$$ By
observing the root supersystem of $\mathbf{B}(m,n)$, $\Phi_{12}\neq
\phi$ in general.

\begin{lemma} Let $\mathfrak{g}$ be a basic classical Lie superalgebra.
Then for any $\alpha\in \Phi_{0}$ and  $x\in \mathfrak{g}_{\alpha}$,
$x^{p}=0$ in $\mathbf{u}(\mathfrak{g})$.
\end{lemma}
\begin{proof} This should be known, but the author can not find
suitable reference. So we give a short proof here. It is known that
the even part $\mathfrak{g}_{0}$ of a basic classical Lie
superalgebra $\mathfrak{g}$ is a direct sum of some Lie algebras of
types $\mathbf{A}_{n},\mathbf{B}_{n},\mathbf{C}_{n},\mathbf{D}_{n},
\mathbf{G}_{2}$ and $\kappa$. Therefore there is no harm to assume
that $\mathfrak{g}_{0}$ is a simple Lie algebra of type
$\mathbf{A}_{n},\mathbf{B}_{n},\mathbf{C}_{n},\mathbf{D}_{n}$ or
$\mathbf{G}_{2}$. So $\mathfrak{g}_{0}$ is generated by
$\mathfrak{sl}_{2}$-triples $\{e_{i},f_{i},h_{i}|i\in I\}$. Thus
firstly we
 assume that $x=e_{i}$ or $x=f_{i}$ for some $i$. Say, $x=e_{i}$.
Note that $e_{i}$ commutes with all $f_{j}$ unless $j=i$ and in this
case $ad(e_{i})^{3}(f_{i})=0$. So $ad(e_{i})^{p}(f_{j})=0$ for all
$j\in I$. From Serre's relation, $ad(e_{i})^{1-a_{ij}}(e_{j})=0$
where $(a_{ij})_{I\times I}$ is the Cartan matrix of
$\mathfrak{g}_{0}$. This implies $ad(e_{i})^{p}(e_{j})=0$ for all
$j\in I$ since $p> 1-a_{ij}$ by our assumption on $p$. Also clearly
$ad(e_{i})^{p}(h_{j})=0$ for all $j\in I$. By the definition of
restricted Lie algebra, $x^{[p]}$ lies in the center of
$\mathfrak{g}_{0}$ and so $x^{[p]}=0$, which implies $x^{p}=0$ in
$\mathbf{u}(\mathfrak{g})$ too. The case $x=f_{i}$ can be proved
similarly. For general $x\in \mathfrak{g}_{\alpha}$, it is well
known that up to a scalar we can get $x$ by applying the Lie algebra
automorphisms $\tau_{j}:=exp(ad (e_{j}))exp(ad (-f_{j}))exp(ad
(e_{j}))$ iteratively to some $e_{i}$ or $f_{i}$. Thus $x^{[p]}$
lies in the center too.
\end{proof}

There is a filtration on $\mathbf{u}(\mathfrak{g})$ with degrees
$$\deg_{1}(\mathfrak{h})=0,\;\;\deg_{1}(\bigoplus_{\alpha\in
\Phi_{1}}\mathfrak{g}_{\alpha})=1,\;\; \deg_{1}(\bigoplus_{\alpha\in
\Phi_{0}}\mathfrak{g}_{\alpha})=2.$$ The associated graded algebra
is denoted by $\Gr^{1}(\mathbf{u}(\mathfrak{g}))$. It is still a
super cocommutative Hopf algebra. It is not hard to see that there
is a natural projection from $\Gr^{1}(\mathbf{u}(\mathfrak{g}))$ to
$\mathbf{u}(\mathfrak{h})$ and thus there is a subsuperalgebra
$R_{\mathfrak{g}}$ such that
$$\Gr^{1}(\mathbf{u}(\mathfrak{g}))=R_{\mathfrak{g}}\# \mathbf{u}(\mathfrak{h}).$$
Actually,  $R_{\mathfrak{g}}$ is the graded subalgebra generated by
$\bigoplus_{\alpha\in \Phi}\mathfrak{g}_{\alpha}$.

For any set $S$, its cardinal number is denoted by $S^{\#}$. Assume
that $\mathfrak{g}\neq \mathbf{A}(1,1)$. Then by property (ii) of
the root space decomposition, up to scalers there is a unique
nonzero element $x_{\alpha}$ belonging to $\mathfrak{g}_{\alpha}$.

\begin{lemma} Assume that $\mathfrak{g}\neq \mathbf{A}(1,1)$ and let
$x_{\alpha}$ defined as above. Then the graded algebra
$R_{\mathfrak{g}}$ has the following PBW basis consisting of
elements
\begin{equation} x_{\alpha_{1}}^{a_{1}}\cdots x_{\alpha_{r}}^{a_{r}}x_{\beta_{1}}^{b_{1}}
\cdots x_{\beta_{s}}^{b_{s}}x_{\gamma_{1}}^{c_{1}}\cdots
x_{\gamma_{s}}^{c_{t}}
\end{equation}
where $\alpha_{i}\in \Phi_{11}, \beta_{j}\in \Phi_{12},
\gamma_{k}\in \Phi_{0}$,
$r=\Phi_{11}^{\#},s=\Phi_{12}^{\#},t=\Phi_{0}^{\#}-\Phi_{12}^{\#}$
and $0\leq a_{i}< 2, 0\leq b_{j}< 2p, 0\leq c_{k}< p$ for $1\leq
i\leq r, 1\leq j\leq s, 1\leq k\leq t$.
\end{lemma}
\begin{proof} Under the grading $\Gr^{1}$, one can see that $$[\bigoplus_{\alpha\in \Phi}\mathfrak{g}_{\alpha},
\bigoplus_{\alpha\in \Phi}\mathfrak{g}_{\alpha}]\subseteq
\bigoplus_{\alpha\in \Phi_{0}}\mathfrak{g}_{\alpha}.$$ So to show
the conclusion, we can assume that $\bigoplus_{\alpha\in
\Phi}\mathfrak{g}_{\alpha}$ is a Lie subsuperalgebra of
$\mathfrak{g}$. Being living in different root spaces,
$\{[x_{\beta_{j}},x_{\beta_{j}}]|1\leq j\leq s\}$ are linear
independent. So Proposition 2.5 can be applied and thus the set of
elements in (2.4) forms a basis of $\mathbf{u}(\bigoplus_{\alpha\in
\Phi}\mathfrak{g}_{\alpha})$. Clearly such elements are homogeneous
in $R_{\mathfrak{g}}$ and so they also give a basis of
$R_{\mathfrak{g}}$.
\end{proof}

Throughout the following of this subsection, we always assume that
$\mathfrak{g}\neq \mathbf{A}(1,1)$. In order to reduce
$R_{\mathfrak{g}}$ to a familiar algebra, we introduce  another kind
of filtration on $R_{\mathfrak{g}}$. To attack it, the degree of an
element in (2.4)
 is defined to be
$$\deg_{2}(x_{\alpha_{1}}^{a_{1}}\cdots
x_{\gamma_{s}}^{c_{t}})=(a_{1},\ldots,a_{r},b_{1},\ldots,b_{s},c_{1},\ldots,c_{t})\in
\mathbb{N}^{\Phi^{\#}}$$ and totally order the elements (2.4)
lexicographically by setting
$$(1,0,\ldots,0)>\cdots>(0,1,\ldots,0)>\cdots>(0,0,\ldots,1).$$
For convenience and consistence, we set
$\alpha_{r+i}:=\beta_{i}\;(1\leq i\leq s)$ and
$\alpha_{r+s+i}:=\gamma_{i}\;(1\leq i\leq t)$.

\begin{lemma} Under the total order defined above, for all $i< j$,
$$\deg_{2}([x_{\alpha_{i}},x_{\alpha_{j}}])< \deg_{2}(x_{\alpha_{i}}x_{\alpha_{j}})$$
unless $[x_{\alpha_{i}},x_{\alpha_{j}}]=0$.
\end{lemma}
\begin{proof} It is not hard to see that any $x\in \bigoplus_{\alpha\in
\Phi_{0}}$ actually lies in the center of $R_{\mathfrak{g}}$. So to
show the lemma, one can assume that  both $x_{\alpha_{i}}$ and
$x_{\alpha_{i}}$ are odd elements and
$[x_{\alpha_{i}},x_{\alpha_{j}}]\neq 0$. Now
$[x_{\alpha_{i}},x_{\alpha_{j}}]$ lies in $\mathfrak{g}_{0}$
automatically and thus either
$[x_{\alpha_{i}},x_{\alpha_{j}}]=cx_{\alpha_{l}}$ for $l>j$ and
$0\neq c\in \kappa$ or
$[x_{\alpha_{i}},x_{\alpha_{j}}]=d[x_{\alpha_{k}},x_{\alpha_{k}}]$
for some odd element with $[x_{\alpha_{k}},x_{\alpha_{k}}]\neq 0$
and $0\neq d\in \kappa$. In the first case, the conclusion is clear.
In the second case, we still need to consider two cases:
$[x_{\alpha_{i}},x_{\alpha_{i}}]=[x_{\alpha_{j}},x_{\alpha_{j}}]=0$
or either of them is not zero. Also, the first case implies that
$j<k$ by the PBW basis we choose and thus the conclusion is proved.
By property (iii) of the root space decomposition, $\alpha_{i}+
\alpha_{j}$ is still a root and it is equals to $2\alpha_{k}$ by
assumption. Comparing with the root supersystem listed in Section
2.5.4 in \cite{Kac1}, this is happened only in the case
$[x_{\alpha_{i}},x_{\alpha_{i}}]=[x_{\alpha_{j}},x_{\alpha_{j}}]=0$.
\end{proof}

By Lemma 2.12, the above ordering induces a filtration on
$R_{\mathfrak{g}}$. The associated graded algebra is denoted by
$\Gr^{2}(R_{\mathfrak{g}})$. It is generated by
$\{x_{\alpha_{i}}|1\leq i\leq \Phi^{\#}-\Phi_{12}^{\#}\}$ with
relations
\begin{equation}[x_{\alpha_{i}},x_{\alpha_{j}}]=0\;\;\textrm{for}\;i\neq
j,\;\;\;\;x_{\alpha_{i}}^{N_{i}}=0
\end{equation}
where
$$N_{i}=\left \{
\begin{array}{lll} 2, & \;\;\;\;0\leq i\leq \Phi_{11}^{\#}\\
  2p, &
\;\;\;\;\Phi_{11}^{\#}+1\leq i\leq \Phi_{11}^{\#}+\Phi_{12}^{\#}\\
p,&\;\;\;\;\Phi_{11}^{\#}+\Phi_{12}^{\#}+1\leq i\leq
\Phi^{\#}-\Phi_{12}^{\#}.
\end{array}\right. $$

Note that $\Gr^{2}(R_{\mathfrak{g}})$ inherits the action of
$\mathbf{u}(\mathfrak{h})$ from that  on $R_{\mathfrak{g}}$
naturally, define
$$\Gr^{2}(\mathbf{u}(\mathfrak{g})):=\Gr^{2}(R_{\mathfrak{g}})\# \mathbf{u}(\mathfrak{h}).$$

\subsection{Spectral sequences and finite generation}

 We will see in the next section that there are some convergent
spectral sequences associated to the filtrations given in Subsection
2.3. The following lemma, which is essentially used in this paper,
is given in \cite{MPSW} as its Lemma 2.5. Recall that an element
$a\in E_{r}^{p,q}$ is called a \emph{permanent cycle} if
$d_{i}(a)=0$ for all $i\geq r$.

\begin{lemma} \emph{(1)} Let $E_{1}^{p,q}\Rightarrow E_{\infty}^{p+q}$ be a
multiplicative spectral sequence of $\kappa$-algebras concentrated
in the half plane $p+q\geq 0$, and let $A^{\ast,\ast}$ be a bigraded
commutative $\kappa$-algebra concentrated in even
\emph{(}total\emph{)} degrees. Assume that there exists a bigraded
map of algebras $\varphi:A^{\ast,\ast}\to E_{1}^{\ast,\ast}$ such
that

\emph{(i)} $\varphi$ makes $E_{1}^{\ast,\ast}$ into a Noetherian
$A^{\ast,\ast}$-module, and

\emph{(ii)} the image of $A^{\ast,\ast}$ in $E_{1}^{\ast,\ast}$
consists of
permanent cycles.\\[1.5mm]
Then $E_{\infty}^{\ast}$ is a Noetherian module over
\emph{Tot($A^{\ast,\ast}$)}.

\emph{(2)} Let $\tilde{E}_{1}^{p,q}\Rightarrow
\tilde{E}_{\infty}^{p+q}$ be a spectral sequence that is a bigraded
module over the spectral sequence $E^{\ast,\ast}$. Assume that
$\tilde{E}_{1}^{\ast,\ast}$ is a Noetherian module over
$A^{\ast,\ast}$ where $A^{\ast,\ast}$ acts on
$\tilde{E}_{1}^{\ast,\ast}$ via the map $\varphi$. Then
$\tilde{E}_{\infty}^{\ast}$ is a finitely generated
$E_{\infty}^{\ast}$-module.
\end{lemma}

\section{Finite generation}

The following conclusion is one of main results of this paper.

\begin{theorem} Let $\mathfrak{g}$ be one of basic classical Lie
superalgebras over $\kappa$ and $\mathbf{u}(\mathfrak{g})$ its
restricted enveloping algebra. Then

\emph{(1)} the algebra
$\H^{\ast}(\mathbf{u}(\mathfrak{g}),\kappa):=\bigoplus_{i\geq
0}\H^{i}(\mathbf{u}(\mathfrak{g}),\kappa)$ is finitely generated.

\emph{(2)} $\H^{\ast}(\mathbf{u}(\mathfrak{g}),M)$ is a finitely
generated module over $\H^{\ast}(\mathbf{u}(\mathfrak{g}),\kappa)$
for $M$ a finitely generated $\mathbf{u}(\mathfrak{g})$-supermodule.
\end{theorem}

We will divide the proof into two cases: $\mathfrak{g}\neq
\mathbf{A}(1,1)$ or $\mathfrak{g}= \mathbf{A}(1,1)$. The basic idea
of the proof is to modify the procedure developed in \cite{MPSW}
into our cases by applying preliminary results gotten in Section 2.
Firstly, $\mathfrak{g}\neq \mathbf{A}(1,1)$ is assumed until
Subsection 3.4.

\subsection{Cohomology of $\Gr^{2}(\mathbf{u}(\mathfrak{g}))$}

The algebraic structure of $\Gr^{2}(R_{(\mathfrak{g})})$ has been
described clearly in (2.5). Recall that we denote the usual algebra
of superalgebra $A$ by $|A|$. For continuation, we write the
algebraic structure of $|\Gr^{2}(R_{(\mathfrak{g})})|$ again as
follows: it is generated by $\{x_{\alpha_{i}}|1\leq i\leq
\Phi^{\#}-\Phi_{12}^{\#}\}$ with relations
\begin{equation}x_{\alpha_{i}}x_{\alpha_{j}}=\left \{
\begin{array}{ll} -x_{\alpha_{j}}x_{\alpha_{i}}, & \;\;\;\;1\leq i<j \leq \Phi_{1}^{\#}\\
 x_{\alpha_{j}}x_{\alpha_{i}}, &
\;\;\;\;1\leq i<j\;\textrm{and}\;j>\Phi_{1}^{\#},
\end{array}\right.  \;\;x_{\alpha_{i}}^{N_{i}}=0
\end{equation}
where
$$N_{i}=\left \{
\begin{array}{lll} 2, & \;\;\;\;0\leq i\leq \Phi_{11}^{\#}\\
  2p, &
\;\;\;\;\Phi_{11}^{\#}+1\leq i\leq \Phi_{11}^{\#}+\Phi_{12}^{\#}\\
p,&\;\;\;\;\Phi_{11}^{\#}+\Phi_{12}^{\#}+1\leq i\leq
\Phi^{\#}-\Phi_{12}^{\#}.
\end{array}\right. $$

The algebra $|\Gr^{2}(R_{(\mathfrak{g})})|$ is a special case of
so-called \emph{quantum complete intersection algebras}: Let $N$ be
positive integer, and for each $i\in\{1,\ldots,N\}$, $N_{i}$ be an
integer greater than $1$. Let $q_{ij}\in
\kappa^{\ast}=\kappa\backslash \{0\}$ for $1\leq i<j\leq N$. Define
$S$ to be the $\kappa$-algebra generated by $x_{1},\ldots,x_{N}$
subject to the relations
\begin{equation} x_{i}x_{j}=q_{ij}x_{j}x_{i}\;\;\textrm{for all}
\;i<j\;\;\textrm{and}\;\;x_{i}^{N_{i}}=0\;\textrm{for all} \;i.
\end{equation}
$S$ is called a quantum complete intersection algebra. For such $S$,
its cohomology ring $\H^{\ast}(S,\kappa)=\bigoplus_{i\geq
0}\Ext^{i}_{S}(\kappa,\kappa)$ was determined in Section 4 of
\cite{MPSW}. For completeness and consistence of the paper, let us
sketch it.

Let $K_{\bullet}$ be the following complex of free $S$-modules. For
each $N$-tuple $(a_{1},\ldots,a_{N})$ of nonnegative integers, let
$\Psi(a_{1},\ldots,a_{N})$ be a free generator in degree
$a_{1}+\cdots+a_{N}$. Then define
$K_{n}=\oplus_{a_{1}+\cdots+a_{N}=n}S\Psi(a_{1},\ldots,a_{N})$. For
each $i\in\{1,\ldots,N\}$, let $\sigma_{i},\tau_{i}:\;\mathbb{N}\to
\mathbb{N}$ be the function defined by
$$\sigma_{i}(a)=\left \{
\begin{array}{ll} 1, & \;\;\;\;a\;\textrm{is odd}\\
  N_{i}-1, &
\;\;\;\;a\; \textrm{is even},
\end{array}\right. $$
and $\tau_{i}(a)=\sum_{j=1}^{a}\sigma_{i}(a)$ for $a\geq 1$,
$\tau(0)=0$. Let
$$d_{i}(\Psi(a_{1},\ldots,a_{N}))=(\prod_{l<i}(-1)^{a_{l}}q_{li}^{\sigma_{i}(a_{i})\tau_{l}(a_{l})})x_{i}^{\sigma_{i}(a_{i})}
\Psi(a_{1},\ldots,a_{i}-1,\ldots, a_{N})$$ if $a_{i}>0$, and
$d_{i}(\Psi(a_{1},\ldots,a_{N}))=0$ if $a_{i}=0$. Extend each
$d_{i}$ to an $S$-module homomorphism and set
$$d=d_{1}+\cdots+d_{N}.$$
It is shown in Section 4 of \cite{MPSW} that $(K_{\bullet},d)$ is a
resolution of $\kappa$.

From this resolution, one can compute $\Ext^{i}_{S}(\kappa,\kappa)$.
Applying $\Hom_{S}(-,\kappa)$ to $K_{\bullet}$, the induced
differential $d^{\ast}$ is the zero map (since
$x_{i}^{\sigma_{i}(a_{i})}$ is always in the augmentation ideal) and
thus the cohomology is just the complex
$\Hom_{S}(K_{\bullet},\kappa)$. Now let $\xi_{i}\in
\Hom_{S}(K_{2},\kappa),\;\eta_{i}\in \Hom_{S}(K_{1},\kappa)$ be the
functions dual to $\Psi(0,\ldots,2,\ldots, 0)$ (the $2$ in the $i$th
place) and $\Psi(0,\ldots,1,\ldots, 0)$ (the $1$ in the $i$th place)
respectively. The following conclusion is the Theorem 4.1 in
\cite{MPSW}.

\begin{lemma} The algebra $\H^{\ast}(S,\kappa)$ is generated by
$\xi_{i},\eta_{i}$ $(1\leq i\leq N)$ with  $\deg \xi_{i}=2$ and
$\deg \eta_{i}=1$, subject to the relations
$$\xi_{i}\xi_{j}=q_{ij}^{N_{i}N_{j}}\xi_{j}\xi_{i},\;\;\eta_{i}\xi_{j}=q_{ji}^{N_{j}}\xi_{j}\eta_{i},\;\;
\eta_{i}\eta_{j}=-q_{ji}\eta_{j}\eta_{i}$$ where
$q_{ij}=q_{ji}^{-1}$ if $i>j$.
\end{lemma}

For any two nonnegative integers $m,n$, define an algebra
$\wedge(m|n)$ as follows. It is generated by
$\eta_{1},\ldots,\eta_{m+n}$ with relations
$$\eta_{i}\eta_{j}=\left \{
\begin{array}{ll} \eta_{j}\eta_{i}, & \;\;\;\;1\leq i<j \leq m\\
 -\eta_{j}\eta_{i}, &
\;\;\;\;1\leq i<j\;\textrm{and}\;j>m,
\end{array}\right.  \;\;\eta_{i}^{2}=0.
$$

\begin{proposition} Let $\mathfrak{g}$ be a basic classical Lie
superalgebra different from $\mathbf{A}(1,1)$ and $\Phi$ its root
supersystem. Then
$$\H^{\ast}(|\Gr^{2}(R_{\mathfrak{g}})|,\kappa)\cong \kappa [\xi_{1},\ldots,\xi_{m+n}]\otimes \wedge(m|n)$$
where $m=\Phi_{1}^{\#}, n=\Phi_{0}^{\#}-\Phi_{12}^{\#}$ and $\deg
\xi_{i}=2,\;\deg \eta_{i}=1$.
\end{proposition}
\begin{proof} It is a direct consequence of Lemma 3.2 and the
definition of $|\Gr^{2}(R_{\mathfrak{g}})|$.
\end{proof}

\begin{proposition}  Let $\mathfrak{g}$ be a basic classical Lie
superalgebra different from $\mathbf{A}(1,1)$. Fix notions as above.
Then

\emph{(1)}
$\H^{\ast}(|\Gr^{2}(\mathbf{u}(\mathfrak{g}))|,\kappa)\cong
\H^{\ast}(|\Gr^{2}(R_{\mathfrak{g}})|,\kappa)^{\mathbf{u(\mathfrak{h})}}$
where the action of $\mathbf{u}(\mathfrak{h})$ on
$\H^{\ast}(|\Gr^{2}(R_{\mathfrak{g}})|,\kappa)$ is given through
\begin{equation} h\cdot \xi_{i}=-N_{i}\alpha_{i}(h)\xi_{i},\;\; h\cdot
\eta_{i}=-\alpha_{i}(h)\eta_{i},
\end{equation}
for $1\leq i\leq \Phi^{\#}-\Phi_{12}^{\#}$ and $h\in
\mathbf{u}(\mathfrak{h})$.

\emph{(2)} $\H^{\ast}(\Gr^{2}(\mathbf{u}(\mathfrak{g})),\kappa)\cong
\H^{\ast}(|\Gr^{2}(R_{\mathfrak{g}})|,\kappa)^{\mathbf{u(\mathfrak{h})}\otimes
 \kappa\mathbb{Z}_{2}}$ where the action of
$\kappa\mathbb{Z}_{2}=\kappa\langle g|g^{2}=1\rangle$ on
$\H^{\ast}(|\Gr^{2}(R_{\mathfrak{g}})|,\kappa)$ is given through
\begin{equation} g\cdot \xi_{i}=\xi_{i},\;\; g\cdot
\eta_{i}=\left \{
\begin{array}{ll} -\eta_{i}, & \;\;\;\;i\leq \Phi_{1}^{\#}\\
\eta_{i}, & \;\;\;\;i> \Phi_{1}^{\#}.
\end{array}\right.
\end{equation}
\end{proposition}
\begin{proof} (1) To give the action  $\mathbf{u}(\mathfrak{h})$ on
$\H^{\ast}(|\Gr^{2}(R_{\mathfrak{g}})|,\kappa)$, we explain
$\xi_{i}, \eta_{i}$ and $h\in \mathbf{u}(\mathfrak{h})$ as chain
maps $K_{\bullet}\to K_{\bullet}$. Then action is given by forming
the commutators of compositions of these chain maps.  In fact,
$\xi_{i}, \eta_{i}$ has been explained as chain maps in \cite{MPSW}
and they are described as follows:
$$\xi_{i}(\Psi(a_{1},\ldots,a_{N}))=\prod_{i<l}q_{il}^{N_{i}\tau_{l}(a_{l})}\Psi(a_{1},\ldots,a_{i}-2,\ldots,a_{N}),$$
$$\eta_{i}(\Psi(a_{1},\ldots,a_{N}))=cx_{i}^{\sigma_{i}(a_{i})-1}\Psi(a_{1},\ldots,a_{i}-1,\ldots,a_{N})$$
where $c=\prod_{l<i}q_{li}^{(\sigma_{i}(a_{i})-1)\tau_{l}(a_{l})}
\prod_{i<l}(-1)^{a_{l}}q_{il}^{\tau_{l}(a_{l})}$ and
$N=\Phi^{\#}-\Phi_{12}^{\#}$. Now let $h$ be an element in
$\mathbf{u}(\mathfrak{h})$. Then $h\cdot \Psi(0,\ldots,1,\ldots,0)$
(the $1$ is in the $i$th place) should equal to
$\alpha_{i}(h)\Psi(0,\ldots,1,\ldots,0)$ (since one can regard
$\Psi(0,\ldots,1,\ldots,0)$ as the generator $x_{\alpha_{i}}$).
Extend it to higher items and one can verify directly the following
extension of $\mathbf{u}(\mathfrak{h})$ on $K_{\bullet}$ indeed
commutes with the differentials:
$$h\cdot \Psi(a_{1},\ldots,a_{N})=\sum_{l=1}^{N}\tau_{l}(a_{l})\alpha_{l}(h) \Psi(a_{1},\ldots,a_{N})$$
for $h\in \mathbf{u}(\mathfrak{h})$ and $a_{1},\ldots,a_{N}\geq 0$.
Then the induced action of $\mathbf{u}(\mathfrak{h})$ on generators
$\xi_{i},\eta_{i}$ is given by
$$h\cdot \xi_{i}=h\xi_{i}-\xi_{i}h=-N_{i}\alpha_{i}(h)\xi_{i},\;\; h\cdot
\eta_{i}=h\eta_{i}-\eta_{i}h=-\alpha_{i}(h)\eta_{i}$$ for $h\in
\mathbf{u}(\mathfrak{h})$.

As $\mathbf{u}(\mathfrak{h})$ is a commutative semisimple algebra,
we indeed have
$$\Ext^{i}_{|\Gr^{2}(\mathbf{u}({\mathfrak{g}}))|}(\kappa,\kappa)=
\Ext^{i}_{|\Gr^{2}(R_{\mathfrak{g}})|\#
\mathbf{u}(\mathfrak{h})}(\kappa,\kappa)\cong
\Ext^{i}_{|\Gr^{2}(R_{\mathfrak{g}})|}(\kappa,\kappa)^{\mathbf{u(\mathfrak{h})}}$$
for $i\geq 0$ (one can prove this fact similarly by applying the
methods used in the proof of Lemma 2.7). Thus
$\H^{\ast}(|\Gr^{2}(\mathbf{u}(\mathfrak{g}))|,\kappa)\cong
\H^{\ast}(|\Gr^{2}(R_{\mathfrak{g}})|,\kappa)^{\mathbf{u(\mathfrak{h})}}$
now.

(2) By Lemma 2.7 and (1),
$$ \H^{\ast}(\Gr^{2}(\mathbf{u}(\mathfrak{g})),\kappa)\cong
\H^{\ast}(|\Gr^{2}(\mathbf{u}(\mathfrak{g}))|,\kappa)^{
\kappa\mathbb{Z}_{2}}\cong
\H^{\ast}(|\Gr^{2}(R_{\mathfrak{g}})|,\kappa)^{\mathbf{u(\mathfrak{h})}\otimes
 \kappa\mathbb{Z}_{2}}.$$ Similar to (1), the following action of
 $\kappa\mathbb{Z}_{2}$ on $K_{\bullet}$ commutes with the
 differentials:
$$g\cdot \Psi(a_{1},\ldots,a_{N})=\prod_{l=1}^{\Phi_{1}^{\#}}(-1)^{\tau_{l}(a_{l})}\Psi(a_{1},\ldots,a_{N}).$$
This induces the action
$$g\cdot \xi_{i}=g\xi_{i}g^{-1}=\left \{
\begin{array}{ll} (-1)^{N_{i}}\xi_{i},&\;i\leq \Phi_{1}^{\#}\\
\xi_{i}, & \;i> \Phi_{1}^{\#},
\end{array}\right.\;\; g\cdot
\eta_{i}=g\eta_{i}g^{-1}=\left \{
\begin{array}{ll} -\eta_{i},&\;i\leq \Phi_{1}^{\#}\\
\eta_{i}, & \;i> \Phi_{1}^{\#}.
\end{array}\right.$$
By the definition of $N_{i}$ in (3.1), it is an even when $i\leq
\Phi_{1}^{\#}$.
\end{proof}

\subsection{Cohomology of $\Gr^{1}(\mathbf{u}(\mathfrak{g}))$}

For a basic classical Lie superalgebra $\mathfrak{g}$, its
enveloping algebra is denoted by $U(\mathfrak{g})$. As the case of
$\mathbf{u}(\mathfrak{g})$, define
$$\deg_{1}(\mathfrak{h}):=0,\;\;\deg_{1}(\bigoplus_{\alpha\in
\Phi_{1}}\mathfrak{g}_{\alpha}):=1,\;\;
\deg_{1}(\bigoplus_{\alpha\in \Phi_{0}}\mathfrak{g}_{\alpha}):=2.$$
Then we will get a filtration on $U(\mathfrak{g})$ and associated
graded algebra
$$\Gr^{1}(U(\mathfrak{g}))=\tilde{R}_{\mathfrak{g}}\# U(\mathfrak{h})$$
similarly.

\begin{lemma} Assume that $\mathfrak{g}\neq \mathbf{A}(1,1)$ and let
$x_{\alpha}$ defined as in Lemma $2.11$. Then the graded algebra
$\tilde{R}_{\mathfrak{g}}$ has the following PBW basis consisting of
elements
\begin{equation} x_{\alpha_{1}}^{a_{1}}\cdots x_{\alpha_{r}}^{a_{r}}x_{\beta_{1}}^{b_{1}}
\cdots x_{\beta_{s}}^{b_{s}}x_{\gamma_{1}}^{c_{1}}\cdots
x_{\gamma_{s}}^{c_{t}}
\end{equation}
where $\alpha_{i}\in \Phi_{11}, \beta_{j}\in \Phi_{12},
\gamma_{k}\in \Phi_{0}$,
$r=\Phi_{11}^{\#},s=\Phi_{12}^{\#},t=\Phi_{0}^{\#}-\Phi_{12}^{\#}$
and $0\leq a_{i}< 2,  b_{j}, c_{k}\in \mathbb{N}$ for $1\leq i\leq
r, 1\leq j\leq s, 1\leq k\leq t$.
\end{lemma}
\begin{proof} Similar to that of Lemma 2.11.
\end{proof}

Also, we set $\alpha_{r+i}:=\beta_{i}\;(1\leq i\leq s)$ and
$\alpha_{r+s+i}:=\gamma_{i}\;(1\leq i\leq t)$ for convenience and
consistence. Clearly,
$$R_{\mathfrak{g}}\cong \tilde{R}_{\mathfrak{g}}/(x_{\alpha_{i}}^{N_{i}},\;1\leq i\leq \Phi^{\#}-\Phi_{12}^{\#})$$
where $N_{i}$ is defined the same as in (2.5). Define
$N:=\Phi^{\#}-\Phi_{12}^{\#}$ and for any
$\mathbf{a}=(a_{1},\ldots,a_{N})\in \mathbb{N}^{N}$ with $0\leq
a_{i}<2\;(1\leq i\leq \Phi_{11}^{\#})$, denote the corresponding PBW
basis element $x_{\alpha_{1}}^{a_{1}}\cdots x_{\alpha_{N}}^{a_{N}}$
by $\mathbf{x^{a}}$ for short.

Our next aim is to give some elements of
$\H^{2}(|R_{\mathfrak{g}}|,\kappa)$. Recall
$|\tilde{R}_{\mathfrak{g}}|^{+}$ is the augmentation ideal of
$|\tilde{R}_{\mathfrak{g}}|$. Now for each $i\in \{1,\ldots,N\}$,
define
$\tilde{\xi}_{\alpha_{i}}:|\tilde{R}_{\mathfrak{g}}|^{+}\otimes
|\tilde{R}_{\mathfrak{g}}|^{+}\to \kappa$ by
$$\tilde{\xi}_{\alpha_{i}}(\mathbf{x^{a}},\mathbf{x^{b}})=c_{\alpha_{i}}$$
where $c_{\alpha_{i}}$ is the coefficient of
$x_{\alpha_{i}}^{N_{i}}$ in the product
$\mathbf{x^{a}}\mathbf{x^{b}}$ as a linear combination of PBW basis
elements. By its definition, $\tilde{\xi}_{\alpha_{i}}$ is
associative on $|\tilde{R}_{\mathfrak{g}}|^{+}$ and thus it may be
extended to a normalized two-cocycle on
$|\tilde{R}_{\mathfrak{g}}|$. We next show that
$\tilde{\xi}_{\alpha_{i}}$ factors through the quotient map
$\pi:|\tilde{R}_{\mathfrak{g}}| \to |R_{\mathfrak{g}}|$ to give a
nonzero two-cocycle on $|R_{\mathfrak{g}}|$. To attack this, we need
show the $\tilde{\xi}_{\alpha_{i}}(\mathbf{x^{a}},\mathbf{x^{b}})=0$
whenever $\mathbf{x^{a}}$ or $\mathbf{x^{b}}$ is in the kernel of
the quotient map $\pi$. Suppose $\mathbf{x^{a}}\in\Ker\pi$, which
implies that $a_{j}\geq N_{j}$ for some $j$. By the proof of Lemma
2.10, $x_{\alpha_{j}}^{N_{j}}$ lies in the center of
$U(\mathfrak{g})$ and so
$\mathbf{x^{a}}=x_{\alpha_{j}}^{N_{j}}\mathbf{x^{c}}$ for some
$\mathbf{c}\in \mathbb{N}^{N}$. Then
$\tilde{\xi}_{\alpha_{i}}(\mathbf{x^{a}},\mathbf{x^{b}})=
\tilde{\xi}_{\alpha_{i}}(x_{\alpha_{j}}^{N_{j}}\mathbf{x^{c}},\mathbf{x^{b}})$
is the the coefficient of $x_{\alpha_{i}}^{N_{i}}$ in the product
$x_{\alpha_{j}}^{N_{j}}\mathbf{x^{c}}\mathbf{x^{b}}$. It is zero
now: If $j=i$, then since $\mathbf{x^{b}}\in
|\tilde{R}_{\mathfrak{g}}|^{+}$, this product cannot have a nonzero
coefficient for $x_{\alpha_{i}}^{N_{i}}$. If $j\neq i$, the same
conclusion is true since $x_{\alpha_{j}}^{N_{j}}$ is always a factor
of $\mathbf{x^{a}}\mathbf{x^{b}}$. One can show the result similarly
in the case $\mathbf{x^{b}}\in \Ker \pi$.

Choose the section $\widetilde{-}: |R_{\mathfrak{g}}|\to
|\tilde{R}_{\mathfrak{g}}|$ of $\pi$ which just sent the PBW basis
elements in $R_{\mathfrak{g}}$, given in Lemma 2.11, to the same
elements in $\tilde{R}_{\mathfrak{g}}$, described in Lemma 3.5.
Since $\tilde{\xi}_{\alpha_{i}}$ factors through
$\pi:|\tilde{R}_{\mathfrak{g}}| \to |R_{\mathfrak{g}}|$, we may
define $\hat{\xi}_{\alpha_{i}}:|R_{\mathfrak{g}}|^{+}\otimes
|R_{\mathfrak{g}}|^{+}\to \kappa$ by
$$\hat{\xi}_{\alpha_{i}}(\mathbf{x^{a}},\mathbf{x^{b}}):=
\tilde{\xi}_{\alpha_{i}}(\tilde{\mathbf{x}}^{\mathbf{a}},\tilde{\mathbf{x}}^{\mathbf{b}})$$
where
$\tilde{\mathbf{x}}^{\mathbf{a}},\tilde{\mathbf{x}}^{\mathbf{b}}$
are defined via the section $\widetilde{-}$.

\begin{proposition} The set
$\{\hat{\xi}_{\alpha_{i}}|i=1,\ldots,N\}$ represents a linear
independent subset of $\H^{2}(|R_{\mathfrak{g}}|,\kappa)$.
\end{proposition}
\begin{proof} At first, let us show that every $\hat{\xi}_{\alpha_{i}}$ is a
2-cocycle. For this, it is enough to show that it is associative,
that is, for any three PBW basis elements
$\mathbf{{x}^{a}},\mathbf{{x}^{b}},\mathbf{{x}^{c}}$, we have
$\hat{\xi}_{\alpha_{i}}(\mathbf{{x}^{a}}\mathbf{{x}^{b}},\mathbf{{x}^{c}})=
\hat{\xi}_{\alpha_{i}}(\mathbf{{x}^{a}},\mathbf{{x}^{b}}\mathbf{{x}^{c}})$.
Since $\pi$ is an algebra homomorphism, we have
$\tilde{\mathbf{x}}^{\mathbf{a}}\tilde{\mathbf{x}}^{\mathbf{b}}=\mathbf{\widetilde{{x}^{a}x^{b}}}+y$
and
$\tilde{\mathbf{x}}^{\mathbf{b}}\tilde{\mathbf{x}}^{\mathbf{c}}=\mathbf{\widetilde{{x}^{b}x^{c}}}+z$
for $y,z\in \Ker\pi$. Therefore,
\begin{eqnarray*}
\hat{\xi}_{\alpha_{i}}(\mathbf{{x}^{a}}\mathbf{{x}^{b}},\mathbf{{x}^{c}})&=&
\tilde{\xi}_{\alpha_{i}}(\mathbf{\widetilde{{x}^{a}x^{b}}},\tilde{\mathbf{x}}^{\mathbf{c}})\\
&=&\tilde{\xi}_{\alpha_{i}}(\tilde{\mathbf{x}}^{\mathbf{a}}\tilde{\mathbf{x}}^{\mathbf{b}}-y,\tilde{\mathbf{x}}^{\mathbf{c}})=
\tilde{\xi}_{\alpha_{i}}(\tilde{\mathbf{x}}^{\mathbf{a}}\tilde{\mathbf{x}}^{\mathbf{b}},\tilde{\mathbf{x}}^{\mathbf{c}})\\
&=&\tilde{\xi}_{\alpha_{i}}(\tilde{\mathbf{x}}^{\mathbf{a}},\tilde{\mathbf{x}}^{\mathbf{b}}\tilde{\mathbf{x}}^{\mathbf{c}})\\
&=&\tilde{\xi}_{\alpha_{i}}(\tilde{\mathbf{x}}^{\mathbf{a}},\mathbf{\widetilde{{x}^{b}x^{c}}}+z)=
\tilde{\xi}_{\alpha_{i}}(\tilde{\mathbf{x}}^{\mathbf{a}},\mathbf{\widetilde{{x}^{b}x^{c}}})\\
&=&\hat{\xi}_{\alpha_{i}}(\mathbf{{x}^{a}},\mathbf{{x}^{b}}\mathbf{{x}^{c}}).
\end{eqnarray*}
Next, let us show that they are linear independent in
$\H^{2}(|R_{\mathfrak{g}}|,\kappa)$. It is equivalent to show that
for any linear combination
$f=\sum_{i=1}^{N}c_{i}\hat{\xi}_{\alpha_{i}}$, if it is a coboundary
then every $c_{i}=0$. Assume that $f=\partial h$ for some
$h:|R_{\mathfrak{g}}|^{+}\to \kappa$. Then
$$c_{i}=f(x_{\alpha_{i}},x_{\alpha_{i}}^{N_{i}-1})=\partial h(x_{\alpha_{i}},x_{\alpha_{i}}^{N_{i}-1})
=-h(x_{\alpha_{i}}^{N_{i}})=0$$ since $x_{\alpha_{i}}^{N_{i}}=0$ in
$|R_{\mathfrak{g}}|$ by Lemma 2.10.
\end{proof}

See Section 6 in \cite{MW} for the definitions of such elements in
the case of pointed Hopf algebras. We are now in the position to
prove the following theorem.

\begin{theorem} The algebra
$\H^{\ast}(|\Gr^{1}(\mathbf{u}(\mathfrak{g}))|,\kappa)$ is finitely
generated. If $M$ is a finitely generated
$|\Gr^{1}(\mathbf{u}(\mathfrak{g}))|$-module, then
$\H^{\ast}(|\Gr^{1}(\mathbf{u}(\mathfrak{g}))|,M)$ is a finitely
generated module over
$\H^{\ast}(|\Gr^{1}(\mathbf{u}(\mathfrak{g}))|,\kappa)$.
\end{theorem}
\begin{proof} By Lemma 2.12, there is a filtration on
$|R_{\mathfrak{g}}|$ and results a graded algebra
$|\Gr^{2}(R_{\mathfrak{g}})|$. As the filtration is finite, there is
a convergent spectral sequence associated to the filtration by 5.4.1
in \cite{We}:
\begin{equation} E_{1}^{s,t}=\H^{s+t}(\Gr^{2}_{(s)}(|R_{\mathfrak{g}}|),\kappa)
\Rightarrow \H^{s+t}(|R_{\mathfrak{g}}|,\kappa).
\end{equation}
Since the PBW basis elements (2.4) are eigenvectors for
$\mathbf{u}(\mathfrak{h})$, the action of $\mathbf{u}(\mathfrak{h})$
on $|R_{\mathfrak{g}}|$ preserves the filtration and we further get
a spectral sequence converging to the cohomology of
$|R_{\mathfrak{g}}\#
\mathbf{u}(\mathfrak{h})|=|\Gr^{1}(\mathbf{u}(\mathfrak{g}))|$:
\begin{equation} \H^{s+t}(\Gr^{2}_{(s)}(|R_{\mathfrak{g}}|),\kappa)^{\mathbf{u}(\mathfrak{h})}
\Rightarrow
\H^{s+t}(|R_{\mathfrak{g}}|,\kappa)^{\mathbf{u}(\mathfrak{h})}\cong
\H^{s+t}(|\Gr^{1}(\mathbf{u}(\mathfrak{g}))|,\kappa),
\end{equation}
where the isomorphism $``\cong"$ can be proved similarly just as in
the proof of Proposition 3.4. We can replace $\kappa$ by $M$ in
(3.6), (3.7) to get convergent spectral sequences with  coefficients
in $M$.

By Proposition 3.6, we have some elements $\hat{\xi}_{\alpha_{i}}$
in $\H^{2}(|R_{\mathfrak{g}}|,\kappa)$. We wish to related the
functions $\hat{\xi}_{\alpha_{i}}$ to elements on the $E_{1}$-page
of the spectral sequence (3.6). In fact, one can copy the arguments
stating before Lemma 5.1 in \cite{MPSW} and can assume that
$\hat{\xi}_{\alpha_{i}}\in E_{1}^{c,2-c}\cong
\H^{2}(|\Gr^{2}(R_{\mathfrak{g}})|,\kappa)$ for some $c\in
\mathbb{N}$. Since $\hat{\xi}_{\alpha_{i}}\in
\H^{2}(|R_{\mathfrak{g}}|,\kappa)$, they are permanent cycles. Now,
by Proposition 3.3, $\H^{2}(|\Gr^{2}(R_{\mathfrak{g}})|,\kappa)$ is
indeed spanned by $\xi_{i}$ for $1\leq i\leq
N=\Phi^{\#}-\Phi_{12}^{\#}$.\\[1.2mm]
\textbf{Claim 1.} \emph{In
$\H^{2}(|\Gr^{2}(R_{\mathfrak{g}})|,\kappa)$,
$\xi_{i}=\hat{\xi}_{\alpha_{i}}$.} (The proof of this claim is the
same with that of Lemma 5.1 in \cite{MPSW} and thus is omitted.)

Let $B^{\ast,\ast}$ be the bigraded subalgebra of
$E_{1}^{\ast,\ast}$ generated by the elements $\xi_{i}$. By the
claim 1, $B^{\ast,\ast}$ consists of permanent cycles. Let
$A^{\ast,\ast}$ be the subalgebra of $B^{\ast,\ast}$ generated by
$\xi_{i}^{p}$ where $p$ is the characteristic of $\kappa$. By (3.3)
and (3.4) in Proposition 3.4, $\xi_{i}^{p}$ is invariant under the
action of $\mathbf{u}(\mathfrak{h})\otimes \kappa\mathbb{Z}_{2}$.
Therefore, $A^{\ast,\ast}$ is a subalgebra of
$\H^{\ast}(\Gr^{2}(\mathbf{u}(\mathfrak{g})),\kappa)$. Lemma 2.1
implies that $A^{\ast,\ast}$ is commutative since it is concentrated
in even (total) degrees.\\[1.2mm]
\textbf{Claim 2.} \emph{$A^{\ast,\ast}$ satisfies the hypotheses of
Lemma \emph{2.13}.} To show it, it is enough to show that
$E_{1}^{\ast,\ast}$ is a finitely generated module over
$A^{\ast,\ast}$. Proposition 3.3 implies $E_{1}^{\ast,\ast}\cong
\H^{\ast}(|\Gr^{2}(R_{\mathfrak{g}})|,\kappa)$ is generated by
$\xi_{i}$ and $\eta_{i}$ where $\eta_{i}^{2}=0$. Hence
$E_{1}^{\ast,\ast}$ is a finitely generated module over
$B^{\ast,\ast}$ which is clearly a finitely generated module over
$A^{\ast,\ast}$. Therefore, the claim is proved.

Thus Lemma 2.13 (1) is applied and so
$\H^{\ast}(|R_{\mathfrak{g}}|,\kappa)$ is a Noetherian
Tot($A^{\ast,\ast}$)-module. Moreover, the action of
$\mathbf{u}(\mathfrak{h})$ on $\H^{\ast}(|R_{\mathfrak{g}}|,\kappa)$
is compatible with the action on $A^{\ast,\ast}$, since the spectral
sequence (3.6) is compatible with the action of
$\mathbf{u}(\mathfrak{h})$. Therefore,
$\H^{\ast}(|\Gr^{1}(\mathbf{u}(\mathfrak{g}))|,\kappa)\cong
\H^{\ast}(|R_{\mathfrak{g}}|,\kappa)^{\mathbf{u}(\mathfrak{h})}$ is
a Noetherian Tot($A^{\ast,\ast}$)-module. Now, Tot($A^{\ast,\ast}$)
is finitely generated since $A^{\ast,\ast}$ is just the polynomial
algebra generated by $\xi_{i}^{p}$. We conclude that
$\H^{\ast}(|\Gr^{1}(\mathbf{u}(\mathfrak{g}))|,\kappa)$ is finitely
generated.

The second statement of the this theorem follows by a direct
application of Lemma 2.13 (2).
\end{proof}

Next result is a direct consequence of Theorem 3.7 and Lemma 2.7.
\begin{corollary}The algebra
$\H^{\ast}(\Gr^{1}(\mathbf{u}(\mathfrak{g})),\kappa)$ is finitely
generated. If $M$ is a finitely generated
$\Gr^{1}(\mathbf{u}(\mathfrak{g}))$-supermodule, then
$\H^{\ast}(\Gr^{1}(\mathbf{u}(\mathfrak{g})),M)$ is a finitely
generated module over
$\H^{\ast}(\Gr^{1}(\mathbf{u}(\mathfrak{g})),\kappa)$.
\end{corollary}

\subsection{Cohomology of $\mathbf{u}(\mathfrak{g})$}

In this subsection, we will give the proof of Theorem 3.1 provided
$\mathfrak{g}\neq \mathbf{A}(1,1)$. Similar to Subsection 3.2, we
have convergent spectral sequences associated the first kind of
filtration given before Lemma 2.11:
\begin{equation} E_{1}^{s,t}=\H^{s+t}(\Gr^{1}_{(s)}(|\mathbf{u}(\mathfrak{g})|),\kappa)
\Rightarrow \H^{s+t}(|\mathbf{u}(\mathfrak{g})|,\kappa),
\end{equation}
\begin{equation}\H^{s+t}(\Gr^{1}_{(s)}(|\mathbf{u}(\mathfrak{g})|),M)
\Rightarrow \H^{s+t}(|\mathbf{u}(\mathfrak{g})|,M),
\end{equation}
for any $|\mathbf{u}(\mathfrak{g})|$-module $M$.

 Previously, we
identify the element $\xi_{i}\in
\H^{2}(|\Gr^{2}(R_{\mathfrak{g}})|,\kappa)$ with the element
$\hat{\xi}_{\alpha_{i}}\in \H^{2}(|R_{\mathfrak{g}}|,\kappa)$. From
this, we know that $\xi_{i}$ is a permanent cycle and
$\H^{\ast}(|R_{\mathfrak{g}}\#
\mathbf{u}(\mathfrak{h})|,\kappa)=\H^{\ast}(|\Gr^{1}(\mathbf{u}(\mathfrak{g}))|,\kappa)$
is finitely generated over the subalgebra generated by all
$\hat{\xi}_{\alpha_{i}}^{p}$. So our next aim is to find an element
$f_{\alpha_{i}}\in \H^{\ast}(|\mathbf{u}(\mathfrak{g})|,\kappa)$
which can be identified with $\hat{\xi}_{\alpha_{i}}^{p}$. If so,
$\hat{\xi}_{\alpha_{i}}^{p}$ will be permanent cycles and Lemma 2.13
can be applied.

For each $i\in\{1,\ldots,\Phi^{\#}-\Phi_{12}^{\#}\}$, let
$\alpha_{i}$ be the corresponding root. For our purpose, we choose a
PBW basis of $U(\mathfrak{g})$, described as in Proposition 2.5 (1),
with requirements: $s_{1}=\Phi_{11}^{\#}, s=\Phi_{1}^{\#}$,
$x_{i}=x_{\alpha_{i}}$ for $1\leq i\leq s$ and
$y_{j}=x_{\alpha_{s+j}}$ for $1\leq j\leq
\Phi_{0}^{\#}-\Phi_{12}^{\#}$ where $x_{\alpha_{k}}$ is defined
before Lemma 2.11. Roughly speaking, we just want the PBW basis
elements given in Lemma 3.5 are still PBW basis elements in the
following discussions. We choose a PBW basis for
$\mathbf{u}(\mathfrak{g})$ with the same requirements as
$U(\mathfrak{g})$. Such PBW basis will be fixed from now on until
the end of this subsection.

 Define a $\kappa$-linear function
$\tilde{f}_{\alpha_{i}}:(|U(\mathfrak{g})|^{+})^{2p}\to \kappa$  as
follows. Let $r_{1},\ldots,r_{2p}$ be PBW basis elements. If all of
them have no factors belonging to $U(\mathfrak{h})$, then
$$\tilde{f}_{\alpha_{i}}(r_{1}\otimes \cdots r_{2p}):=c_{12}c_{34}\cdots c_{2p-1,2p}$$
where $c_{ij}$ is the coefficient of $x_{\alpha_{i}}^{N_{i}}$ in the
product $r_{i}r_{j}$ as a linear combination of PBW basis elements.
And set $\tilde{f}_{\alpha_{i}}$ to be zero whenever there is a
$r_{i}$ which contains a factor living in $U(\mathfrak{h})$.

Similar to Subsection 3.2, we will show that
$\tilde{f}_{\alpha_{i}}$ factors through the quotient
$\pi:U(\mathfrak{g})\to \mathbf{u}(\mathfrak{g})$ to give a map
$(|\mathbf{u}(\mathfrak{g})|^{+})^{2p}\to \kappa$. Note that by the
definition of $\tilde{f}_{\alpha_{i}}$, it is always 0 whenever the
elements of $U(\mathfrak{h})$ appear in a PBW basis element. So we
need only to consider the PBW basis elements totally the same as
that of $\tilde{R}_{\mathfrak{g}}$. So we can apply the same
arguments designed for $\tilde{\xi}_{\alpha_{i}}$ to
$\tilde{f}_{\alpha_{i}}$ and show that it indeed factors through the
quotient map $\pi:U(\mathfrak{g})\to \mathbf{u}(\mathfrak{g})$.
Also, we choose a section $\widetilde{-}:
\mathbf{u}(\mathfrak{g})\to U(\mathfrak{g})$ of the quotient map
$\pi$. Then define $f_{\alpha_{i}}:
(|\mathbf{u}(\mathfrak{g})|^{+})^{2p}\to \kappa$ by setting
$$f_{\alpha_{i}}(r_{1}\otimes \cdots \otimes r_{2p}):=\tilde{f}_{\alpha_{i}}(\tilde{r}_{1}\otimes \cdots \tilde{r}_{2p})$$
for PBW basis elements $r_{1},\ldots,r_{2p}\in
\mathbf{u}(\mathfrak{g})$.

\begin{proposition} The set
$\{f_{\alpha_{i}}|i=1,\ldots,\Phi^{\#}-\Phi_{12}^{\#}\}$ represents
a linear independent subset of
$\H^{2p}(|\mathbf{u}(\mathfrak{g})|,\kappa)$.
\end{proposition}
\begin{proof} The proof is similar to that of Lemma 6.2 in
\cite{MPSW} and Proposition 3.6. For completeness, we write it out.

Firstly, we show that $\tilde{f}_{\alpha_{i}}$ is a $2p$-cocycle on
$|U(\mathfrak{g})|$. Let $r_{0},\ldots,r_{2p}\in
|U(\mathfrak{g})|^{+}$ be PBW basis elements without factors coming
from $U(\mathfrak{h})$. Then
$$\partial(\tilde{f}_{\alpha_{i}})(r_{0}\otimes \cdots \otimes r_{2p})=\sum_{j=0}^{2p-1}(-1)^{i+1}\tilde{f}_{\alpha_{i}}
(r_{0}\otimes \cdots  \otimes r_{j}r_{j+1}\otimes \cdots\otimes
r_{2p}).$$ By the definition of $\tilde{f}_{\alpha_{i}}$, the first
two terms cancel and similarly for all other terms. So
$\partial(\tilde{f}_{\alpha_{i}})=0$.

Now we verify that $f_{\alpha_{i}}$ is a $2p$-cocycle. Also, let
$r_{0},\ldots,r_{2p}\in |\mathbf{u}(\mathfrak{g})|^{+}$ be PBW basis
elements. Then
$$\partial({f}_{\alpha_{i}})(r_{0}\otimes \cdots \otimes r_{2p})=\sum_{j=0}^{2p-1}(-1)^{i+1}{f}_{\alpha_{i}}
(r_{0}\otimes \cdots  \otimes r_{j}r_{j+1}\otimes \cdots\otimes
r_{2p}).$$ Using the same methods as in the proof of Proposition
3.6, we have
\begin{eqnarray*}
{f}_{\alpha_{i}}(r_{0}r_{1}\otimes r_{2}\otimes \cdots  \otimes
r_{2p})&=& \tilde{f}_{\alpha_{i}}(\widetilde{r_{0}r_{1}}\otimes
\tilde{r}_{2}\otimes \cdots  \otimes \tilde{r}_{2p})\\
&=&\tilde{f}_{\alpha_{i}}(\tilde{r}_{0}\tilde{r}_{1}\otimes
\tilde{r}_{2}\otimes \cdots  \otimes
\tilde{r}_{2p})\\
&=&\tilde{f}_{\alpha_{i}}(\tilde{r}_{0}\otimes
\tilde{r}_{1}\tilde{r}_{2}\otimes \cdots  \otimes \tilde{r}_{2p})\\
&=&\tilde{f}_{\alpha_{i}}(\tilde{r}_{0}\otimes
\widetilde{r_{1}{r}_{2}}\otimes \cdots  \otimes \tilde{r}_{2p})\\
&=&{f}_{\alpha_{i}}(r_{0}\otimes r_{1}r_{2}\otimes \cdots \otimes
r_{2p}).
\end{eqnarray*}
Similarly, we have $${f}_{\alpha_{i}} (r_{0}\otimes \cdots  \otimes
r_{j}r_{j+1}\otimes \cdots\otimes r_{2p})={f}_{\alpha_{i}}
(r_{0}\otimes \cdots  \otimes r_{j+1}r_{j+2}\otimes \cdots\otimes
r_{2p})$$ for $j=0,\ldots,2p-2$. So $\partial({f}_{\alpha_{i}})=0$.

Now assume that $\sum_{i}c_{i}{f}_{\alpha_{i}}=\partial h$ for some
$h\in \Hom_{\kappa}((|\mathbf{u}(\mathfrak{g})|^{+})^{\otimes 2p-1},
\kappa)$. Then for each $i$,
\begin{eqnarray*}
c_{i}&=&(\sum_{j}c_{j}{f}_{\alpha_{j}})(x_{\alpha_{i}}\otimes
x_{\alpha_{i}}^{N_{i}-1}\otimes \cdots \otimes x_{\alpha_{i}}\otimes
x_{\alpha_{i}}^{N_{i}-1})\\
&=& (\partial h)(x_{\alpha_{i}}\otimes
x_{\alpha_{i}}^{N_{i}-1}\otimes \cdots \otimes x_{\alpha_{i}}\otimes
x_{\alpha_{i}}^{N_{i}-1})\\
&=&\sum \pm h(x_{\alpha_{i}}\otimes x_{\alpha_{i}}^{N_{i}-1}\otimes
\cdots \otimes x_{\alpha_{i}}^{N_{i}}\otimes \cdots \otimes
x_{\alpha_{i}}\otimes x_{\alpha_{i}}^{N_{i}-1})\\
&=&0
\end{eqnarray*}
since $x_{\alpha_{i}}^{N_{i}}=0$ in $\mathbf{u}(\mathfrak{g})$ by
Lemma 2.10.

\end{proof}

\noindent{\bf Proof of Theorem 3.1 in case $\mathfrak{g}\neq
\mathbf{A}(1,1)$.} The functions ${f}_{\alpha_{i}}$ correspond to
their counterpart $\hat{\xi}_{\alpha_{i}}^{p}$ defined on
$|\Gr^{1}(\mathbf{u}(\mathfrak{g}))|$, in the $E_{1}$-page of the
spectral sequence (3.8), by observing that they are the same
functions at the level of chain complex (2.3) where we need replace
$|\mathbf{u}(\mathfrak{g})|^{+}$ by
$|\Gr^{1}(\mathbf{u}(\mathfrak{g}))|^{+}$. Thus Proposition 3.9
implies that the function $\hat{\xi}_{\alpha_{i}}^{p}$ is a
permanent cycle. Now we have known that $E_{1}^{\ast,\ast}\cong
\H^{\ast}(|\Gr^{1}(\mathbf{u}(\mathfrak{g}))|,\kappa)$ is finitely
generated over the subalgebra $A^{\ast,\ast}$ generated by all
$\hat{\xi}_{\alpha_{i}}^{p}$ (see the proof of Theorem 3.7). Thus
$A^{\ast,\ast}$ satisfies the conditions of Lemma 2.13 and thus
$\H^{\ast}(|\mathbf{u}(\mathfrak{g})|,\kappa)$ is a Noetherian
Tot$(A^{\ast,\ast})$-module and thus finitely generated. By Lemma
2.7, the first part of Theorem 3.1 is proved. The second part can be
prove similarly by applying Lemma 2.13 (2) and Lemma 2.7.

\subsection{The case $\mathfrak{g}=\mathbf{A}(1,1)$}
We deal with the case $\mathfrak{g}=\mathbf{A}(1,1)$ in a bigger
context: Those basic classical Lie superalgebras with $\Phi_{12}$
being empty. By the descriptions of root supersystems given in
Section 2.5.4 in \cite{Kac1}, this includes all basic classical Lie
superalgebras except $\mathbf{B}(m,n)$ and $\mathbf{G}(3)$. For such
Lie superalgebras, we have a nice filtration on them.

We give a notion at first. For a coalgebra $C$ and $D\subseteq C$ a
subcoalgebra of $C$, define
$$\wedge^{0}D:=D,\;\;\wedge^{1} D:=\Delta^{-1}(C\otimes D+D\otimes C),$$
$$\wedge^{i}D:=\wedge^{1}(\wedge^{i-1} D)=\Delta^{-1}(C\otimes \wedge^{i-1} D
+\wedge^{i-1} D\otimes C)$$ for $i\geq 2$. If $D$ contains the
coradical $C_{0}$ of $C$, by definition $C_{0}$ is the sum of all
simple subcoalgebras of $C$, then $D\subseteq \wedge D\subseteq
\wedge^{2} D\subseteq \cdots$ will give a filtration of $C$. See
Chapter 5 in \cite{Mon} for details.

Let $\mathfrak{g}$ be a basic classical Lie superalgebra with
$\Phi_{12}=\phi$. Then $\mathbf{u}(\mathfrak{g})$ is a finite
dimensional super cocommutative Hopf algebra and its coradical is
$\kappa$. Define
$$F^{i}\mathbf{u}(\mathfrak{g}):=\wedge^{i}\mathbf{u}(\mathfrak{h})$$
for $i\geq 0$ and then this gives a filtration of
$\mathbf{u}(\mathfrak{g})$. The associated graded algebra is denoted
by $\gr(\mathbf{u}(\mathfrak{g}))$. It is a superalgebra naturally.
For any $\alpha\in \Phi$, we fix a basis $b_{\alpha}$ of
$\mathfrak{g}_{\alpha}$. By taking the union of such $b_{\alpha}$,
we get a basis of $\bigoplus_{\alpha\in\Phi}\mathfrak{g}_{\alpha}$.
Denote this basis by $\{x_{1},\ldots,x_{m},x_{m+1},\ldots,x_{m+n}\}$
and assume that $x_{i}\in
\bigoplus_{\alpha\in\Phi_{1}}\mathfrak{g}_{\alpha}$ for $1\leq i\leq
m$ while $x_{i}\not\in
\bigoplus_{\alpha\in\Phi_{1}}\mathfrak{g}_{\alpha}$ for $m< i\leq
m+n$ .

\begin{lemma} $\gr (\mathbf{u}(\mathfrak{g}))\cong S_{\mathfrak{g}}\#
\mathbf{u}(\mathfrak{h})$ where $S_{\mathfrak{g}}$ is generated by
$x_{1},\ldots,x_{m+n}$ with relations
\begin{equation}x_{i}x_{j}=\left \{
\begin{array}{ll} -x_{j}x_{i}, & \;1\leq i<j\leq m\\
x_{j}x_{i}, &\;1\leq i<j,\;j>m,
\end{array}\right.\;\;\;\;x_{i}^{n_{i}}=0,\end{equation}
where $n_{i}=\left \{
\begin{array}{ll} 2, & \;1\leq i\leq m\\
p, & \;m<i\leq m+n.
\end{array}\right.$
\end{lemma}
\begin{proof} Here the action of $\mathbf{u}(\mathfrak{h})$ on $S_{\mathfrak{g}}$ is gotten through
extending the actions of $\mathfrak{h}$ on
$\bigoplus_{\alpha\in\Phi}\mathfrak{g}_{\alpha}$ naturally. By the
definition of the coproduct of $\mathbf{u}(\mathfrak{g})$,
$\bigoplus_{\alpha\in\Phi}\mathfrak{g}_{\alpha}\subset
\wedge^{1}\mathbf{u}(\mathfrak{h})$. So $[x_{i},x_{j}]\in
\wedge^{1}\mathbf{u}(\mathfrak{h})$. This implies we have
$$[x_{i},x_{j}]=0$$
in $\gr(\mathbf{u}(\mathfrak{g}))$. It is direct to show that every
$x_{i}^{n_{i}}$ is still a primitive element and so
$x_{i}^{n_{i}}\in \wedge^{1}\mathbf{u}(\mathfrak{h})$ too.
Therefore, $x_{i}^{n_{i}}=0$ in  $\gr(\mathbf{u}(\mathfrak{g}))$.
Now all relations in (3.10) are fulfilled. By comparing the
dimensions, we indeed get the desire isomorphism.
\end{proof}

\noindent{\bf Proof of Theorem 3.1 in case $\Phi_{12}=\phi$.} Since
$\mathbf{u}(\mathfrak{g})$ is finite dimensional, then the
filtration $F^{0}\mathbf{u}(\mathfrak{g})\subset
F^{1}\mathbf{u}(\mathfrak{g})\subset \cdots$ is finite, that is,
there is $n\in \mathbb{N}$ such that
$F^{n}\mathbf{u}(\mathfrak{g})=\mathbf{u}(\mathfrak{g})$. So have a
convergent spectral sequence
\begin{equation} E_{1}^{s,t}=\H^{s+t}(\gr_{(s)}(|\mathbf{u}(\mathfrak{g})|),\kappa)
\Rightarrow \H^{s+t}(|\mathbf{u}(\mathfrak{g})|,\kappa).
\end{equation}
By Lemma 3.10, $|\gr \mathbf{u}(\mathfrak{g})|\cong
|S_{\mathfrak{g}}\#\mathbf{u}(\mathfrak{h})|$. Now it is clear
$|S_{\mathfrak{g}}|$ is a quantum complete intersection algebra (see
the second paragraph of Subsection 3.1). Thus it cohomology algebra
is clear by Lemma 3.2. Actually, similar to Proposition 3.3, we have
$$\H^{\ast}(|S_{\mathfrak{g}}|,\kappa)\cong k[\xi_{1},\ldots,\xi_{m+n}]\otimes \wedge(m|n)$$
with $m=\Phi_{1}^{\#}, n=\Phi_{0}^{\#}$ and $\deg \xi_{i}=2, \deg
\eta_{i}=1$. Also one can get that $\xi_{i}^{p}\in
\H^{\ast}(|S_{\mathfrak{g}}|,\kappa)^{\mathbf{u}(\mathfrak{h})}\cong\H^{\ast}(|\gr(\mathbf{u}{(\mathfrak{g}}))|,\kappa)$.
By applying the same discussions used in the proof of Claim 2 in
that of Theorem 3.7, $E_{1}^{\ast,\ast}$ is finitely generated over
the subalgebra generated by all $\xi_{i}^{p}$. So Lemma 2.13 (1) can
be applied if we can show all $\xi_{i}^{p}$ are permanent cycles. In
fact, we can define $f_{i}\in
\H^{2p}(|\mathbf{u}(\mathfrak{g})|,\kappa)$ through the same way as
that of $f_{\alpha_{i}}$ (see Proposition 3.9) and get $f_{i}$
corresponds to its counterpart $\xi_{i}^{p}$ defined on
$|\gr(\mathbf{u}(\mathfrak{g}))|$. Therefore, every $\xi_{i}^{p}$ is
a permanent cycle and thus
$\H^{\ast}(|\mathbf{u}(\mathfrak{g})|,\kappa)$ is a finitely
generated algebra. Using Lemma 2.7, we know that
$\H^{\ast}(\mathbf{u}(\mathfrak{g}),\kappa)$ is also finitely
generated as an algebra.

The second part of the theorem can be proved by applying Lemma 2.13
(2) and Lemma 2.7 now.

\begin{remark} \emph{To show the theorem, we cannot apply the filtration developed
in this subsection to Lie superalgebras $\mathbf{B}(m,n),
\mathbf{G}(3)$ directly since otherwise more nilpotent elements will
be created. On the contrary, the two kinds of filtration given in
Section $2$ can be applied to $\mathbf{A}(1,1)$ and indeed
$\Gr^{2}(\mathbf{A}(1,1))=\gr(\mathbf{A}(1,1))$. But in the case of
$\mathfrak{g}=\mathbf{A}(1,1)$, it is possible that
$\dim_{\kappa}\mathfrak{g}_{\alpha}\geq 2$ and so the notation
$x_{\alpha}$ has no meaning now. Therefore, if we want to deal with
all basic classical Lie superalgebras in a unified way (that is, by
using two kinds of filtration ), the notations and descriptions will
be too delicate to grasp the main line.}
\end{remark}

\section{Support varieties and representation type of Lie superalgebras}

In this section, we will recall the definition of the support
variety of a module and give its relation with the complexity of
this module. As a consequence, we will show that only
$|\mathbf{u}(\mathfrak{sl}_{2})|$,
$|\mathbf{u}(\mathfrak{osp}(1|2))|$ are tame and the others
$|\mathbf{u}(\mathfrak{g})|$ are all wild (see Section 5 for
explicit description of $\mathfrak{osp}(1|2)$).

Let $\mathfrak{g}$ be a basic classical Lie superalgebra and $N$ a
finitely generated left $\mathbf{u}(\mathfrak{g})$-supermodule. By
Corollary 2.8 and Theorem 3.1,
$\H^{ev}(\mathbf{u}(\mathfrak{g}),\kappa)$ is a finitely generated
commutative algebra and $\H^{\ast}(\mathbf{u}(\mathfrak{g}),N)$ is a
finitely generated
$\H^{ev}(\mathbf{u}(\mathfrak{g}),\kappa)$-module. In particular,
for any finitely generated $\mathbf{u}(\mathfrak{g})$-supermodule
$M$, $\Ext^{\ast}_{\mathbf{u}(\mathfrak{g})}(M,M):=\bigoplus_{i\geq
0} \H^{i}_{\mathbf{u}(\mathfrak{g})}(M,M)\cong \bigoplus_{i\geq
0}\H^{i}(\mathbf{u}(\mathfrak{g}), M^{\ast}\otimes M)$ is finitely
generated over $\H^{ev}(\mathbf{u}(\mathfrak{g}),\kappa)$ where
$M^{\ast}$ is the dual $\mathbf{u}(\mathfrak{g})$-module of $M$. Let
$I_{M}$ be the annihilator of action of
$\H^{ev}(\mathbf{u}(\mathfrak{g}),\kappa)$ on
$\Ext^{\ast}_{\mathbf{u}(\mathfrak{g})}(M,M)$. The
\emph{cohomological support variety} of $M$ is defined to be
$$\mathcal{V}_{\mathbf{u}(\mathfrak{g})}(M):=Z(I_{M})\subset \textrm{Maxspec}(\H^{ev}(\mathbf{u}(\mathfrak{g}),\kappa)).$$
Note that we can regard $M$ as a
$\mathbf{u}(\mathfrak{g})\#\kappa\mathbb{Z}_{2}$-module  by Lemma
2.2.

Let $A$ be an associative algebra, $M$ an $A$-module with minimal
projective resolution
$$\cdots \rightarrow P_{n}\rightarrow P_{n-1}\rightarrow \cdots \rightarrow P_{0}\rightarrow M\rightarrow 0.$$
Then the \emph{complexity} of $M$ is defined to be the integer
$$\C_{A}(M):=\textrm{min}\{c\in \mathbb{N}_{0}\cup \infty\;| \exists \lambda>0:\textrm{dim}_{k}P_{n}\leq
\lambda n^{c-1},\;\forall \;n\geq 1\}.$$

\begin{lemma}Let $\mathfrak{g}$ be a basic classical Lie superalgebra and $M$ a
finitely generated left $\mathbf{u}(\mathfrak{g})$-supermodule. Then
$$\dim \mathcal{V}_{\mathbf{u}(\mathfrak{g})}(M)=\emph{\C}_{\mathbf{u}(\mathfrak{g})\#\kappa\mathbb{Z}_{2}}(M).$$
\end{lemma}
\begin{proof} By definition, $\H^{ev}(\mathbf{u}(\mathfrak{g}),\kappa)=\bigoplus_{i\geq 0}
\Ext^{2i}_{\mathbf{u}(\mathfrak{g})\#\kappa\mathbb{Z}_{2}}(\kappa,\kappa)$
and now $\mathbf{u}(\mathfrak{g})\#\kappa\mathbb{Z}_{2}$ is an
ordinary finite dimensional Hopf algebra. So this lemma is just a
corollary of Proposition 2.3 in \cite{FW}.
\end{proof}

Recall the finite dimensional associative algebras over an
algebraically closed field $\kappa$ can be divided into three
classes (see \cite{Dr}): A finite-dimensional algebra $A$ is said to
be of \emph{finite representation type} provided there are finitely
many non-isomorphic indecomposable $A$-modules. $A$ is of \emph{tame
type} or $A$ is a \emph{tame} algebra if $A$ is not of finite
representation type, whereas for any dimension $d>0$, there are
finite number of $A$-$\kappa[T]$-bimodules $M_{i}$ which are free of
finite rank as right $\kappa[T]$-modules such that all but a finite
number of indecomposable $A$-modules of dimension $d$ are isomorphic
to $M_{i}\otimes_{\kappa[T]}\kappa[T]/(T-\lambda)$ for $\lambda\in
k$. We say that $A$ is of \emph{wild type} or $A$ is a \emph{wild}
algebra if there is a finitely generated $A$-$\kappa\langle
X,Y\rangle $-bimodule $B$ which is free as a right $\kappa\langle
X,Y\rangle $-module such that the functor $B\otimes_{\kappa\langle
X,Y\rangle}-\;\;$ from $\kappa\langle X,Y\rangle$-mod, the category
of finitely generated $\kappa\langle X,Y\rangle$-modules, to
$A$-mod, the category of finitely generated $A$-modules, preserves
indecomposability and reflects isomorphisms.

The following result is a direct consequence of Proposition 3.2 in
Chapter VI of \cite{ARS}.

\begin{lemma} Let $A$ be a superalgebra and assume that characteristic of $\kappa$ is not $2$. Then $|A|$ and
$A\# \kappa\mathbb{Z}_{2}$ have the same representation type.
\end{lemma}

\begin{remark} For a finite dimensional superalgebra $A$, one also
can define its representation type in the super world, that is, in
the category of supermodules with even homomorphisms. By Lemma $4.2$
and Lemma $2.2$, the representation type of $|A|$ as an ordinary
algebra is indeed the same with that of $A$ when we consider it as a
superalgebra. So to consider the representation type of a
superalgebra $A$, it is enough to consider that of its underline
algebra $|A|$.
\end{remark}

The following conclusion is also needed.

\begin{lemma} If there is a finite dimensional $\mathbf{u}(\mathfrak{g})\#
\kappa\mathbb{Z}_{2}$-module $M$ such that
$\emph{\C}_{\mathbf{u}(\mathfrak{g})\# \kappa\mathbb{Z}_{2}}(M)\geq
3$, then $\mathbf{u}(\mathfrak{g})\# \kappa\mathbb{Z}_{2}$ is wild.
\end{lemma}
\begin{proof} Let $H$ be an arbitrary finite dimensional Hopf
algebra such that $\C_{H}(N)\geq 3$ for some $H$-module $N$. Then
Theorem 3.1 in \cite{FW} implies that $H$ is wild provided
$\H^{\ast}(H,\kappa)$ is finitely generated and $\H^{\ast}(H,N')$ is
a Noetherian module over $\H^{\ast}(H,\kappa)$ for any finite
dimensional $H$-module $N'$. So the lemma is proved due to our
Theorem 3.1.
\end{proof}

\begin{theorem} Let $\mathfrak{g}=\mathfrak{g}_{0}\oplus \mathfrak{g}_{1}$ be a basic classical Lie
superalgebra over $\kappa$. Then $|\mathbf{u}(\mathfrak{g})|$ is
wild except $\mathfrak{g}=\mathfrak{sl}_{2}$ or
$\mathfrak{g}=\mathfrak{osp}(1|2)$ or $\mathfrak{g}=\mathbf{C}(2)$.
Both $|\mathbf{u}(\mathfrak{sl}_{2})|$ and
$|\mathbf{u}(\mathfrak{ops}(1|2))|$ are tame.
\end{theorem}
\begin{proof} The proof is base on the estimation of the number $\C_{\mathbf{u}(\mathfrak{g})\#
\kappa\mathbb{Z}_{2}}(\kappa)$. By Proposition 2.1 in \cite{FW}, we
have $$\C_{\mathbf{u}({\mathfrak{g}_{0}})}(\kappa)\leq
\C_{\mathbf{u}(\mathfrak{g})\# \kappa\mathbb{Z}_{2}}(\kappa).$$
Owing to (1.4) in \cite{FB},
$\mathcal{V}_{\mathbf{u}({\mathfrak{g}_{0}})}(\kappa)$ can be
identified with
$$\mathscr{V}_{\mathbf{u}({\mathfrak{g}_{0}})}(\kappa):=\{x\in
\mathfrak{g}_{0}|x^{[p]}=0\}\cup \{0\}.$$ Now we have known that
$\mathfrak{g}_{0}$ is a direct sum of simple Lie algebras of type
$\mathbf{A}_{n},\mathbf{B}_{n},\mathbf{C}_{n},\mathbf{D}_{n},\mathbf{G}_{2}$
or $\kappa$. By Lemma 2.10,
$$\dim \mathcal{V}_{\mathbf{u}({\mathfrak{g}_{0}})}(\kappa)=\dim \mathscr{V}_{\mathbf{u}({\mathfrak{g}_{0}})}(\kappa)\geq 3$$
except $\mathfrak{g}_{0}=\mathfrak{sl}_{2}$ or
$\mathfrak{g}_{0}=\mathfrak{sl}_{2}\oplus \kappa$. Thus Lemma 4.1
implies that $\C_{\mathbf{u}(\mathfrak{g})\#
\kappa\mathbb{Z}_{2}}(\kappa)\geq
\C_{\mathbf{u}({\mathfrak{g}_{0}})}(\kappa)=\dim
\mathcal{V}_{\mathbf{u}({\mathfrak{g}_{0}})}(\kappa)\geq 3$ unless
$\mathfrak{g}_{0}=\mathfrak{sl}_{2}$ or
$\mathfrak{g}_{0}=\mathfrak{sl}_{2}\oplus \kappa$. The latter only
appear in the case $\mathfrak{g}=\mathbf{C}(2)$. So now it is not
hard to see that in the rest list of basic classical Lie
superalgebras only $\mathfrak{sl}_{2}$ and $\mathfrak{osp}(1|2)$
satisfy its even part is $\mathfrak{sl}_{2}$. By applying Lemma 4.4,
the first part of theorem is proved.

For the second part, it in known that
$\mathbf{u}(\mathfrak{sl}_{2})$ is tame (see for example
\cite{Far1}). The algebra $|\mathbf{u}(\mathfrak{osp}(1|2))|$ is
proved to be a tame algebra by Farnsteiner in the  Example in
Section 4 of \cite{Far}.
\end{proof}

\begin{conjecture} The algebra $|\mathbf{u}(\mathbf{C}(2))|$ is a
wild algebra.
\end{conjecture}

\section{Restricted representations of $\mathfrak{osp}(1|2)$}
Comparing with the case $\mathfrak{sl}_{2}$, we know a little about
the representations of $\mathbf{u}(\mathfrak{osp}(1|2))$. In the
last section of the paper, we want to determine all finite
dimensional representations of $\mathbf{u}(\mathfrak{osp}(1|2))$
inspired that fact that $|\mathbf{u}(\mathfrak{osp}(1|2))|$ is tame.
To do it, the representation theory of $\mathfrak{sl}_{2}$ and the
theory of Frobenius extensions are need. In this section, we only
need $p\neq 2$.

\subsection{$\mathfrak{sl}_{2}$ case} In this subsection, the
restricted simples and projectives of $\mathfrak{sl}_{2}$ are
summarized. Recall the restricted enveloping algebra
$\mathbf{u}(\mathfrak{sl}_{2})$ of $\mathfrak{sl}_{2}$ is generated
by $e,f,h$ with relations
$$[h,e]=2e,\;\;[h,f]=-2f,\;\;[e,f]=h,\;\;h^{p}=h,\;\;e^{p}=f^{p}=0.$$
For any $0\leq \lambda\leq p-1$, we define a $\lambda+1$-dimensional
$\mathbf{u}(\mathfrak{sl}_{2})$-module $V_{0}^{\lambda}$ as follows.
This module has a basis $v_{0},v_{1},\ldots,v_{\lambda}$ and the
actions of the generators are given by the following rules
\begin{equation}
hv_{i}=(\lambda-2i)v_{i},\;\;ev_{i}=-i(\lambda+1-i)v_{i-1},\;\;fv_{i}=-v_{i+1}
\end{equation}
where $i=0,1,\ldots,\lambda$ and $v_{-1}=v_{\lambda+1}=0$. It is
well-known that $\{V_{0}^{\lambda}|0\leq \lambda\leq p-1\}$ forms a
complete non-redundant list of simple
$\mathbf{u}(\mathfrak{sl}_{2})$-modules and $V_{0}^{p-1}$ is
projective, which is called a \emph{Steinberg module }in general.

It is convenient to use a graphical representation for the
structures of modules. Every vertex stands for a vector from our
chosen basis; arrows and dotted ones show the actions of $e$ and $f$
respectively. The example below is for $p=3, \lambda=2$.
\begin{figure}[hbt]
\begin{picture}(50,60)(0,0)

 \put(0,60){\makebox(0,0){$\cdot$}}
 \put(2,55){\vector(0,-1){15}}
 \put(-2,55){\vector(0,1){0}}
 \put(-2,48){\makebox(0,0){$\vdots$}}

 \put(0,35){\makebox(0,0){$\cdot$}}
 \put(2,30){\vector(0,-1){15}}
 \put(-2,30){\vector(0,1){0}}
 \put(-2,23){\makebox(0,0){$\vdots$}}

  \put(0,10){\makebox(0,0){$\cdot$}}
\end{picture}
\end{figure}

Also for any $0\leq \lambda\leq p-2$, we define the module
$P_{0}^{p-2-\lambda}$ by the following rules. The basis of
$P_{0}^{p-2-\lambda}$ is $\{b_{i},a_{i},x_{j},y_{j}|0\leq i\leq
p-2-\lambda, \;0\leq j\leq \lambda\}$ and the actions of $h,e,f$ are
given by:

\begin{eqnarray}
hb_{i}&=&(p-2-\lambda-2i)b_{i},\;\; fb_{i}=-b_{i+1}\\
eb_{i}&=&-i(p-\lambda-1-i)b_{i-1}+a_{i-1};
\end{eqnarray}
\begin{eqnarray}
hy_{j}=(2j-\lambda)y_{j},ey_{j}=y_{j+1},&&hx_{j}=(\lambda-2j)x_{j},fx_{j}=-x_{j+1}\\
fy_{j}=-j(j-\lambda-1)y_{j-1};&&ex_{j}=-j(\lambda+1-j)x_{j-1};
\end{eqnarray}
\begin{equation}
ha_{i}=(p-2-\lambda-2i)a_{i},\;\;ea_{i}=-i(p-\lambda-1-i)a_{i-1},\;\;fa_{i}=-a_{i+1},
\end{equation}
where $b_{p-1-\lambda}=y_{\lambda},
\;y_{\lambda+1}=a_{p-2-\lambda},\;
x_{\lambda+1}=a_{0},\;b_{-1}=x_{\lambda}$.

The graphical description of the $P_{0}^{p-2-\lambda}$ (for
$p=5,\lambda=1$) is indicated as follows.

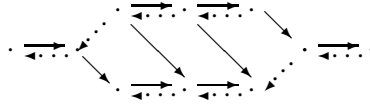
\begin{figure}[hbt]
\begin{picture}(100,40)(0,0)

 \put(0,30){\makebox(0,0){$\cdot$}}
 \put(5,32){\vector(1,0){15}}
 \put(5,28){\vector(-1,0){0}}
  \put(17,28){\makebox(0,0){$\cdots$}}
  \put(5,25){\vector(1,-1){20}}
  \put(30,25){\vector(1,-1){20}}

  \put(25,30){\makebox(0,0){$\cdot$}}
 \put(30,32){\vector(1,0){15}}
 \put(30,28){\vector(-1,0){0}}
  \put(42,28){\makebox(0,0){$\cdots$}}
   \put(50,30){\makebox(0,0){$\cdot$}}

 \put(55,30){\vector(1,-1){10}} \put(70,15){\makebox(0,0){$\cdot$}}
  \put(70,15){\makebox(0,0){$\cdot$}}
 \put(75,17){\vector(1,0){15}}
 \put(75,13){\vector(-1,0){0}}
  \put(87,13){\makebox(0,0){$\cdots$}} \put(95,15){\makebox(0,0){$\cdot$}}

\put(65,10){\makebox(0,0){$\cdot$}}\put(59,4){\makebox(0,0){$\cdot$}}\put(55,0){\vector(-1,-1){0}}
\put(62,7){\makebox(0,0){$\cdot$}}

 \put(0,0){\makebox(0,0){$\cdot$}}
 \put(5,2){\vector(1,0){15}}
 \put(5,-2){\vector(-1,0){0}}
  \put(17,-2){\makebox(0,0){$\cdots$}}

  \put(25,0){\makebox(0,0){$\cdot$}}
 \put(30,2){\vector(1,0){15}}
 \put(30,-2){\vector(-1,0){0}}
  \put(42,-2){\makebox(0,0){$\cdots$}}
   \put(50,0){\makebox(0,0){$\cdot$}}

   \put(-5,25){\makebox(0,0){$\cdot$}}\put(-11,19){\makebox(0,0){$\cdot$}}\put(-15,15){\vector(-1,-1){0}}
\put(-8,22){\makebox(0,0){$\cdot$}}

 \put(-13,13){\vector(1,-1){10}}

 \put(-40,15){\makebox(0,0){$\cdot$}}
 \put(-35,17){\vector(1,0){15}}
 \put(-35,13){\vector(-1,0){0}}
  \put(-23,13){\makebox(0,0){$\cdots$}}
  \put(-15,15){\makebox(0,0){$\cdot$}}
\end{picture}
\end{figure}

It should be known that $\{P_{0}^{p-2-\lambda},V_{0}^{p-1}|0 \leq
\lambda \leq p-2\}$ forms a complete list of indecomposable
projective $\mathbf{u}(\mathfrak{sl}_{2})$-modules up to
isomorphism. For safety, one also can duplicate the proof of Lemma
2.2.6 in \cite{Xiao} to show this fact.

  One also can use the following easy
way to represent the structure of $P_{0}^{p-2-\lambda}$ where we use
$\bullet$ or $\circ$ to denote the composition factors of
$P_{0}^{p-2-\lambda}$.

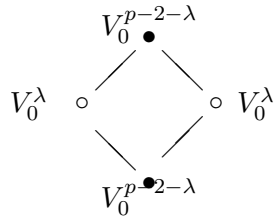
\begin{figure}[hbt]
\begin{picture}(150,70)(0,0)
 \put(30,60){\makebox(0,0){$\bullet$}}  \put(30,65){\makebox(0,0){$V_{0}^{p-2-\lambda}$}}
 \put(35,55){\line(1,-1){15}}
 \put(55,35){\makebox(0,0){$\circ$}}\put(70,35){\makebox(0,0){$V_{0}^{\lambda}$}}
 \put(50,30){\line(-1,-1){15}} \put(5,35){\makebox(0,0){$\circ$}}
\put(-15,35){\makebox(0,0){$V_{0}^{\lambda}$}}
 \put(25,10){\line(-1,1){15}}
 \put(30,5){\makebox(0,0){$\bullet$}}\put(25,55){\line(-1,-1){15}}
\put(30,0){\makebox(0,0){$V_{0}^{p-2-\lambda}$}}

\end{picture}
\end{figure}
From this, it is not hard to see $V_{0}(\lambda)$ and
$V_{0}(p-2-\lambda)$ belongs to the same block ${B}_{0}(\lambda)$
for any $0\leq \lambda\leq p-2$ and there are exactly
$\frac{p+1}{2}$ blocks. Also, one can compute the endomorphism ring
$\End_{\mathbf{u}(\mathfrak{sl}_{2})}(P_{0}^{\lambda}\oplus
P_{0}^{p-2-\lambda})$ out to get the basic algebra of the
$B_{0}(\lambda)$ now. In fact, we will give such computations for
$|\mathbf{u}(\mathfrak{osp}(1|2))|$ and the readers can recover the
block structures of $\mathbf{u}(\mathfrak{sl}_{2})$ from our
computations easily.

\begin{remark} \emph{Since the notions such as $V^{\lambda},P^{\lambda}$,
etc. will be used for $|\mathbf{u}(\mathfrak{osp}(1|2))|$, we add
the subscript $0$ to each notion and get
$V_{0}^{\lambda},P_{0}^{\lambda}$, etc. denoting the corresponding
concepts appeared in classical case,
$\mathbf{u}(\mathfrak{sl}_{2})$}.
\end{remark}

\subsection{Frobenius extensions} Let $R$ be a ring and $S\subseteq
R$ a subring. Suppose that $\alpha$ is an automorphism of $S$. If
$M$ is an $S$-module, we let $_{\alpha}M$ denote the $S$-module with
a new action defined by $s\ast m:=\alpha(s)m$. We say $R$ is an
$\alpha$-\emph{Frobenius extension} of $S$ if

(i) $R$ is a finitely generated projective $S$-module, and

(ii) there exists an isomorphism $\varphi:R\rightarrow
\Hom_{S}(R,\;_{\alpha}S)$ of $(R,S)$-bimodules. More on Frobenius
extensions and their applications can be found in \cite{BF,Far}. For
our purpose, the following serval concepts are needed.

Given an endomorphism $\beta$ of $S$, a \emph{$\beta$-associative
form} from $R$ to $S$ is a biadditive map $\langle,\rangle: R\times
R\rightarrow S$ such that
$$(a)\;\langle sx,y\rangle=s\langle x, y\rangle,\;\;(b) \langle x,ys\rangle
=\langle x,y\rangle \beta(s),\;\;(c) \langle xr,y\rangle=\langle
x,ry\rangle$$ for all $s\in S$ and $r,x,y\in R$.

Let $\langle,\rangle: R\times R\rightarrow S$ be an
$\alpha^{-1}$-associative form. We say two subsets
$\{x_{1},\ldots,x_{n}\}$ and $\{y_{1},\ldots,y_{n}\}$ of $R$ form a
\emph{dual projective pair} relative to $\langle,\rangle$ if
$$r=\sum_{i=1}^{n}y_{i}\alpha(\langle x_{i},r\rangle)=\sum_{i=1}^{n}\langle r, y_{i}\rangle x_{i}\;\;\textrm{for}\;\textrm{all}\;r\in R.$$

Recall Theorem 1.1 in \cite{BF} states that $R$ is an
$\alpha$-Frobenius extension of $S$ if and only if there is an
$\alpha^{-1}$-associative from $\langle,\rangle$ from $R$ to $S$
relative to which a dual projective pair $\{x_{1},\ldots,x_{n}\}$,
$\{y_{1},\ldots,y_{n}\}$ exists. Now let $R:S$ be an
$\alpha$-Frobenius extension and consider two $R$-modules $M,N$.
Then there is a dual projective pair $\{x_{1},\ldots,x_{n}\}$,
$\{y_{1},\ldots,y_{n}\}$. The mapping
$\textrm{Tr}_{[R:S]}:\Hom_{S}(M,\;_{\alpha}N)\rightarrow
\Hom_{R}(M,N)$, which is defined by
$$\textrm{Tr}_{[R:S]}(f)(m)=\sum_{i=1}^{n}y_{i}f(x_{i}m),\;\;\textrm{for }f\in \Hom_{S}(M,\;_{\alpha}N)\; \textrm{and}\; m\in M$$
is usually called the \emph{trace map}.

The following lemma will give us a connection between
$|\mathbf{u}(\mathfrak{osp}(1|2))|$-modules and
$\mathbf{u}(\mathfrak{sl}_{2})$-modules. To describe it, we fix a
notation firstly.  Let $R$ be a ring and $M,\;N$ two $R$-modules. If
$M$ is a direct summand of $N$ as a $R$-module, then we denote it by
$M|N$.

\begin{lemma} Let $M$ be an
$|\mathbf{u}(\mathfrak{osp}(1|2))|$-module, then
$$M||\mathbf{u}(\mathfrak{osp}(1|2))|\otimes _{\mathbf{u}(\mathfrak{sl}_{2})} M.$$
\end{lemma}

\begin{proof} Define $|\mathbf{u}(\mathfrak{osp}(1|2))|\otimes _{\mathbf{u}(\mathfrak{sl}_{2})} M\rightarrow M$
by $a\otimes m\mapsto am$ for $a\in
|\mathbf{u}(\mathfrak{osp}(1|2))|$ and $m\in M$. Clearly, $\varphi$
is an $|\mathbf{u}(\mathfrak{osp}(1|2))|$-epimorphism.

Recall
$|\mathbf{u}(\mathfrak{osp}(1|2))|:\mathbf{u}(\mathfrak{sl}_{2})$ is
an $id$-Frobenius extension. Let $x=e_{21}-e_{13},y=e_{31}+e_{12}$
where  $e_{ij}$ is the unit matrix with $1$ in the $i,j$-entry and
$0$ otherwise. Then the dual projective pair is
$x_{1}=1,x_{2}=x,x_{3}=y,x_{4}=xy+1-[x,y];\;
y_{1}=xy,y_{2}=y,y_{3}=-x,y_{4}=1$. It is straightforward to show
that $\sum_{i=1}^{4}y_{i}x_{i}=1$. For details, see the Example in
page 423 of \cite{BF}.

Define $\psi: M\rightarrow |\mathbf{u}(\mathfrak{osp}(1|2))|\otimes
_{\mathbf{u}(\mathfrak{sl}_{2})} M$ by $m\mapsto 1\otimes m$ for
$m\in M$. It is a morphism of
$\mathbf{u}(\mathfrak{sl}_{2})$-modules. Therefore, the trace map
$\textrm{Tr}_{[|\mathbf{u}(\mathfrak{osp}(1|2))|:\mathbf{u}(\mathfrak{sl}_{2})]}(\psi)$
of $\psi$ is an $|\mathbf{u}(\mathfrak{osp}(1|2))|$-morphism from
$M$ to $|\mathbf{u}(\mathfrak{osp}(1|2))|\otimes
_{\mathbf{u}(\mathfrak{sl}_{2})} M$. By definition,
$$\textrm{Tr}_{[|\mathbf{u}(\mathfrak{osp}(1|2))|:\mathbf{u}(\mathfrak{sl}_{2})]}(\psi)(m)=\sum_{i=1}^{4}y_{i}\otimes x_{i}m$$
for $m\in M$. Then
$\varphi\textrm{Tr}_{[|\mathbf{u}(\mathfrak{osp}(1|2))|:\mathbf{u}(\mathfrak{sl}_{2})]}(\psi)(m)=\sum_{i=1}^{4}y_{i}x_{i}m=m$
for $m\in M$. Consequently,
$$\varphi\textrm{Tr}_{[|\mathbf{u}(\mathfrak{osp}(1|2))|:\mathbf{u}(\mathfrak{sl}_{2})]}(\psi)=id_{M}$$
and thus $M||\mathbf{u}(\mathfrak{osp}(1|2))|\otimes
_{\mathbf{u}(\mathfrak{sl}_{2})} M$.\end{proof}

\subsection{Simples, Projectives and Blocks}

In this subsection, the structures of simple modules, projective
modules and the basic algebras of blocks of
$|\mathbf{u}(\mathfrak{osp}(1|2))|$ are given. As a byproduct, its
Auslander-Reiten quiver is determined.

\subsubsection{Simples and Verma modules} As usual, for a Lie
superalgebra $\mathfrak{g}$, its even (resp. odd) part is denoted by
$\mathfrak{g}_{0}$ (resp. $\mathfrak{g}_{1}$) and
$\mathfrak{g}=\mathfrak{g}_{0}\oplus \mathfrak{g}_{1}$. Recall that
$\mathfrak{g}=\mathfrak{osp}(1|2)$ consists of $3\times 3$ matrices
in the following $(1|2)$-block form
$$\left [ \begin{array}{cccc} 0&v&u\\u&a&b\\-v&c&-a
\end{array}\right]$$
for $a,b,c,u,v\in \kappa$. The even subalgebra
$\mathfrak{osp}(1|2)_{0}$, which is isomorphic to
$\mathfrak{sl}_{2}$, is generated by
$$e=\left [ \begin{array}{cccc} 0&0&0\\0&0&1\\0&0&0
\end{array}\right],\;\;h=\left [ \begin{array}{cccc}
0&0&0\\0&1&0\\0&0&-1
\end{array}\right],\;\;f=\left [ \begin{array}{cccc} 0&0&0\\0&0&0\\0&1&0
\end{array}\right].$$
A basis for the odd part $\mathfrak{osp}(1|2)_{1}$ is given by
$$E=\left [ \begin{array}{cccc} 0&0&1\\1&0&0\\0&0&0
\end{array}\right],\;\;F=\left [ \begin{array}{cccc} 0&1&0\\0&0&0\\-1&0&0
\end{array}\right].$$ The commutation relations of these basis are
collected as follows:
$$[h,E]=E,\;\;[h,F]=-F,\;\;[h,e]=2e,\;\;[h,f]=-2f,$$
$$[e,E]=0,\;\;[e,F]=-E,\;\;[e,f]=h,$$
$$[f,E]=-F,\;\;[f,F]=0,$$
$$[E,E]=2e,\;\;[E,F]=h,\;\;[F,F]=-2f.$$
It is not hard to see that the restricted enveloping algebra
$\mathbf{u}(\mathfrak{osp}(1|2))$ of $\mathfrak{osp}(1|2)$ is
generated by even element $h$ and odd elements $E,F$ with relations
$$EF+FE=h,\;\;hE-Eh=E,\;\;hF-Fh=-F,$$
$$E^{2p}=F^{2p}=0,\;\;h^{p}=h.$$

The structures of simple modules and Verma modules have been given
in a more general context in \cite{WZ}. Let's recall them. For any
$0\leq \lambda\leq p-1$, we define a $2\lambda+1$-dimensional
$|\mathbf{u}(\mathfrak{osp}(1|2))|$-module $V^{\lambda}$ as follows.
This module has a basis $v_{0},v_{1},\ldots,v_{2\lambda}$ and the
actions of the generators are given by the following rules
\begin{equation}
hv_{i}=(\lambda-i)v_{i},\;\;Ev_{i}=\left \{
\begin{array}{ll} -\frac{i}{2} v_{i-1},& \;\; \textrm{if}\; i\; \textrm{is even}\\(\lambda-\frac{i-1}{2})v_{i-1}, &
\;\;\textrm{if}\; i\; \textrm{is odd},
\end{array}\right.\;\;Fv_{i}=v_{i+1}
\end{equation}
where $i=0,1,\ldots,2\lambda$ and $v_{-1}=v_{2\lambda+1}=0$. By
Proposition 6.3 in \cite{WZ}, $\{V^{\lambda}|0\leq \lambda\leq
p-1\}$ forms a complete non-redundant list of simple
$|\mathbf{u}(\mathfrak{osp}(1|2))|$-modules. The graphical
representation for $V^{1}$ is indicated as follows. Similar to the
case of $\mathfrak{sl}_{2}$, arrows and dotted ones show the actions
of $E$ and $F$ respectively.

\begin{figure}[hbt]
\begin{picture}(50,60)(0,0)

 \put(0,55){\makebox(0,0){$\cdot$}}
 \put(2,50){\vector(0,-1){15}}
 \put(-2,50){\vector(0,1){0}}
 \put(-2,43){\makebox(0,0){$\vdots$}}

  \put(0,30){\makebox(0,0){$\cdot$}}

\put(2,25){\vector(0,-1){15}}
 \put(-2,25){\vector(0,1){0}}
 \put(-2,18){\makebox(0,0){$\vdots$}}

  \put(0,5){\makebox(0,0){$\cdot$}}
\end{picture}
\end{figure}

Let $\mathbf{u}_{+}$ and $\mathbf{u}_{-}$ be the subalgebras of
$|\mathbf{u}(\mathfrak{osp}(1|2))|$ generated by $h,E$ and $h,F$
respectively. Also, for any $0\leq \lambda\leq p-1$, we have the
Verma modules $W^{\lambda}$ and $\tilde{W}^{\lambda}$ which are free
over $\mathbf{u}_{+}$ and $\mathbf{u}_{-}$ respectively. They are
given by the following rules.\\
$W^{\lambda}:$
\begin{equation}
hv_{i}=(\lambda-i)v_{i},\;\;Ev_{i}=\left \{
\begin{array}{ll} -\frac{i}{2} v_{i-1},& \;\; \textrm{if}\; i\; \textrm{is even}\\(\lambda-\frac{i-1}{2})v_{i-1}, &
\;\;\textrm{if}\; i\; \textrm{is odd},
\end{array}\right.\;\;Fv_{i}=v_{i+1},
\end{equation}
$\tilde{W}^{\lambda}:$
\begin{equation}
hv_{i}=(i-\lambda)v_{i},\;\;Ev_{i}=v_{i+1},\;\;Fv_{i}=\left \{
\begin{array}{ll} \frac{i}{2} v_{i-1},& \;\; \textrm{if}\; i\; \textrm{is even}\\(\frac{i-1}{2}-\lambda)v_{i-1}, &
\;\;\textrm{if}\; i\; \textrm{is odd},
\end{array}\right.
\end{equation}
where $i=0,1,\ldots,2p-1$ and $v_{-1}=v_{2p}=0$. For
$p=3,\lambda=1$, their graphical representations are as follows.

\begin{figure}[hbt]
\begin{picture}(150,20)(0,0)

 \put(-55,10){\makebox(0,0){$W^{1}:$}}

 \put(20,10){\makebox(0,0){$\cdot$}}
 \put(30,12){\makebox(0,0){$\cdots$}}
 \put(40,12){\vector(1,0){0}}
 \put(40,8){\vector(-1,0){15}}
 \put(45,10){\makebox(0,0){$\cdot$}}

\put(45,10){\makebox(0,0){$\cdot$}}
 \put(55,12){\makebox(0,0){$\cdots$}}
 \put(65,12){\vector(1,0){0}}
 \put(65,8){\vector(-1,0){15}}
 \put(70,10){\makebox(0,0){$\cdot$}}

 \put(80,10){\makebox(0,0){$\cdots$}} \put(90,10){\vector(1,0){0}} \put(95,10){\makebox(0,0){$\cdot$}}

\put(95,10){\makebox(0,0){$\cdot$}}
 \put(105,12){\makebox(0,0){$\cdots$}}
 \put(115,12){\vector(1,0){0}}
 \put(115,8){\vector(-1,0){15}}
 \put(120,10){\makebox(0,0){$\cdot$}}

 \put(120,10){\makebox(0,0){$\cdot$}}
 \put(130,12){\makebox(0,0){$\cdots$}}
 \put(140,12){\vector(1,0){0}}
 \put(140,8){\vector(-1,0){15}}
 \put(145,10){\makebox(0,0){$\cdot$}}

\end{picture}
\end{figure}

\begin{figure}[hbt]
\begin{picture}(150,20)(0,0)

 \put(-55,10){\makebox(0,0){$\tilde{W}^{1}:$}}

 \put(20,10){\makebox(0,0){$\cdot$}}
 \put(25,12){\vector(1,0){15}}
 \put(25,8){\vector(-1,0){0}}
 \put(37,8){\makebox(0,0){$\cdots$}}
 \put(45,10){\makebox(0,0){$\cdot$}}

 \put(45,10){\makebox(0,0){$\cdot$}}
 \put(50,12){\vector(1,0){15}}
 \put(50,8){\vector(-1,0){0}}
 \put(62,8){\makebox(0,0){$\cdots$}}
 \put(70,10){\makebox(0,0){$\cdot$}}

\put(75,10){\vector(1,0){15}}

 \put(95,10){\makebox(0,0){$\cdot$}}
 \put(100,12){\vector(1,0){15}}
 \put(100,8){\vector(-1,0){0}}
 \put(112,8){\makebox(0,0){$\cdots$}}
 \put(120,10){\makebox(0,0){$\cdot$}}
 \put(120,10){\makebox(0,0){$\cdot$}}
 \put(125,12){\vector(1,0){15}}
 \put(125,8){\vector(-1,0){0}}
 \put(137,8){\makebox(0,0){$\cdots$}}
 \put(145,10){\makebox(0,0){$\cdot$}}

\end{picture}
\end{figure}

Clearly, all Verma modules have dimensions $2p$ and we have the
following non-split extensions:
$$0\rightarrow V^{p-1-\lambda}\rightarrow W^{\lambda}\rightarrow V^{\lambda}\rightarrow 0,$$
$$0\rightarrow V^{p-1-\lambda}\rightarrow \tilde{W}^{\lambda}\rightarrow V^{\lambda}\rightarrow 0,$$
for $0\leq \lambda\leq p-1$.

\begin{remark} \emph{Contrast to the $\mathfrak{sl}_{2}$ case, the Verma
modules $W^{\frac{p-1}{2}},\tilde{W}^{\frac{p-1}{2}}$ are special.
Now,
$\Hom_{|\mathbf{u}(\mathfrak{osp}(1|2))|}(W^{\frac{p-1}{2}},W^{\frac{p-1}{2}})\cong
\Hom_{|\mathbf{u}(\mathfrak{osp}(1|2))|}(\tilde{W}^{\frac{p-1}{2}},\tilde{W}^{\frac{p-1}{2}})
\\\cong \kappa[x]/(x^{2})$ while all Verma modules of
$\mathbf{u}(\mathfrak{sl}_{2})$ are bricks, that is, their
endomorphism rings are isomorphic to $\kappa$.}
\end{remark}

\subsubsection{Projective modules} Inspired by the case of
$\mathfrak{sl}_{2}$ and the work given by Xiao \cite{Xiao}, we
define the following modules, which will be shown to form a complete
list of indecomposable projective
$|\mathbf{u}(\mathfrak{osp}(1|2))|$-modules.

 For any $0\leq
\lambda\leq p-1$, we define an
$|\mathbf{u}(\mathfrak{osp}(1|2))|$-module, denoted by
$P^{p-1-\lambda}$, by the following rules. As a space, it has a
basis consisting of $\{b_{i},a_{i},x_{j},y_{j}|0\leq i\leq
2p-2-2\lambda, \;0\leq j\leq 2\lambda\}$. The actions of $h,E,F$ are
given by:

\begin{eqnarray}
hb_{i}&=&(p-1-\lambda-i)b_{i},\;\; Fb_{i}=b_{i+1},\\
Eb_{i}&=&\left \{
\begin{array}{ll} -\frac{i}{2} b_{i-1}+a_{i-1},& \;\; \textrm{if}\; i\; \textrm{is even}\\
(p-1-\lambda-\frac{i-1}{2})b_{i-1}-a_{i-1}, & \;\;\textrm{if}\; i\;
\textrm{is odd};
\end{array}\right.
\end{eqnarray}

\begin{eqnarray}
hy_{j}&=&(j-\lambda)y_{j}, \;\;Ey_{j}=y_{j+1},\\
Fy_{j}&=&\left \{
\begin{array}{ll} \frac{j}{2} y_{j-1},& \;\; \textrm{if}\; j\; \textrm{is even}\\(\frac{j-1}{2}-\lambda)y_{j-1}, &
\;\;\textrm{if}\; j\; \textrm{is odd};
\end{array}\right.
\end{eqnarray}
\begin{eqnarray}
hx_{j}&=&(\lambda-j)x_{j},\;\;Fx_{j}=x_{j+1},\\
Ex_{j}&=&\left \{
\begin{array}{ll} -\frac{j}{2} x_{j-1},& \;\; \textrm{if}\; j\; \textrm{is even}\\(\lambda-\frac{j-1}{2})x_{j-1}, &
\;\;\textrm{if}\; j\; \textrm{is odd};
\end{array}\right.
\end{eqnarray}

\begin{eqnarray}
ha_{i}&=&(p-1-\lambda-i)a_{i},\;\; Fa_{i}=a_{i+1},\\
Ea_{i}&=&\left \{
\begin{array}{ll} -\frac{i}{2} a_{i-1},& \;\; \textrm{if}\; i\; \textrm{is even}\\
(p-1-\lambda-\frac{i-1}{2})a_{i-1}, & \;\;\textrm{if}\; i\;
\textrm{is odd},
\end{array}\right.
\end{eqnarray}
where $b_{2p-1-2\lambda}=y_{2\lambda},
\;y_{2\lambda+1}=a_{2p-1-2\lambda},\;
x_{2\lambda+1}=a_{0},\;b_{-1}=x_{2\lambda}$.

The graphical description of the $P^{p-1-\lambda}$ (for
$p=3,\lambda=1$) is indicated as follows.

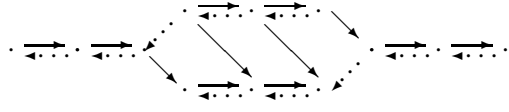
\begin{figure}[hbt]
\begin{picture}(100,40)(0,0)

 \put(0,30){\makebox(0,0){$\cdot$}}
 \put(5,32){\vector(1,0){15}}
 \put(5,28){\vector(-1,0){0}}
  \put(17,28){\makebox(0,0){$\cdots$}}
  \put(5,25){\vector(1,-1){20}}
  \put(30,25){\vector(1,-1){20}}

  \put(25,30){\makebox(0,0){$\cdot$}}
 \put(30,32){\vector(1,0){15}}
 \put(30,28){\vector(-1,0){0}}
  \put(42,28){\makebox(0,0){$\cdots$}}
   \put(50,30){\makebox(0,0){$\cdot$}}

 \put(55,30){\vector(1,-1){10}} \put(70,15){\makebox(0,0){$\cdot$}}
  \put(70,15){\makebox(0,0){$\cdot$}}
 \put(75,17){\vector(1,0){15}}
 \put(75,13){\vector(-1,0){0}}
  \put(87,13){\makebox(0,0){$\cdots$}} \put(95,15){\makebox(0,0){$\cdot$}}
  \put(100,17){\vector(1,0){15}}
 \put(100,13){\vector(-1,0){0}}
  \put(112,13){\makebox(0,0){$\cdots$}} \put(120,15){\makebox(0,0){$\cdot$}}

\put(65,10){\makebox(0,0){$\cdot$}}\put(59,4){\makebox(0,0){$\cdot$}}\put(55,0){\vector(-1,-1){0}}
\put(62,7){\makebox(0,0){$\cdot$}}

 \put(0,0){\makebox(0,0){$\cdot$}}
 \put(5,2){\vector(1,0){15}}
 \put(5,-2){\vector(-1,0){0}}
  \put(17,-2){\makebox(0,0){$\cdots$}}

  \put(25,0){\makebox(0,0){$\cdot$}}
 \put(30,2){\vector(1,0){15}}
 \put(30,-2){\vector(-1,0){0}}
  \put(42,-2){\makebox(0,0){$\cdots$}}
   \put(50,0){\makebox(0,0){$\cdot$}}

   \put(-5,25){\makebox(0,0){$\cdot$}}\put(-11,19){\makebox(0,0){$\cdot$}}\put(-15,15){\vector(-1,-1){0}}
\put(-8,22){\makebox(0,0){$\cdot$}}

 \put(-13,13){\vector(1,-1){10}}

 \put(-40,15){\makebox(0,0){$\cdot$}}
 \put(-35,17){\vector(1,0){15}}
 \put(-35,13){\vector(-1,0){0}}
  \put(-23,13){\makebox(0,0){$\cdots$}}
  \put(-15,15){\makebox(0,0){$\cdot$}}

   \put(-65,15){\makebox(0,0){$\cdot$}}
 \put(-60,17){\vector(1,0){15}}
 \put(-60,13){\vector(-1,0){0}}
  \put(-48,13){\makebox(0,0){$\cdots$}}
  \put(-40,15){\makebox(0,0){$\cdot$}}
\end{picture}
\end{figure}

\begin{proposition} The set $\{P^{\lambda}|0\leq \lambda\leq p-1\}$
gives a complete list of non-isomorphic indecomposable projective
$|\mathbf{u}(\mathfrak{osp}(1|2))|$-modules. All of them have
dimensions $4p$.
\end{proposition}
\begin{proof} The second statement is obvious. For any $0\leq \lambda\leq p-1$, $P^{p-1-\lambda}$ is
clearly indecomposable and its head is isomorphic to
$V^{p-1-\lambda}$. Owing to the classification of simple
$|\mathbf{u}(\mathfrak{osp}(1|2))|$-modules, the conclusion is
proved provided that we can show $P^{p-1-\lambda}$ is projective.
Actually, from $E^{2}=e$ and $F^{2}=-f$ in
$|\mathbf{u}(\mathfrak{osp}(1|2))|$, one can write the actions of
$e,f$ on the basis given in (5.10)-(5.17) directly:
 $$eb_{i}=\left \{\begin{array}{ll}
-\frac{i}{2}(p-\lambda-\frac{i}{2}) b_{i-2}+(p-\lambda)a_{i-2},&
  \textrm{if}\; i\; \textrm{is even}\\
 -\frac{i-1}{2}(p-\lambda-1-\frac{i-1}{2}) b_{i-2}+(p-1-\lambda)a_{i-2}, &\textrm{if}\; i\;
\textrm{is odd},
\end{array}\right. fb_{i}=-b_{i+2};$$
$$ea_{i}=\left \{\begin{array}{ll} -\frac{i}{2}(p-\lambda-\frac{i}{2}) a_{i-2},&
 \;\; \textrm{if}\; i\; \textrm{is even}\\
 -\frac{i-1}{2}(p-\lambda-1-\frac{i-1}{2}) a_{i-2}, & \;\;\textrm{if}\; i\;
\textrm{is odd},
\end{array}\right. fa_{i}=-a_{i+2};$$
$$ey_{j}=y_{j+2},\;\;fy_{j}=\left \{\begin{array}{ll} -\frac{j}{2}(\frac{j}{2}-\lambda-1) y_{j-2},&
 \;\; \textrm{if}\; j\; \textrm{is even}\\
 -\frac{j-1}{2}(\frac{j-1}{2}-\lambda) y_{j-2}, & \;\;\textrm{if}\; j\;
\textrm{is odd},
\end{array}\right.$$
$$ex_{j}=\left \{\begin{array}{ll} -\frac{j}{2}(\lambda+1-\frac{j}{2}) x_{j-2},&
 \;\; \textrm{if}\; j\; \textrm{is even}\\
 -\frac{j-1}{2}(\lambda-\frac{j-1}{2}) x_{j-2}, & \;\;\textrm{if}\; j\;
\textrm{is odd},
\end{array}\right., fx_{j}=-x_{j+2}$$
for $0\leq i\leq 2p-2\lambda-2$ and $0\leq j\leq 2\lambda$. Denote
the restriction of $P^{p-1-\lambda}$ to
$\mathbf{u}(\mathfrak{sl}_{2})$ by
$P^{p-1-\lambda}|_{\mathbf{u}(\mathfrak{sl}_{2})}$. Then, it is not
hard to see that
$$P^{p-1-\lambda}|_{\mathbf{u}(\mathfrak{sl}_{2})}\cong P_{0}^{p-1-\lambda}\oplus P_{0}^{p-2-\lambda}$$
if $\lambda\neq 0,p-1$, and
$$P^{0}|_{\mathbf{u}(\mathfrak{sl}_{2})}\cong P_{0}^{0}\oplus 2V_{0}^{p-1},\;
P^{p-1}|_{\mathbf{u}(\mathfrak{sl}_{2})}\cong P_{0}^{p-2}\oplus
2V_{0}^{p-1}.$$ As a conclusion, the restriction
$P^{p-1-\lambda}|_{\mathbf{u}(\mathfrak{sl}_{2})}$ is projective for
all $0\leq \lambda\leq p-1$. Therefore, Lemma 5.2 implies
$P^{p-1-\lambda}$ itself is projective.
\end{proof}

An indecomposable projective module corresponds to an extension of
Verma modules. Indeed, for any $0\leq \lambda \leq p-1$ we have the
following non-split exact sequences
$$0\rightarrow W^{\lambda}\rightarrow P^{p-1-\lambda}\rightarrow W^{p-1-\lambda}\rightarrow 0,$$
$$0\rightarrow \tilde{W}^{\lambda}\rightarrow P^{p-1-\lambda}\rightarrow \tilde{W}^{p-1-\lambda}\rightarrow 0.$$
This verifies and strengthens the Proposition 6.3 (iii) in
\cite{WZ}, which states $P^{\lambda}$ has a Verma filtration with
$W^{\lambda}$ and $W^{p-1-\lambda}$ as subquotients.

\subsubsection{Blocks and Auslander-Reiten quivers} By above
proposition and the structures of projective modules, we know that
only $V^{\lambda}$ and $V^{p-1-\lambda}$ are composition factors of
$P^{\lambda}$ for $0\leq \lambda \leq p-1$. Thus there are exactly
$\frac{p+1}{2}$ blocks and $V^{\lambda},V^{p-1-\lambda}$ belong to
the same block for $0\leq \lambda \leq \frac{p-1}{2}$. In
particular, the block containing $V^{\frac{p-1}{2}}$ is primary,
that is, it has only one simple module. Our next aim is to describe
the basic algebras of these blocks using quivers with relations. For
more on quivers and related terminologies, see \cite{ARS}.

Take an $\lambda\in \{0,1,\ldots,p-1\}$. By the standard methods
using in representation theory of finite dimensional algebras
\cite{ARS}, the basic algebra of the block containing $V^{\lambda}$
is isomorphic to
$$\End_{|\mathbf{u}(\mathfrak{osp}(1|2))|}(P^{\lambda}\oplus P^{p-1-\lambda})$$
if $\lambda\neq \frac{p-1}{2}$ and isomorphic to
$$\End_{|\mathbf{u}(\mathfrak{osp}(1|2))|}(P^{
\frac{p-1}{2}})$$ otherwise.

Define  $\Lambda_{2}$ to be the algebra given by the following
quiver and relations

\begin{figure}[hbt]
\begin{picture}(200,60)(0,0)

 \put(0,30){\makebox(0,0){$\cdot$}}
  \put(5,32){\vector(1,0){30}} \put(20,35){\makebox(0,0){$x_{1}$}}
  \put(5,42){\vector(1,0){30}}\put(20,45){\makebox(0,0){$y_{1}$}}

  \put(40,30){\makebox(0,0){$\cdot$}} \put(35,28){\vector(-1,0){30}} \put(20,22){\makebox(0,0){$x_{2}$}}
  \put(35,18){\vector(-1,0){30}}\put(20,12){\makebox(0,0){$y_{2}$}}

\put(80,35){\makebox(0,0){$x_{i}x_{j}=y_{i}y_{j}$}}
\put(180,30){\makebox(0,0){$\textrm{for}\; 1\leq i\neq j\leq 2,$}}
\put(90,25){\makebox(0,0){$x_{i}y_{j}=y_{i}x_{j}=0$}}
\end{picture}
\end{figure}
and $\Lambda_{1}$  given by

\begin{figure}[hbt]
\begin{picture}(150,50)(30,0)
\put(25,25){\circle{40}}
\put(5,25){\vector(0,1){0.01}}\put(25,0){\makebox(0,0){$x$}}
\put(44,25){\makebox(0,0){$\cdot$}} \put(65,25){\circle{40}}
\put(85,25){\vector(0,-1){0.01}}\put(65,0){\makebox(0,0){$y$}}

\put(165,25){\makebox(0,0){$x^{2}=y^{2},\;xy=yx=0.$}}
\end{picture}
\end{figure}

\begin{lemma} Assume that $0\leq \lambda\leq p-1$.

\emph{(1)} If $\lambda\neq\frac{p-1}{2}$, then
$\End_{|\mathbf{u}(\mathfrak{osp}(1|2))|}(P^{\lambda}\oplus
P^{p-1-\lambda})\cong \Lambda_{2}$.

\emph{(2)}
$\End_{|\mathbf{u}(\mathfrak{osp}(1|2))|}(P^{\frac{p-1}{2}})\cong
\Lambda_{1}$.
\end{lemma}
\begin{proof} We only prove (1) since (2) can be proved similarly.
For (1), we can represent projective modules $P^{p-1-\lambda}$ and
$P^{\lambda}$ by using the following graphs:

\begin{figure}[hbt]
\begin{picture}(60,70)(20,0)

\put(-55,65){\makebox(0,0){$P^{p-1-\lambda}:$}}
 \put(0,60){\makebox(0,0){$\bullet$}}  \put(0,65){\makebox(0,0){$V^{p-1-\lambda}$}}
 \put(5,55){\line(1,-1){15}}
 \put(25,35){\makebox(0,0){$\circ$}}\put(35,35){\makebox(0,0){$V^{\lambda}$}}
 \put(20,30){\line(-1,-1){15}} \put(-25,35){\makebox(0,0){$\circ$}}
\put(-35,35){\makebox(0,0){$V^{\lambda}$}}
 \put(-5,10){\line(-1,1){15}}
 \put(0,5){\makebox(0,0){$\bullet$}}\put(-5,55){\line(-1,-1){15}}
\put(0,0){\makebox(0,0){$V^{p-1-\lambda}$}}

\put(85,65){\makebox(0,0){$P^{\lambda}:$}}

\put(130,60){\makebox(0,0){$\circ$}}
\put(130,65){\makebox(0,0){$V^{\lambda}$}}
 \put(135,55){\line(1,-1){15}}
 \put(155,35){\makebox(0,0){$\bullet$}}\put(170,35){\makebox(0,0){$V^{p-1-\lambda}$}}
 \put(150,30){\line(-1,-1){15}} \put(105,35){\makebox(0,0){$\bullet$}}
\put(85,35){\makebox(0,0){$V^{p-1-\lambda}$}}
 \put(125,10){\line(-1,1){15}}
 \put(130,5){\makebox(0,0){$\circ$}}\put(125,55){\line(-1,-1){15}}
\put(130,0){\makebox(0,0){$V^{\lambda}$}}
\end{picture}
\end{figure}
From this, one can see that there are exactly two non-trivial linear
independent $|\mathbf{u}(\mathfrak{osp}(1|2))|$-morphisms from
$P^{p-1-\lambda}$ to $P^{\lambda}$:\\

\begin{figure}[hbt]
\begin{picture}(150,30)(0,0)
\put(-40,15){\makebox(0,0){$x_{1}:\;P^{p-1-\lambda}\rightarrow $}}
\put(0,25){\makebox(0,0){$\bullet$}} \put(5,20){\line(1,-1){15}}
\put(25,0){\makebox(0,0){$\circ$}}

\put(90,15){\makebox(0,0){$y_{1}:\;P^{p-1-\lambda}\rightarrow $}}
\put(145,25){\makebox(0,0){$\bullet$}}
\put(140,20){\line(-1,-1){15}} \put(120,0){\makebox(0,0){$\circ$}}
\end{picture}
\end{figure}

Similarly, we also have two non-trivial linear independent
$|\mathbf{u}(\mathfrak{osp}(1|2))|$-morphisms $x_{2},y_{2}$ from
$P^{\lambda}$ to $P^{p-1-\lambda}$:

\begin{picture}(0,35)(-100,0)

\put(-30,15){\makebox(0,0){$x_{2}:\;P^{\lambda}\rightarrow $}}
\put(25,25){\makebox(0,0){$\circ$}} \put(20,20){\line(-1,-1){15}}
\put(0,0){\makebox(0,0){$\bullet$}}

\put(70,15){\makebox(0,0){$y_{2}:\;P^{\lambda}\rightarrow $}}
\put(100,25){\makebox(0,0){$\circ$}} \put(105,20){\line(1,-1){15}}
\put(125,0){\makebox(0,0){$\bullet$}}
\end{picture}

Clearly, such maps indeed generate
$\End_{|\mathbf{u}(\mathfrak{osp}(1|2))|}(P^{\lambda}\oplus
P^{p-1-\lambda})$ and exactly satisfy the relations in the
definition of $\Lambda_{2}$.
\end{proof}

Summarizing, we have proved the following.

\begin{proposition} Let $\kappa$ be an algebraically closed field of
characteristic $p>2$, and $\mathbf{u}(\mathfrak{osp}(1|2))$ the
restricted enveloping algebra of Lie superalgebra
$\mathfrak{osp}(1|2)$ over $\kappa$. Then
 \emph{\item(1)} The
algebra $|\mathbf{u}(\mathfrak{osp}(1|2))|$ has $p$ isomorphism
classes of simple modules, i.e. $V^{\lambda}$ for $0\leq \lambda\leq
p-1$.
 \emph{\item(2)} The algebra
$|\mathbf{u}(\mathfrak{osp}(1|2))|$ has $\frac{p+1}{2}$ blocks.
\emph{\item(3)} The block containing $V^{\frac{p-1}{2}}$ is primary
and its basic algebra is isomorphic to $\Lambda_{1}$.
\emph{\item(4)} For any $0\leq \lambda< \frac{p-1}{2}$, the simple
modules $V^{\lambda}$ and $V^{p-1-\lambda}$ belong to the same
block, denoted by $B(\lambda)$, whose basic algebra is isomorphic to
$\Lambda_{2}$.
\end{proposition}

\begin{remark} \emph{(1)} \emph{Let $A$ be an artin algebra. Operating it by its dual
$D(A):=\Hom_{\kappa}(A,\kappa)$, one can get a new algebra $T(A)$,
called the \emph{trivial extension} of $A$. By definition, the
underlying vector space of $T(A)=A\oplus D(A)$ and the
multiplication is given by $$(a,d)(a',d')=(aa',da'+ad')$$ for
$a,a'\in A, \;d,d'\in D(A)$ by noting $D(A)$ is an $A$-$A$-bimodule
in an obvious way. It is not hard to see that $\Lambda_{2}$ is
indeed the trivial extension of the Kronecker algebra, that is, the
path algebra of the quiver}
\begin{figure}[hbt]
\begin{picture}(30,15)(0,0)

\put(0,0){\makebox(0,0){$\cdot$}}
 \put(5,5){\vector(1,0){15}}

 \put(25,0){\makebox(0,0){$\cdot$}}
 \put(5,-5){\vector(1,0){15}}

\end{picture}
\end{figure}

\emph{(2)}\emph{ By the Example in Section 4 in \cite{Far},
$|\mathbf{u}(\mathfrak{osp}(1|2))|$ is a tame algebra, which is a
direct consequence of our results now. Moreover, one can see that
all
 blocks of $|\mathbf{u}(\mathfrak{osp}(1|2))|$ are tame. This is not
 the case for $\mathbf{u}(\mathfrak{sl}_{2})$, which has exactly one
 block of finite representation type.}
\end{remark}

The categories of finite dimensional representations over algebras
$\Lambda_{1}$ and $\Lambda_{2}$ had been well understood. Recall a
graph is called a \emph{tube} if it is isomorphic to
$\mathbb{Z}A_{\infty}/n$ for some positive integer $n$ and $n$ is
called the \emph{rank} of this tube. A rank $1$ tube is said to be
\emph{homogeneous}. For details about Auslander-Reiten quivers and
translation quivers, see Chapter VII in \cite{ARS} and \cite{Ri}.
The Auslander-Reiten quiver of $\Lambda_{1}$ can be drawn as
follows.

\begin{figure}[hbt]
\begin{picture}(220,80)(0,-10)
\put(-10,40){\makebox(0,0){$\ddots$}}

\put(0,30){\makebox(0,0){$\cdot$}}
 \put(2,35){\vector(1,1){15}}
\put(4,33){\vector(1,1){15}}
\put(22,52){\makebox(0,0){$\cdot$}}\put(24,54){\vector(1,0){15}}
\put(45,54){\makebox(0,0){$\cdot$}}\put(50,54){\vector(1,0){15}}
\put(45,60){\makebox(0,0){$P$}}
 \put(24,48){\vector(1,-1){15}}
\put(28,52){\vector(1,-1){15}} \put(45,30){\makebox(0,0){$\cdot$}}

\put(45,30){\makebox(0,0){$\cdot$}}
 \put(47,35){\vector(1,1){15}}
\put(49,33){\vector(1,1){15}} \put(67,52){\makebox(0,0){$\cdot$}}
 \put(69,48){\vector(1,-1){15}}
\put(73,52){\vector(1,-1){15}} \put(90,30){\makebox(0,0){$\cdot$}}

\put(90,30){\makebox(0,0){$\cdot$}}
 \put(92,35){\vector(1,1){15}}
\put(94,33){\vector(1,1){15}} \put(112,52){\makebox(0,0){$\cdot$}}
\put(122,47){\makebox(0,0){$\ddots$}}

\put(170,70){\makebox(0,0){$\vdots$}}

\put(170,60){\makebox(0,0){$\cdot$}} \put(168,55){\vector(0,-1){15}}
\put(170,35){\makebox(0,0){$\cdot$}} \put(172,40){\vector(0,1){15}}

\put(170,35){\makebox(0,0){$\cdot$}} \put(168,30){\vector(0,-1){15}}
\put(170,10){\makebox(0,0){$\cdot$}} \put(172,15){\vector(0,1){15}}

\put(200,70){\makebox(0,0){$\vdots$}}

\put(200,60){\makebox(0,0){$\cdot$}} \put(198,55){\vector(0,-1){15}}
\put(200,35){\makebox(0,0){$\cdot$}} \put(202,40){\vector(0,1){15}}

\put(200,35){\makebox(0,0){$\cdot$}} \put(198,30){\vector(0,-1){15}}
\put(200,10){\makebox(0,0){$\cdot$}} \put(202,15){\vector(0,1){15}}

\put(220,35){\makebox(0,0){$\cdots$}}
\put(150,35){\makebox(0,0){$\cdots$}}

\put(170,-10){\makebox(0,0){$\textrm{\Small A}\;
\mathbb{P}^{1}\kappa\; \textrm{\Small family of homogeneous
tubes}$}}
\end{picture}
\end{figure}

The Auslander-Reiten quiver of $\Lambda_{2}$ is just the double of
that of $\Lambda_{1}$.

\subsection{Finite dimensional indecomposable modules}

Inspired by the forms of the Auslander-Reiten quivers of the basic
algebras of its blocks and the familiar representation theory of
$\Lambda_{2}$ and $\Lambda_{1}$, we can construct all the
indecomposable representations of
$|\mathbf{u}(\mathfrak{osp}(1|2))|$ now.

\subsubsection{$V^{\lambda}(n)$ and $\tilde{V}^{\lambda}(n)$} For any
positive integer $n$ and $0\leq \lambda \leq p-1$, the basis of
$V^{\lambda}(n)$ is
$$\{a_{u}(m-1),e_{v}(m)|0\leq m\leq n,\;0\leq u\leq2p-2\lambda-2,\;0\leq v\leq 2\lambda\}$$
with actions given by
\begin{eqnarray*}
&&he_{v}(m)=(\lambda-v)e_{v}(m),\;\;Fe_{v}(m)=e_{v+1}(m),\\
&&Ee_{v}(m)=\left \{
\begin{array}{ll} -\frac{v}{2} e_{v-1}(m)+\delta_{v0}a_{2p-2\lambda-2}(m),& \;\; \textrm{if}\; v\; \textrm{is even}\\
(\lambda-\frac{v-1}{2})e_{v-1}(m),& \;\;\textrm{if}\; v\; \textrm{is
odd},
\end{array}\right.\\
&&ha_{u}(m-1)=(p-1-\lambda-u)a_{u}(m-1),\;\;Fa_{u}(m-1)=a_{u+1}(m-1),\\
&&Ea_{u}(m-1)=\left \{
\begin{array}{ll} -\frac{u}{2} a_{u-1}(m-1),& \; \textrm{if}\; u\; \textrm{is even}\\
(p-1-\lambda-\frac{u-1}{2})a_{u-1}(m-1),& \;\textrm{if}\; u\;
\textrm{is odd},
\end{array}\right.
\end{eqnarray*}
where $a_{u}(-1)=a_{u}(n)=0$ for $0\leq u\leq2p-2\lambda-2$,
$a_{-1}(m-1)=a_{2p-2\lambda-1}(m-1)=0$ for $1\leq m\leq n$ and
$e_{2\lambda+1}(m)=a_{0}(m-1)$ for $1\leq m\leq n$. The following is
the graphical description of $V^{\lambda}(n)$ in the case
$n=1,\lambda=1$ and $p=3$:

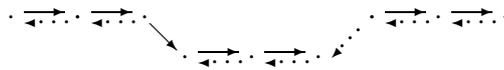
\begin{figure}[hbt]
\begin{picture}(100,40)(0,0)

 \put(70,15){\makebox(0,0){$\cdot$}}
  \put(70,15){\makebox(0,0){$\cdot$}}
 \put(75,17){\vector(1,0){15}}
 \put(75,13){\vector(-1,0){0}}
  \put(87,13){\makebox(0,0){$\cdots$}} \put(95,15){\makebox(0,0){$\cdot$}}
  \put(100,17){\vector(1,0){15}}
 \put(100,13){\vector(-1,0){0}}
  \put(112,13){\makebox(0,0){$\cdots$}} \put(120,15){\makebox(0,0){$\cdot$}}

\put(65,10){\makebox(0,0){$\cdot$}}\put(59,4){\makebox(0,0){$\cdot$}}\put(55,0){\vector(-1,-1){0}}
\put(62,7){\makebox(0,0){$\cdot$}}

 \put(0,0){\makebox(0,0){$\cdot$}}
 \put(5,2){\vector(1,0){15}}
 \put(5,-2){\vector(-1,0){0}}
  \put(17,-2){\makebox(0,0){$\cdots$}}

  \put(25,0){\makebox(0,0){$\cdot$}}
 \put(30,2){\vector(1,0){15}}
 \put(30,-2){\vector(-1,0){0}}
  \put(42,-2){\makebox(0,0){$\cdots$}}
   \put(50,0){\makebox(0,0){$\cdot$}}

 \put(-13,13){\vector(1,-1){10}}

 \put(-40,15){\makebox(0,0){$\cdot$}}
 \put(-35,17){\vector(1,0){15}}
 \put(-35,13){\vector(-1,0){0}}
  \put(-23,13){\makebox(0,0){$\cdots$}}
  \put(-15,15){\makebox(0,0){$\cdot$}}

   \put(-65,15){\makebox(0,0){$\cdot$}}
 \put(-60,17){\vector(1,0){15}}
 \put(-60,13){\vector(-1,0){0}}
  \put(-48,13){\makebox(0,0){$\cdots$}}
  \put(-40,15){\makebox(0,0){$\cdot$}}
\end{picture}
\end{figure}

For any positive integer $n$ and $0\leq \lambda \leq p-1$, the basis
of $\tilde{V}^{\lambda}(n)$ is
$$\{a_{u}(m-1),e_{v}(m)|0\leq m\leq n,\;0\leq u\leq2p-2\lambda-2,\;0\leq v\leq 2\lambda\}$$
with actions given by
\begin{eqnarray*}
&&he_{v}(m)=(\lambda-v)e_{v}(m),\;\;Fe_{v}(m)=e_{v+1}(m),\\
&&Ee_{v}(m)=\left \{
\begin{array}{ll} -\frac{v}{2} e_{v-1}(m),& \;\; \textrm{if}\; v\; \textrm{is even}\\
(\lambda-\frac{v-1}{2})e_{v-1}(m),& \;\;\textrm{if}\; v\; \textrm{is
odd},
\end{array}\right.\\
&&ha_{u}(m-1)=(p-1-\lambda-u)a_{u}(m-1),\;\;Fa_{u}(m-1)=a_{u+1}(m-1),\\
&&Ea_{u}(m-1)=\left \{
\begin{array}{ll} -\frac{u}{2} a_{u-1}(m-1)+\delta_{u0}e_{2\lambda}(m),& \; \textrm{if}\; u\; \textrm{is even}\\
(p-1-\lambda-\frac{u-1}{2})a_{u-1}(m-1),& \;\textrm{if}\; u\;
\textrm{is odd},
\end{array}\right.
\end{eqnarray*}
where $a_{u}(-1)=a_{u}(n)=0$ for $0\leq u\leq2p-2\lambda-2$,
$e_{-1}(m-1)=e_{2\lambda+1}(m-1)=0$ for $0\leq m\leq n$ and
$e_{0}(m-1)=a_{2p-2\lambda-1}(m-1)$ for $1\leq m\leq n$. The
following is the graphical description of $\tilde{V}^{\lambda}(n)$
in the case $n=1,\lambda=1$ and $p=3$:

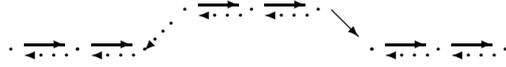
\begin{figure}[hbt]
\begin{picture}(100,40)(0,0)

 \put(0,30){\makebox(0,0){$\cdot$}}
 \put(5,32){\vector(1,0){15}}
 \put(5,28){\vector(-1,0){0}}
  \put(17,28){\makebox(0,0){$\cdots$}}

  \put(25,30){\makebox(0,0){$\cdot$}}
 \put(30,32){\vector(1,0){15}}
 \put(30,28){\vector(-1,0){0}}
  \put(42,28){\makebox(0,0){$\cdots$}}
   \put(50,30){\makebox(0,0){$\cdot$}}

 \put(55,30){\vector(1,-1){10}} \put(70,15){\makebox(0,0){$\cdot$}}
  \put(70,15){\makebox(0,0){$\cdot$}}
 \put(75,17){\vector(1,0){15}}
 \put(75,13){\vector(-1,0){0}}
  \put(87,13){\makebox(0,0){$\cdots$}} \put(95,15){\makebox(0,0){$\cdot$}}
  \put(100,17){\vector(1,0){15}}
 \put(100,13){\vector(-1,0){0}}
  \put(112,13){\makebox(0,0){$\cdots$}} \put(120,15){\makebox(0,0){$\cdot$}}

   \put(-5,25){\makebox(0,0){$\cdot$}}\put(-11,19){\makebox(0,0){$\cdot$}}\put(-15,15){\vector(-1,-1){0}}
\put(-8,22){\makebox(0,0){$\cdot$}}

 \put(-40,15){\makebox(0,0){$\cdot$}}
 \put(-35,17){\vector(1,0){15}}
 \put(-35,13){\vector(-1,0){0}}
  \put(-23,13){\makebox(0,0){$\cdots$}}
  \put(-15,15){\makebox(0,0){$\cdot$}}

   \put(-65,15){\makebox(0,0){$\cdot$}}
 \put(-60,17){\vector(1,0){15}}
 \put(-60,13){\vector(-1,0){0}}
  \put(-48,13){\makebox(0,0){$\cdots$}}
  \put(-40,15){\makebox(0,0){$\cdot$}}
\end{picture}
\end{figure}

For $n\geq 1$, the induced Auslander-Reiten sequences are
$$0\rightarrow V^{\lambda}(n)\rightarrow V^{\lambda}(n+1)\oplus V^{\lambda}(n+1)\rightarrow V^{\lambda}(n+2)\rightarrow 0,$$
$$0\rightarrow \tilde{V}^{\lambda}(n+2)\rightarrow \tilde{V}^{\lambda}(n+1)\oplus \tilde{V}^{\lambda}(n+1)\rightarrow \tilde{V}^{\lambda}(n)\rightarrow 0,$$
$$0\rightarrow \tilde{V}^{\lambda}(1)\rightarrow V^{\lambda}\oplus P^{p-1-\lambda}\oplus V^{\lambda} \rightarrow V^{\lambda}(1)\rightarrow 0.$$
Note that $V^{\lambda}(0)=V^{\lambda}=\tilde{V}^{\lambda}(0)$. The
Auslander-Reiten translation is given by
$$\tau V^{\lambda}(n+2)=V^{\lambda}(n),\;\tau \tilde{V}^{\lambda}(n)=\tilde{V}^{\lambda}(n+2),\;\textrm{for}\; n\geq 0
\;\textrm{and}\;\tau V^{\lambda}(1)=\tilde{V}^{\lambda}(1).$$ It is
not hard to see that they indeed give the preprojective component,
showing as follows, of the Auslander-Reiten quiver described after
Remark 5.7.

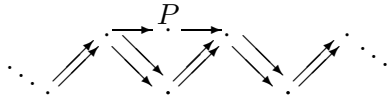
\begin{figure}[hbt]
\begin{picture}(150,60)(0,10)
\put(-10,40){\makebox(0,0){$\ddots$}}

\put(0,30){\makebox(0,0){$\cdot$}}
 \put(2,35){\vector(1,1){15}}
\put(4,33){\vector(1,1){15}}
\put(22,52){\makebox(0,0){$\cdot$}}\put(24,54){\vector(1,0){15}}
\put(45,54){\makebox(0,0){$\cdot$}}\put(50,54){\vector(1,0){15}}
\put(45,60){\makebox(0,0){$P$}}
 \put(24,48){\vector(1,-1){15}}
\put(28,52){\vector(1,-1){15}} \put(45,30){\makebox(0,0){$\cdot$}}

\put(45,30){\makebox(0,0){$\cdot$}}
 \put(47,35){\vector(1,1){15}}
\put(49,33){\vector(1,1){15}} \put(67,52){\makebox(0,0){$\cdot$}}
 \put(69,48){\vector(1,-1){15}}
\put(73,52){\vector(1,-1){15}} \put(90,30){\makebox(0,0){$\cdot$}}

\put(90,30){\makebox(0,0){$\cdot$}}
 \put(92,35){\vector(1,1){15}}
\put(94,33){\vector(1,1){15}} \put(112,52){\makebox(0,0){$\cdot$}}
\put(122,47){\makebox(0,0){$\ddots$}}

\end{picture}
\end{figure}

\subsubsection{$W^{\lambda}(n)$ and $\tilde{W}^{\lambda}(n)$} For any
positive integer $n$ and $0\leq \lambda \leq p-1$, the basis of
$W^{\lambda}(n)$ is
$$\{e_{u}(m)|1\leq m\leq n,\;0\leq u\leq 2p-1\}$$
with actions given by
\begin{eqnarray*}
&&he_{u}(m)=(\lambda-u)e_{u}(m),\;\;Fe_{u}(m)=e_{u+1}(m),\\
&&Ee_{u}(m)=\left \{
\begin{array}{ll} -\frac{u}{2} e_{u-1}(m)+\delta_{u0}e_{2p-1}(m+1),& \;\; \textrm{if}\; u\; \textrm{is even}\\
(\lambda-\frac{u-1}{2})e_{u-1}(m),& \;\;\textrm{if}\; u\; \textrm{is
odd},
\end{array}\right.
\end{eqnarray*}
where $e_{u}(n+1)=0$ for $0\leq u\leq 2p-1$, $e_{2p}(m)=0$ for
$1\leq m\leq n$. The following is the graphical description of
$W^{\lambda}(n)$ in the case $n=2,\lambda=1$ and $p=3$:

\begin{figure}[hbt]
\begin{picture}(100,40)(0,0)

 \put(0,30){\makebox(0,0){$\cdot$}}
 \put(5,32){\vector(1,0){15}}
 \put(5,28){\vector(-1,0){0}}
  \put(17,28){\makebox(0,0){$\cdots$}}

  \put(25,30){\makebox(0,0){$\cdot$}}
 \put(30,32){\vector(1,0){15}}
 \put(30,28){\vector(-1,0){0}}
  \put(42,28){\makebox(0,0){$\cdots$}}
   \put(50,30){\makebox(0,0){$\cdot$}}

 \put(55,30){\vector(1,-1){10}} \put(70,15){\makebox(0,0){$\cdot$}}
  \put(70,15){\makebox(0,0){$\cdot$}}
 \put(75,17){\vector(1,0){15}}
 \put(75,13){\vector(-1,0){0}}
  \put(87,13){\makebox(0,0){$\cdots$}} \put(95,15){\makebox(0,0){$\cdot$}}
  \put(100,17){\vector(1,0){15}}
 \put(100,13){\vector(-1,0){0}}
  \put(112,13){\makebox(0,0){$\cdots$}} \put(120,15){\makebox(0,0){$\cdot$}}

   \put(-5,25){\makebox(0,0){$\cdot$}}\put(-11,19){\makebox(0,0){$\cdot$}}\put(-15,15){\vector(-1,-1){0}}
\put(-8,22){\makebox(0,0){$\cdot$}}

 \put(-40,15){\makebox(0,0){$\cdot$}}
 \put(-35,17){\vector(1,0){15}}
 \put(-35,13){\vector(-1,0){0}}
  \put(-23,13){\makebox(0,0){$\cdots$}}
  \put(-15,15){\makebox(0,0){$\cdot$}}

   \put(-65,15){\makebox(0,0){$\cdot$}}
 \put(-60,17){\vector(1,0){15}}
 \put(-60,13){\vector(-1,0){0}}
  \put(-48,13){\makebox(0,0){$\cdots$}}
  \put(-40,15){\makebox(0,0){$\cdot$}}

 \put(130,30){\makebox(0,0){$\cdot$}}
 \put(135,32){\vector(1,0){15}}
 \put(135,28){\vector(-1,0){0}}
  \put(147,28){\makebox(0,0){$\cdots$}}

  \put(155,30){\makebox(0,0){$\cdot$}}
 \put(160,32){\vector(1,0){15}}
 \put(160,28){\vector(-1,0){0}}
  \put(172,28){\makebox(0,0){$\cdots$}}
   \put(180,30){\makebox(0,0){$\cdot$}}

  \put(122,22){\makebox(0,0){$\cdot$}}    \put(126,26){\makebox(0,0){$\cdot$}}\put(119,19){\makebox(0,0){$\cdot$}}\put(117,17){\vector(-1,-1){0}}
 \put(124,24){\makebox(0,0){$\cdot$}}

\end{picture}
\end{figure}

 For any positive integer $n$ and
$0\leq \lambda \leq p-1$, the basis of $\tilde{W}^{\lambda}(n)$ is
$$\{f_{u}(m)|1\leq m\leq n,\;0\leq u\leq 2p-1\}$$
with actions given by
\begin{eqnarray*}
&&hf_{u}(m)=(\lambda-u)f_{u}(m),\;\;Ef_{u}(m)=f_{u+1}(m),\\
&&Ff_{u}(m)=\left \{
\begin{array}{ll} \frac{u}{2} f_{u-1}(m)+\delta_{u0}f_{2p-1}(m-1),& \;\; \textrm{if}\; u\; \textrm{is even}\\
(\frac{u-1}{2}-\lambda)f_{u-1}(m),& \;\;\textrm{if}\; u\; \textrm{is
odd},
\end{array}\right.
\end{eqnarray*}
where $f_{u}(0)=0$ for $0\leq u\leq 2p-1$, $f_{2p}(m)=0$ for $1\leq
m\leq n$. The following is the graphical description of
$\tilde{W}^{\lambda}(n)$ in the case $n=2,\lambda=1$ and $p=3$:
\begin{figure}[hbt]
\begin{picture}(100,30)(0,0)

 \put(70,15){\makebox(0,0){$\cdot$}}
  \put(70,15){\makebox(0,0){$\cdot$}}
 \put(75,17){\vector(1,0){15}}
 \put(75,13){\vector(-1,0){0}}
  \put(87,13){\makebox(0,0){$\cdots$}} \put(95,15){\makebox(0,0){$\cdot$}}
  \put(100,17){\vector(1,0){15}}
 \put(100,13){\vector(-1,0){0}}
  \put(112,13){\makebox(0,0){$\cdots$}} \put(120,15){\makebox(0,0){$\cdot$}}

\put(65,10){\makebox(0,0){$\cdot$}}\put(59,4){\makebox(0,0){$\cdot$}}\put(55,0){\vector(-1,-1){0}}
\put(62,7){\makebox(0,0){$\cdot$}}

 \put(0,0){\makebox(0,0){$\cdot$}}
 \put(5,2){\vector(1,0){15}}
 \put(5,-2){\vector(-1,0){0}}
  \put(17,-2){\makebox(0,0){$\cdots$}}

  \put(25,0){\makebox(0,0){$\cdot$}}
 \put(30,2){\vector(1,0){15}}
 \put(30,-2){\vector(-1,0){0}}
  \put(42,-2){\makebox(0,0){$\cdots$}}
   \put(50,0){\makebox(0,0){$\cdot$}}

 \put(-13,13){\vector(1,-1){10}}

 \put(-40,15){\makebox(0,0){$\cdot$}}
 \put(-35,17){\vector(1,0){15}}
 \put(-35,13){\vector(-1,0){0}}
  \put(-23,13){\makebox(0,0){$\cdots$}}
  \put(-15,15){\makebox(0,0){$\cdot$}}

   \put(-65,15){\makebox(0,0){$\cdot$}}
 \put(-60,17){\vector(1,0){15}}
 \put(-60,13){\vector(-1,0){0}}
  \put(-48,13){\makebox(0,0){$\cdots$}}
  \put(-40,15){\makebox(0,0){$\cdot$}}

\put(134,0){\makebox(0,0){$\cdot$}}
 \put(139,2){\vector(1,0){15}}
 \put(139,-2){\vector(-1,0){0}}
  \put(151,-2){\makebox(0,0){$\cdots$}}

  \put(159,0){\makebox(0,0){$\cdot$}}
 \put(164,2){\vector(1,0){15}}
 \put(164,-2){\vector(-1,0){0}}
  \put(176,-2){\makebox(0,0){$\cdots$}}
   \put(184,0){\makebox(0,0){$\cdot$}}
 \put(121,13){\vector(1,-1){10}}

\end{picture}
\end{figure}

For $n\geq 1$, the induced Auslander-Reiten sequences are
$$0\rightarrow W^{\lambda}(n)\rightarrow W^{\lambda}(n+1)\oplus W^{\lambda}(n-1)\rightarrow W^{\lambda}(n)\rightarrow 0,$$
$$0\rightarrow \tilde{W}^{\lambda}(n)\rightarrow \tilde{W}^{\lambda}(n+1)\oplus \tilde{W}^{\lambda}(n-1)\rightarrow \tilde{W}^{\lambda}(n)\rightarrow 0.$$
Here we define $W^{\lambda}(0)=\tilde{W}^{\lambda}(0)=0$. The
Auslander-Reiten translation is given by
$$\tau W^{\lambda}(n)=W^{\lambda}(n),\;\tau \tilde{W}^{\lambda}(n)=\tilde{W}^{\lambda}(n)\;\textrm{for}\;n\geq 1.$$

\subsubsection{$T^{\lambda}(s,n)$} For  $n\in \mathbb{Z}^{+}$,
$0\leq \lambda \leq p-1$ and $s=(s_{1},s_{2})\in \kappa^{\ast}\times
\kappa^{\ast}$, the basis of ${T}^{\lambda}(s,n)$ is
$$\{e_{u}(m),\hat{e}_{u}(m)|1\leq m\leq n,\;0\leq u\leq2p-1\}$$
with actions given by
\begin{eqnarray*}
&&he_{u}(m)=(\lambda-u)e_{u}(m),\;\;Fe_{u}(m)=e_{u+1}(m),\\
&&Ee_{u}(m)=\left \{
\begin{array}{ll} -\frac{u}{2} e_{u-1}(m)+s_{1}\delta_{u0}\hat{e}_{2p-1}(m)+\delta_{u0}\hat{e}_{2p-1}(m-1),& \;\; \textrm{if}\; u\; \textrm{is even}\\
(\lambda-\frac{u-1}{2})e_{u-1}(m),& \;\;\textrm{if}\; u\; \textrm{is
odd},
\end{array}\right.\\
&&h\hat{e}_{u}(m)=(\lambda-u)\hat{e}_{u}(m),\;\;F\hat{e}_{u}(m)=\hat{e}_{u+1}(m),\\
&&E\hat{e}_{u}(m)=\left \{
\begin{array}{ll} -\frac{u}{2} \hat{e}_{u-1}(m)+s_{2}\delta_{u0}e_{2p-1}(m)+\delta_{u0}e_{2p-1}(m-1),& \;\; \textrm{if}\; u\; \textrm{is even}\\
(\lambda-\frac{u-1}{2})\hat{e}_{u-1}(m),& \;\;\textrm{if}\; u\;
\textrm{is odd},
\end{array}\right.
\end{eqnarray*}
where $\kappa^{\ast}:=\kappa\setminus \{0\}$,
$e_{u}(0)=\hat{e}_{u}(0)=0$ for $0\leq u\leq2p-1$, and
$e_{2p}(m)=\hat{e}_{2p}(m)=0$ for $1\leq m\leq n$. The following is
the graphical description of $T^{\lambda}(s,n)$ in the case
$n=2,\lambda=1$ and $p=3$:

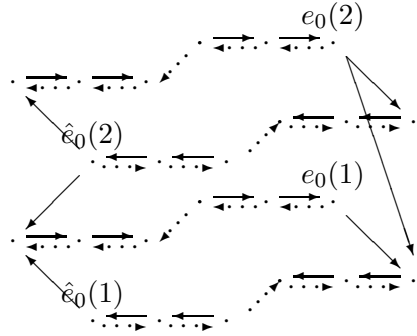
\begin{figure}[hbt]
\begin{picture}(100,130)(30,0)

 \put(70,115){\makebox(0,0){$\cdot$}}
  \put(70,115){\makebox(0,0){$\cdot$}}
 \put(75,117){\vector(1,0){15}}
 \put(75,113){\vector(-1,0){0}}
  \put(87,113){\makebox(0,0){$\cdots$}} \put(95,115){\makebox(0,0){$\cdot$}}
  \put(100,117){\vector(1,0){15}}
 \put(100,113){\vector(-1,0){0}}
  \put(112,113){\makebox(0,0){$\cdots$}} \put(120,115){\makebox(0,0){$\cdot$}}
  \put(120,125){\makebox(0,0){${ e_{0}(2)}$}} \put(125,110){\vector(1,-1){20}}
  \put(125,110){\vector(1,-3){25}}
\put(65,110){\makebox(0,0){$\cdot$}}\put(59,104){\makebox(0,0){$\cdot$}}\put(55,100){\vector(-1,-1){0}}
\put(62,107){\makebox(0,0){$\cdot$}}
 \put(0,100){\makebox(0,0){$\cdot$}}
 \put(5,102){\vector(1,0){15}}
 \put(5,98){\vector(-1,0){0}}
  \put(17,98){\makebox(0,0){$\cdots$}}
  \put(25,100){\makebox(0,0){$\cdot$}}
 \put(30,102){\vector(1,0){15}}
 \put(30,98){\vector(-1,0){0}}
  \put(42,98){\makebox(0,0){$\cdots$}}
   \put(50,100){\makebox(0,0){$\cdot$}}

    \put(70,55){\makebox(0,0){$\cdot$}}
  \put(70,55){\makebox(0,0){$\cdot$}}
 \put(75,57){\vector(1,0){15}}
 \put(75,53){\vector(-1,0){0}}
  \put(87,53){\makebox(0,0){$\cdots$}} \put(95,55){\makebox(0,0){$\cdot$}}
  \put(100,57){\vector(1,0){15}}
 \put(100,53){\vector(-1,0){0}}
  \put(112,53){\makebox(0,0){$\cdots$}} \put(120,55){\makebox(0,0){$\cdot$}}
  \put(120,65){\makebox(0,0){${ e_{0}(1)}$}}\put(125,50){\vector(1,-1){20}}
\put(65,50){\makebox(0,0){$\cdot$}}\put(59,44){\makebox(0,0){$\cdot$}}\put(55,40){\vector(-1,-1){0}}
\put(62,47){\makebox(0,0){$\cdot$}}
 \put(0,40){\makebox(0,0){$\cdot$}}
 \put(5,42){\vector(1,0){15}}
 \put(5,38){\vector(-1,0){0}}
  \put(17,38){\makebox(0,0){$\cdots$}}
  \put(25,40){\makebox(0,0){$\cdot$}}
 \put(30,42){\vector(1,0){15}}
 \put(30,38){\vector(-1,0){0}}
  \put(42,38){\makebox(0,0){$\cdots$}}
   \put(50,40){\makebox(0,0){$\cdot$}}

    \put(100,85){\makebox(0,0){$\cdot$}}
  \put(100,85){\makebox(0,0){$\cdot$}}
 \put(120,87){\vector(-1,0){15}}
 \put(120,83){\vector(1,0){0}}
  \put(110,83){\makebox(0,0){$\cdots$}} \put(125,85){\makebox(0,0){$\cdot$}}
  \put(145,87){\vector(-1,0){15}}
 \put(145,83){\vector(1,0){0}}
  \put(135,83){\makebox(0,0){$\cdots$}} \put(150,85){\makebox(0,0){$\cdot$}}
  \put(95,80){\makebox(0,0){$\cdot$}}\put(89,74){\makebox(0,0){$\cdot$}}\put(100,85){\vector(1,1){0}}
\put(92,77){\makebox(0,0){$\cdot$}}
 \put(30,70){\makebox(0,0){$\cdot$}}\put(30,80){\makebox(0,0){$\hat{e}_{0}(2)$}}\put(25,75){\vector(-1,1){20}}
 \put(25,65){\vector(-1,-1){20}}
 \put(50,72){\vector(-1,0){15}}
 \put(50,68){\vector(1,0){0}}
  \put(40,68){\makebox(0,0){$\cdots$}}
  \put(55,70){\makebox(0,0){$\cdot$}}
 \put(75,72){\vector(-1,0){15}}
 \put(75,68){\vector(1,0){0}}
  \put(65,68){\makebox(0,0){$\cdots$}}
   \put(80,70){\makebox(0,0){$\cdot$}}

      \put(100,25){\makebox(0,0){$\cdot$}}
  \put(100,25){\makebox(0,0){$\cdot$}}
 \put(120,27){\vector(-1,0){15}}
 \put(120,23){\vector(1,0){0}}
  \put(110,23){\makebox(0,0){$\cdots$}} \put(125,25){\makebox(0,0){$\cdot$}}
  \put(145,27){\vector(-1,0){15}}
 \put(145,23){\vector(1,0){0}}
  \put(135,23){\makebox(0,0){$\cdots$}} \put(150,25){\makebox(0,0){$\cdot$}}
  \put(95,20){\makebox(0,0){$\cdot$}}\put(89,14){\makebox(0,0){$\cdot$}}\put(100,25){\vector(1,1){0}}
\put(92,17){\makebox(0,0){$\cdot$}}
 \put(30,10){\makebox(0,0){$\cdot$}}\put(30,20){\makebox(0,0){$\hat{e}_{0}(1)$}}\put(25,15){\vector(-1,1){20}}
 \put(50,12){\vector(-1,0){15}}
 \put(50,8){\vector(1,0){0}}
  \put(40,8){\makebox(0,0){$\cdots$}}
  \put(55,10){\makebox(0,0){$\cdot$}}
 \put(75,12){\vector(-1,0){15}}
 \put(75,8){\vector(1,0){0}}
  \put(65,8){\makebox(0,0){$\cdots$}}
   \put(80,10){\makebox(0,0){$\cdot$}}
\end{picture}
\end{figure}

For $n\geq 1$, the induces Auslander-Reiten sequence is
$$0\rightarrow T^{\lambda}(s,n)\rightarrow T^{\lambda}(s,n+1)\oplus T^{\lambda}(s,n-1)\rightarrow T^{\lambda}(s,n)\rightarrow 0,$$
where $T^{\lambda}(s,0)=0$. And the Auslander-Reiten translation is
given by
$$\tau T^{\lambda}(s,n)=T^{\lambda}(s,n)\;\;\textrm{for}\;n\geq 1.$$

\begin{proposition} Let $s=(s_{1},s_{2}),t=(t_{1},t_{2})\in \kappa^{\ast}\times \kappa^{\ast}$. Then $ T^{\lambda}(s,n)\cong  T^{\lambda}(t,m)$
if and only if $m=n$ and $\frac{s_{1}}{t_{1}}=\frac{t_{2}}{s_{2}}$.
\end{proposition}
\begin{proof}
$``\Rightarrow"$ Comparing the dimensions of both modules, we have
$m=n$. As a case to explain our understanding, there is no harm to
assume that $T^{\lambda}(s,1)\cong T^{\lambda}(t,1)$. Denote the
basis of $T^{\lambda}(t,1)$ by $\{e_{u}(1)',\hat{e}_{u}(1)'|0\leq
u\leq 2p-1\}$ and the isomorphism from $T^{\lambda}(s,1)$ to
$T^{\lambda}(t,1)$ by $\varphi$. It is not hard to see that one can
assume that
$$\varphi(e_{u})=c_{1}e'_{u},\;\;\varphi(\hat{e}_{u}(1))=c_{2}\hat{e}'_{u}(1),$$
or
$$\varphi(e_{u}(1))=c_{1}\hat{e}'_{u}(1),\;\;\varphi(\hat{e}_{u}(1))=c_{2}e'_{u}(1),$$ for some
$c_{1},c_{2}\in \kappa^{\ast}$. In the first case, by
$\varphi(Ee_{0}(1))=E\varphi(e_{0}(1))$ and
$\varphi(E\hat{e}_{0}(1))=E\varphi(\hat{e}_{0}(1))$, we have
$$s_{1}c_{2}=c_{1}t_{1},\;\;s_{2}c_{1}=c_{2}t_{2}$$
which implies
$\frac{s_{1}}{t_{1}}=\frac{c_{1}}{c_{2}}=\frac{t_{2}}{s_{2}}.$ In
the second case, also from $\varphi(Ee_{0}(1))=E\varphi(e_{0}(1))$
and $\varphi(E\hat{e}_{0}(1))=E\varphi(\hat{e}_{0}(1))$, one can
show that $$s_{1}c_{2}=c_{1}t_{2},\;\;s_{2}c_{1}=c_{2}t_{1}$$ and so
$\frac{s_{1}}{t_{2}}=\frac{c_{1}}{c_{2}}=\frac{t_{1}}{s_{2}}$.

$``\Leftarrow"$ Conversely, define
$$\varphi:\;T^{\lambda}(s,1)\rightarrow T^{\lambda}(t,1),\;\;e_{u}(1)\mapsto \frac{s_{1}}{t_{1}} e'_{u}(1),\;
\hat{e}_{u}(1)\mapsto \hat{e}'_{u}(1)$$ for $0\leq u\leq 2p-1$. It
is direct to show that $\varphi$ is a morphism and bijective. Thus
$T^{\lambda}(s,1)\cong T^{\lambda}(t,1)$. From the Ausanlder-Reiten
sequences we constructed,   $T^{\lambda}(s,n)\cong T^{\lambda}(t,n)$
for any $n\geq 1$.

\end{proof}

Not that $\frac{s_{1}}{t_{1}}=\frac{t_{2}}{s_{2}}$ is equivalent to
$s_{1}s_{2}=t_{1}t_{2}$. For any $c\in \kappa^{\ast}$, define
$T^{\lambda}_{c}(n)$ to be any one of $T^{\lambda}(s,n)$ satisfying
$s_{1}s_{2}=c$. Proposition 5.8 implies that
$\{T^{\lambda}_{c}(n)|c\in \kappa^{\ast},n\geq 1\}$ forms a complete
set of representatives of modules
$\{T^{\lambda}(s,n)|s=(s_{1},s_{2})\in \kappa^{\ast}\times
\kappa^{\ast}, n\geq 1\}$ for any fixed  $\lambda\in
\{0,1,\ldots,2p-1\}$.

\begin{remark} \emph{Similar to the case of  5.4.1, 5.4.2,
 one also can define an indecomposable module  $\tilde{T}^{\lambda}(s,n)$ for  $n\in \mathbb{Z}^{+}$,  $0\leq \lambda \leq p-1$ and
$s=(s_{1},s_{2})\in \kappa^{\ast}\times \kappa^{\ast}$. Similarly,
the basis of ${\tilde{T}}^{\lambda}(s,n)$ is
$$\{f_{u}(m),\hat{f}_{u}(m)|1\leq m\leq n,\;0\leq u\leq2p-1\}$$
with actions given by}
\begin{eqnarray*}
&&hf_{u}(m)=(u-\lambda)e_{u}(m),\;\;Ef_{u}(m)=f_{u+1}(m),\\
&&Ff_{u}(m)=\left \{
\begin{array}{ll} \frac{u}{2} f_{u-1}(m)+s_{1}\delta_{u0}\hat{f}_{2p-1}(m)+\delta_{u0}\hat{f}_{2p-1}(m-1),& \;\; \textrm{if}\; u\; \textrm{is even}\\
(\frac{u-1}{2}-\lambda)f_{u-1}(m),& \;\;\textrm{if}\; u\; \textrm{is
odd},
\end{array}\right.\\
&&h\hat{f}_{u}(m)=(u-\lambda)\hat{f}_{u}(m),\;\;E\hat{f}_{u}(m)=\hat{f}_{u+1}(m),\\
&&F\hat{f}_{u}(m)=\left \{
\begin{array}{ll} \frac{u}{2} \hat{f}_{u-1}(m)+s_{2}\delta_{u0}f_{2p-1}(m)+\delta_{u0}f_{2p-1}(m-1),& \;\; \textrm{if}\; u\; \textrm{is even}\\
(\frac{u-1}{2}-\lambda)\hat{f}_{u-1}(m),& \;\;\textrm{if}\; u\;
\textrm{is odd},
\end{array}\right.
\end{eqnarray*}
\emph{where $f_{u}(0)=\hat{f}_{u}(0)=0$ for $0\leq u\leq2p-1$, and
$f_{2p}(m)=\hat{f}_{2p}(m)=0$ for $1\leq m\leq n$. The following is
the graphical description of $\tilde{T}^{\lambda}(s,n)$ in the case
$n=2,\lambda=1$ and $p=3$:}\\

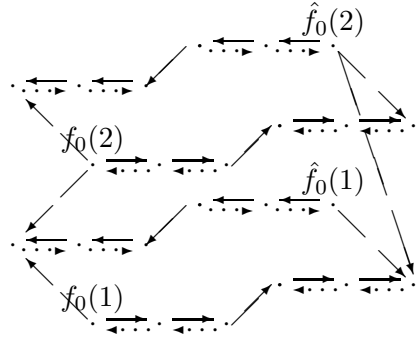
\begin{figure}[hbt]
\begin{picture}(100,130)(30,0)

 \put(70,115){\makebox(0,0){$\cdot$}}
  \put(70,115){\makebox(0,0){$\cdot$}}
 \put(90,117){\vector(-1,0){15}}
 \put(90,113){\vector(1,0){0}}
  \put(80,113){\makebox(0,0){$\cdots$}} \put(95,115){\makebox(0,0){$\cdot$}}
  \put(115,117){\vector(-1,0){15}}
 \put(115,113){\vector(1,0){0}}
  \put(105,113){\makebox(0,0){$\cdots$}} \put(120,115){\makebox(0,0){$\cdot$}}
  \put(120,125){\makebox(0,0){${ \hat{f}_{0}(2)}$}} \put(122,113){\line(1,-1){10}}
  \put(137,98){\vector(1,-1){10}}\put(122,113){\line(1,-3){8}} \put(132,83){\line(1,-3){8}}
  \put(142,53){\vector(1,-3){8}}
  \put(65,115){\vector(-1,-1){15}}
 \put(0,100){\makebox(0,0){$\cdot$}}
 \put(20,102){\vector(-1,0){15}}
 \put(20,98){\vector(1,0){0}}
  \put(10,98){\makebox(0,0){$\cdots$}}
  \put(25,100){\makebox(0,0){$\cdot$}}
 \put(45,102){\vector(-1,0){15}}
 \put(45,98){\vector(1,0){0}}
  \put(35,98){\makebox(0,0){$\cdots$}}
   \put(50,100){\makebox(0,0){$\cdot$}}

   \put(70,55){\makebox(0,0){$\cdot$}}
  \put(70,55){\makebox(0,0){$\cdot$}}
 \put(90,57){\vector(-1,0){15}}
 \put(90,53){\vector(1,0){0}}
  \put(80,53){\makebox(0,0){$\cdots$}} \put(95,55){\makebox(0,0){$\cdot$}}
  \put(115,57){\vector(-1,0){15}}
 \put(115,53){\vector(1,0){0}}
  \put(105,53){\makebox(0,0){$\cdots$}} \put(120,55){\makebox(0,0){$\cdot$}}
  \put(120,65){\makebox(0,0){${ \hat{f}_{0}(1)}$}}\put(122,53){\line(1,-1){10}}
  \put(137,38){\vector(1,-1){10}}
  \put(65,55){\vector(-1,-1){15}}
 \put(0,40){\makebox(0,0){$\cdot$}}
 \put(20,42){\vector(-1,0){15}}
 \put(20,38){\vector(1,0){0}}
  \put(10,38){\makebox(0,0){$\cdots$}}
  \put(25,40){\makebox(0,0){$\cdot$}}
 \put(45,42){\vector(-1,0){15}}
 \put(45,38){\vector(1,0){0}}
  \put(35,38){\makebox(0,0){$\cdots$}}
   \put(50,40){\makebox(0,0){$\cdot$}}

    \put(100,85){\makebox(0,0){$\cdot$}}
  \put(100,85){\makebox(0,0){$\cdot$}}
 \put(105,87){\vector(1,0){15}}
 \put(105,83){\vector(-1,0){0}}
  \put(117,83){\makebox(0,0){$\cdots$}} \put(125,85){\makebox(0,0){$\cdot$}}
  \put(130,87){\vector(1,0){15}}
 \put(130,83){\vector(-1,0){0}}
  \put(142,83){\makebox(0,0){$\cdots$}} \put(150,85){\makebox(0,0){$\cdot$}}
 \put(82,70){\vector(1,1){15}}
 \put(30,70){\makebox(0,0){$\cdot$}}\put(30,80){\makebox(0,0){$f_{0}(2)$}}\put(28,72){\line(-1,1){10}}
 \put(15,85){\vector(-1,1){10}}\put(28,68){\line(-1,-1){10}}
 \put(15,55){\vector(-1,-1){10}}
 \put(35,72){\vector(1,0){15}}
 \put(35,68){\vector(-1,0){0}}
  \put(47,68){\makebox(0,0){$\cdots$}}
  \put(55,70){\makebox(0,0){$\cdot$}}
 \put(60,72){\vector(1,0){15}}
 \put(60,68){\vector(-1,0){0}}
  \put(72,68){\makebox(0,0){$\cdots$}}
   \put(80,70){\makebox(0,0){$\cdot$}}

         \put(100,25){\makebox(0,0){$\cdot$}}
  \put(100,25){\makebox(0,0){$\cdot$}}
 \put(105,27){\vector(1,0){15}}
 \put(105,23){\vector(-1,0){0}}
  \put(117,23){\makebox(0,0){$\cdots$}} \put(125,25){\makebox(0,0){$\cdot$}}
  \put(130,27){\vector(1,0){15}}
 \put(130,23){\vector(-1,0){0}}
  \put(142,23){\makebox(0,0){$\cdots$}} \put(150,25){\makebox(0,0){$\cdot$}}
 \put(82,10){\vector(1,1){15}}
 \put(30,10){\makebox(0,0){$\cdot$}}\put(30,20){\makebox(0,0){$f_{0}(1)$}}\put(28,12){\line(-1,1){10}}
 \put(15,25){\vector(-1,1){10}}
 \put(35,12){\vector(1,0){15}}
 \put(35,8){\vector(-1,0){0}}
  \put(47,8){\makebox(0,0){$\cdots$}}
  \put(55,10){\makebox(0,0){$\cdot$}}
 \put(60,12){\vector(1,0){15}}
 \put(60,8){\vector(-1,0){0}}
  \put(72,8){\makebox(0,0){$\cdots$}}
   \put(80,10){\makebox(0,0){$\cdot$}}
\end{picture}
\end{figure}

\emph{But they will not provide new modules, which is a result of
our next theorem.}\end{remark}

Combing the Auslander-Reite sequences constructed for
$W^{\lambda}(n),\tilde{W}^{\lambda}(n),\\T^{\lambda}(s,n)$ and
Proposition 5.8, we have $\mathbb{P}^{1}\kappa$ family of
homogeneous tubes now for any fixed $\lambda\in
\{0,1,\ldots,2p-1\}$. Comparing with the Auslander-Reiten quiver
constructed after Remark 5.7, we can summarize our works up to now
into the following result.

\begin{theorem} The modules

\emph{(1)} $P^{\lambda}$ for $0\leq \lambda\leq p-1$,

 \emph{(2)}
$V^{\lambda}(n),\tilde{V}^{\lambda}(n)$ for $0\leq \lambda\leq p-1$,
$n\geq 0$,

\emph{(3)} $W^{\lambda}(n),\tilde{W}^{\lambda}(n)$ for $0\leq
\lambda\leq p-1$, $n\geq 1$, and

 \emph{(4)} $T^{\lambda}_{c}(n)$ for
$c\in \kappa^{\ast}$, $0\leq \lambda\leq p-1$, $n\geq 1$,

\noindent form a complete list of all finite dimensional
indecomposable modules of $|\mathbf{u}(\mathfrak{osp}(1|2))|$ up to
isomorphism.
\end{theorem}

It is not hard to see that all indecomposable
$|\mathbf{u}(\mathfrak{osp}(1|2))|$-modules we constructed are
indeed supermodules naturally. For example, for the simple module
$V^{\lambda}$ defined as in (5.7), we can set
$v_{0},v_{2},\ldots,v_{2\lambda}$ to be even elements while
$v_{1},v_{3},\ldots,v_{2\lambda-1}$ be odd elements. From this,
$V^{\lambda}$ is an $\mathbf{u}(\mathfrak{osp}(1|2))$-supermodule.
Similar for other modules.

So now we can assume all modules in Theorem 5.10 are
$\mathbf{u}(\mathfrak{osp}(1|2))$-supermodules. Let $\Pi$ be the
\emph{parity change functor}, by definition it just interchanges the
$\mathbb{Z}_{2}$-grading of a supermodule.

\begin{corollary}The modules

\emph{(1)} $P^{\lambda}, \Pi(P^{\lambda})$ for $0\leq \lambda\leq
p-1$,

\emph{(2)}
$V^{\lambda}(n),\tilde{V}^{\lambda}(n),\Pi(V^{\lambda}(n)),\Pi(\tilde{V}^{\lambda}(n))$
for $0\leq \lambda\leq p-1$, $n\geq 0$,

 \emph{(3)}
$W^{\lambda}(n),\tilde{W}^{\lambda}(n),
\Pi(W^{\lambda}(n)),\Pi(\tilde{W}^{\lambda}(n))$ for $0\leq
\lambda\leq p-1$, $n\geq 1$, and

\emph{(4)} $T^{\lambda}_{c}(n),\Pi(T^{\lambda}_{c}(n))$ for $c\in
\kappa^{\ast}$, $0\leq \lambda\leq p-1$, $n\geq 1$,

\noindent form a complete list of all finite dimensional
indecomposable supermodules of $\mathbf{u}(\mathfrak{osp}(1|2))$ up
to isomorphism.
\end{corollary}
\begin{proof} It is enough to show every indecomposable supermodule
$M$ is indeed indecomposable as an
$|\mathbf{u}(\mathfrak{osp}(1|2))|$-module. Assume now
$M=M_{1}\oplus M_{2}$ in the
$|\mathbf{u}(\mathfrak{osp}(1|2))|$-module category. Let
$\pi:\bigoplus_{i\in I} P_{i}\to M$ be the projective cover of $M$
in the category of supermodules. Here we assume every $P_{i}$ is
indecomposable as a supermodule. By Proposition 12.2.12 in
\cite{Kle} and our description of projective
$|\mathbf{u}(\mathfrak{osp}(1|2))|$-modules, all $P_{i}$ are indeed
 indecomposable projective
 $|\mathbf{u}(\mathfrak{osp}(1|2))|$-modules. So $\bigoplus_{i\in I}
P_{i}\to M$ is also a surjection as
$|\mathbf{u}(\mathfrak{osp}(1|2))|$-modules. Therefore, we can
assume that there is subset $J\subset I$ such that $
\pi(\bigoplus_{i\in J}P_{i})=M_{1}$ and so $M_{1}$ is a supermodule.
Similarly, $M_{2}$ is a supermodule too. Thus $M=M_{1}$ or
$M=M_{2}$.
\end{proof}

\section*{Acknowledgements} The author is supported by
Japan Society for the Promotion of Science under the item ``JSPS
Postdoctoral Fellowship for Foreign Researchers" and Grant-in-Aid
for Foreign JSPS Fellow. I would gratefully acknowledge JSPS. The
work is also supported by Natural Science Foundation (No. 10801069).
I would like thank Professor A. Masuoka for stimulating discussions
and his encouragements.


\begin{thebibliography}{20}
\bibitem{Alp} J. Alperin, Periodicity in groups, Illinois J. Math.
21(1977) 776-783.

 \bibitem{ARS} M. Auslander, I. Reiten,  S. Smal$\phi$,
Representation theory of artin algebras, Cambridge University Press,
1995.

 \bibitem{BF} A. Bell, R. Farnsteiner, On the theory of Frobenius extensions and its application to Lie superalgebras,
  Trans. Amer. Math. Soc. 335(1993), no. 1, 407-424.

  \bibitem{BKN} B. Boe, J. Kujawa, D. Nakano, Cohomology and support
  varieties for Lie superalgebras, to appear in Trans. Amer. Math.
  Soc.  arXiv: math/0609363.

 \bibitem{B} J. Brundan, Modular representations of
 the supergroup $Q(n)$ Part II, Pacific J. Math. 224(2006), no. 1, 65-90.

 \bibitem{BK1} J. Brundan, A. Kleshchev, Projective representations
 of summetric groups via Sergeev duality, Math. Z. 239(2002), 27-68.

 \bibitem{BK2} J. Brundan, A. Kleshchev, Modular representations of
 the supergroup $Q(n)$, I, J. Algebra 260(2003), 64-98.

 \bibitem{BKu} J. Brundan, J. Kujawa, A new proof of the Mullineux conjecture,
  J. Algebraic Combin. 18(2003), no. 1, 13-39.

\bibitem{Ca1} J. Carlson, The varieties and cohomology ring of a
module, J. Algebra 85(1983), 104-143.

\bibitem{Ca2} J. Carlsom, The variety of an indecomposable module is
connected, Invent. Math. 77(1984), 291-299.

 \bibitem{DM} P. Deligne, J. Morgan, Note on supersymmetry (following
 Joseph Bernstein), Quantum fields and strings; a course for
 mathematicians, Vol. 1 (Princeton, NJ, 1996/1997), 41-96, Amer.
 Math. Soc., Providence, RI, 1999.

 \bibitem{Dr} Y. Drozd, Tame and wild matrix problems,
Representations and Quadratic Forms. Inst.Math.,
Acad.Sciences.Ukrainian SSR, Kiev 1979, 39-74. Amer.Math.Soc.
Transl. 128(1986), 31-55.

\bibitem{Far1} R. Farnsteiner,  Block representation type of Frobenius kernels
of smooth groups, J. Reine Angew. Math. 586(2005), 45-69.

 \bibitem{Far} R. Farnsteiner, Note on Frobenius extensions and restricted Lie superalgebras,
  J. Pure Appl. Algebra 108(1996), no. 3, 241-256.

\bibitem{FW} J. Feldvoss, S. Witherspoon, Support varieties and
representation type of small quantum groups, to appear in IMRN.
arXiv: 0910.1383.

\bibitem{FB} E. Friedlander, B. Parshall, Support
varieties for restricted Lie algebras, Invent. Math. 86(1986), no.
3, 553-562.

\bibitem{FS} E. Friedlander, A. Suslin, Cohomology of finite groups
schemes over a field, Invent. Math. 127(2)(1997), 209-270.

\bibitem{Kac1} V. Kac, Lie superalgebras, Adv. in Math. 26(1977),
8-96.

\bibitem{Kle} A. Kleshchev, Linear and projective representations of symmetric groups,
 Cambridge Tracts in Mathematics, 163. Cambridge University Press, Cambridge, 2005.

\bibitem{Ku} J. Kujawa, The Steinberg tensor product theorem for ${\rm GL}(m\vert n)$,
 Representations of algebraic groups, quantum groups, and Lie algebras, 123-132, Contemp. Math., 413,
 Amer. Math. Soc., Providence, RI, 2006.

 \bibitem{Majid} S. Majid, Crossed products by braided groups and bosonization,
J. Algebra 163(1994), 165-190.

\bibitem{Man} Yu. Manin, Gauge field theory and complex geometry,
Grundlehren der mathematischen Wissenschaften 289, Second Edition,
Springer, 1997.

\bibitem{MPSW} M. Mastnak, J. Pevtsova, P. Schauenburg, S.
Witherspoon, Cohomology of finite dimensional pointed Hopf algebras,
Proc. Lond. Math. Soc. (3) 100(2010), no. 2, 377-404.

\bibitem{MW} M. Mastnak, S. Witherspoon, Bialgebra cohomology, pointed Hopf algebras, and deformations,
 J. Pure Appl. Algebra 213(2009), no. 7, 1399-1417.

\bibitem{Mon} S. Montgomery, Hopf Algebras and Their Actions
on Rings. CBMS, Lecture in Math.; Providence, RI, (1993); Vol. 82.

\bibitem{NP} D. Nakano, J. Palmieri, Support varieties for the
Steenrod algebra, Math. Z. 227(1998), 663-684.

\bibitem{Radford} D.Radford, The structure of Hopf algebras with
a projection, J. Algebra 92(1985), 322-347.

\bibitem{Ri} C. Ringel, Tame algebras and integral quadratic forms,
Springer, Lecture Notes in Math. 1099, 1984.

\bibitem{SW} B. Shu, W. Wang, Modular representations of the ortho-symplectic supergroups,
 Proc. Lond. Math. Soc. (3) 96(2008), no. 1, 251-271.

\bibitem{SS} N. Snashall, O. Solberg, Support varieties and Hochschild cohomology rings, Proc. London
Math. Soc. (3) 88(2004), no. 3, 705-732.

\bibitem{So} O. Solberg, Support varieties for modules and
complexes, Trends in representation theory of algebras and related
topics, 239-270, Contemp. Math. 406, Amer. Math. Soc. Providence,
RI, 2006.

 \bibitem{SFB} A. Suslin, E. Friedlander, C. Bendel,
  Support varieties for infinitesimal group schemes, J. Amer. Math. Soc. 10(1997), no. 3, 729-759.

 \bibitem{WZ} W. Wang, L. Zhao, Representations of Lie superalgebras in prime characteristic. I,
  Proc. Lond. Math. Soc. (3) 99(2009), no. 1, 145-167.

  \bibitem{WZ2} W. Wang, L. Zhao, Representations of Lie superalgebras in prime characteristic II: The queer
  series. arXiv:0902.2758

  \bibitem{We} C. Weibel, An introduction to homological algebra,
  Cambridge University Press, 1994.

\bibitem{Xiao} J. Xiao, Finite-dimensional representations of
$U_{t}(\mathfrak{sl}_{2})$ at roots of unity, Can. J. Math.
49(1997), no. 4, 772-787.

\bibitem{Z} C. Zhang, On the simple modules for the restricted Lie superalgebra ${\rm sl}(n\vert 1)$,
 J. Pure Appl. Algebra 213(2009), no. 5, 756-765.
\end{thebibliography}
\end{document}